%% file: main.tex
\documentclass{article}
\usepackage[utf8]{inputenc}

\input{defs}

\usepackage{enumitem}

\usepackage[
  margin=2cm,
  includefoot,
  footskip=30pt,
]{geometry}

\usepackage{thmtools}
\usepackage{thm-restate}
\usepackage{thm-autoref}

\usepackage{pdflscape}
\usepackage{afterpage}

\usepackage{multirow}

\usepackage{booktabs}
\usepackage{threeparttable}

\usepackage{tabularx}
\usepackage{siunitx}
\usepackage{dcolumn}

\newcommand{\dan}[1]{{\bf \color{green} Dan: #1}}
\newcommand{\rasmus}[1]{{\bf \color{olive} Rasmus: #1}}
\newcommand{\yuan}[1]{{\bf \color{olive} Yuan: #1}}
\renewcommand{\dan}[1]{}
\renewcommand{\rasmus}[1]{}
\renewcommand{\yuan}[1]{}

\newcommand{\todo}[1]{{\bf \color{red} TODO: #1}}
\newcommand{\todomed}[1]{{\bf \color{orange} TODOMED: #1}}
\newcommand{\todolow}[1]{{\bf \color{Dandelion} TODOLOW: #1}}
\renewcommand{\todo}[1]{}
\renewcommand{\todomed}[1]{}
\renewcommand{\todolow}[1]{}

\newcommand{\AC}{AC}
\newcommand{\ACSM}{AC2} %

\newcolumntype{s}{D{.}{.}{1.1}}
\allowdisplaybreaks

\title{\vspace{-3em}Robust and Practical Solution of Laplacian Equations\\
by Approximate Elimination%
\thanks{The research leading to these results has received funding from the grant ``Algorithms and complexity for high-accuracy flows and convex optimization'' (no. 200021 204787) of the Swiss National Science Foundation. It was also supported in part by NSF Grant CCF-1562041, ONR Award N00014-16-2374, and a Simons Investigator Award to Daniel Spielman}
}
\author{Gao, Kyng, Spielman}

\author{
Yuan Gao\thanks{Part of this work was conducted as a master's thesis at ETH Zurich. Part of this work was also conducted while the author is at MPI for Informatics} %
 \\ \small CISPA
 \\ \small yuan.gao@cispa.de
\and
Rasmus Kyng\\  \small ETH Zurich \\  \small kyng@inf.ethz.ch 
\and
Daniel A. Spielman  \\  \small Yale University \\  \small spielman@cs.yale.edu
\and
}

\date{May 2022}

\begin{document}

\maketitle
\vspace{-2em}
\begin{abstract} 
 We introduce a new algorithm and software for solving linear equations
 in symmetric diagonally dominant matrices with non-positive
 off-diagonal entries (SDDM matrices), including Laplacian matrices.
 We use preconditioned conjugate gradient to solve the system of
 linear equations.
 Our preconditioner is a variant of the Approximate Cholesky factorization
 of Kyng and Sachdeva (FOCS 2016).
 Our factorization approach is simple: we eliminate matrix
 rows/columns one at a time, and update the remaining entries of the
 matrix by sampling entries to approximate the outcome of complete
 Cholesky factorization.
 Unlike earlier approaches, our sampled entries always
 maintain a connected graph on the neighbors of the eliminated variable.
 Our algorithm comes with a tuning parameter that upper bounds the
 number of samples made per original entry.

 We implement our solver algorithm in Julia, and experimentally evaluate 
 its 
 performance when using 1 or 2 samples for each original
 entry: We refer to these variants as AC and AC2 respectively.
 We investigate the performance of
 these implementations and compare their single-threaded performance
 to that of current state-of-the-art solvers including Combinatorial
 Multigrid (CMG), BoomerAMG-preconditioned Krylov solvers from HyPre and PETSc, Lean
 Algebraic Multigrid (LAMG), and MATLAB's preconditioned conjugate gradient
 with Incomplete Cholesky Factorization (ICC). 
  Our experiments suggest that {\ACSM} and {\AC} attain a
level of robustness and reliability not seen before in solvers for
SDDM linear equations, while retaining good performance across all instances.

Many large-scale evaluations of SDDM linear equation solvers have focused on solving problems on discretized 3D grids.
We conduct evaluations across a much broader class of problems, including all large SDDM matrices from the SuiteSparse collection, as well as a broad array of programmatically generated instances.
Our tests range up to 200 million non-zeros per system of linear equations.
Our experiments show that AC and AC2 obtain good practical
 performance across many different types of SDDM-matrices, 
 and significantly greater reliability than existing solvers.
 AC2 is the only solver that succeeds across all our
 tests, and it uses less than $7.2 \mu s$ time per non-zero across to converge to $10^{-8}$ relative residual error across all our experiments.
 AC is typically 1.5-2 times faster, but fails on one family of problems
 engineered to attack this algorithm.
 In our experiments using general sparse non-zero patterns,
CMG, HyPre, PETSc, and ICC all fail on some instances from a majority of the families tested.
 Across these families, the median running time across different instances of AC2 and AC is comparable to or faster than other solvers.

 We also test the performance of our solvers on a wide array of 
 Poisson problems on 3D grids, including grids with uniform,
 high-contrast, or anisotropic coefficients, and grids from the SPE Benchmark.
 Here, the CMG and HyPre solvers perform best, but
 AC and AC2 achieve worst case and median 
 running times within a factor 4.1 and 6.2 of these respectively.
 
 Our code is public, and we detail precisely the set of tests we run
 and provide a tutorial on how to replicate the tests.
 We hope that others will adopt this suite of tests as a benchmark,
 which we refer to as SDDM2023. 
 Our solver code is available at:
 
 \begin{mdframed}[backgroundcolor=lightgray]
   \centering
   \url{https://github.com/danspielman/Laplacians.jl/}
 \end{mdframed}
 
 Our benchmarking data and tutorial is available at:

 \begin{mdframed}[backgroundcolor=lightgray]
   \centering
       \url{https://rjkyng.github.io/SDDM2023/}
 \end{mdframed}
\end{abstract}

\newpage
\input{intro}

\newpage
\tableofcontents
\newpage

\input{prelims}
\input{algo}

\input{results}

\FloatBarrier

\bibliographystyle{alpha}
\bibliography{refs}

\appendix

\input{appendixTheory}

\input{appendixUnbiased}

\input{appendixEquivalence}

\end{document}

%% file: defs.tex
\usepackage{blindtext}
\usepackage[usenames,dvipsnames,svgnames,table]{xcolor}
\usepackage[inline]{enumitem}
\usepackage{caption}
\usepackage{placeins}
\usepackage{graphicx}

\usepackage{fancybox}
\usepackage{hyperref}

\usepackage{amsmath,amsfonts,amsthm,amssymb,xcolor}

\usepackage{mdframed}

\usepackage{bm} %

\usepackage{nicefrac}

\usepackage{algorithm} %
\usepackage{algorithmicx}
\usepackage[noend]{algpseudocode}

\newtheorem{theorem}{Theorem}[section]
\newtheorem{corollary}[theorem]{Corollary}

\newtheorem{claim}[theorem]{Claim}
\newtheorem{fact}[theorem]{Fact}

\newtheorem*{theorem*}{Theorem}
\newtheorem*{claim*}{Claim}
\newtheorem*{corollary*}{Corollary}
\newtheorem*{conjecture*}{Conjecture}
\newtheorem*{lemma*}{Lemma}
\newtheorem*{thm*}{Theorem}
\newtheorem*{prop*}{Proposition}
\newtheorem*{obs*}{Observation}
\newtheorem*{definition*}{Definition}

\newtheorem*{remark*}{Remark}
\newtheorem*{rec*}{Recommendation}

\theoremstyle{definition}

\newtheorem{remark}[theorem]{Remark}

\newenvironment{fminipage}%
  {\begin{Sbox}\begin{minipage}}%
  {\end{minipage}\end{Sbox}\fbox{\TheSbox}}

\def\defeq{\stackrel{\mathrm{def}}{=}}
\def\setof#1{\left\{#1  \right\}}
\def\sizeof#1{\left|#1  \right|}

\def\floor#1{\left\lfloor #1 \right\rfloor}

\def\abs#1{\left|#1  \right|}

\def\norm#1{\left\| #1 \right\|}

\newcommand{\trp}{\top}

\newcommand\grad{\boldsymbol{\nabla}}

\newcommand\PPhi{\boldsymbol{\Phi}}

\newcommand\PPi{\boldsymbol{\Pi}}

\def\aa{\pmb{\mathit{a}}}
\newcommand\bb{\boldsymbol{\mathit{b}}}

\newcommand\ee{\boldsymbol{\mathit{e}}}

\renewcommand\gg{\boldsymbol{\mathit{g}}}

\renewcommand\ll{\boldsymbol{\mathit{l}}}

\newcommand\yy{\boldsymbol{\mathit{y}}}

\newcommand\xx{\boldsymbol{\mathit{x}}}

\newcommand\veczero{\boldsymbol{0}}
\newcommand\vecone{\boldsymbol{1}}
\newcommand\vecind{\boldsymbol{\chi}}

\renewcommand\AA{\boldsymbol{\mathit{A}}}
\newcommand\BB{\boldsymbol{\mathit{B}}}
\newcommand\CC{\boldsymbol{\mathit{C}}}

\newcommand\II{\boldsymbol{\mathit{I}}}

\newcommand\MM{\boldsymbol{\mathit{M}}}
\newcommand\LL{\boldsymbol{\mathit{L}}}
\newcommand\PP{\boldsymbol{\mathit{P}}}

\renewcommand\SS{\boldsymbol{\mathit{S}}}

\newcommand\CCtil{\boldsymbol{\mathit{\tilde{C}}}}

\newcommand\SStil{\boldsymbol{\mathit{\tilde{S}}}}

\newcommand\matzero{\boldsymbol{0}}

\newcommand\R{\mathbb{R}}

\newcommand{\expct}[2]{\mathop{{}\mathbb{E}}_{#1}\left[#2\right]}

\DeclareMathOperator*{\diag}{diag}

\newcommand{\matlow}{\boldsymbol{\mathit{{\mathcal{L}}}}}

\newcommand{\matdiag}{\boldsymbol{\mathit{{\mathcal{D}}}}}

\newcommand{\matup}{\boldsymbol{\mathit{{\mathcal{U}}}}}

\newcommand{\cliq}[2]{\textsc{Clique}\!\left[{#1}\right]_{#2}}

\newcommand{\vstar}[2]{\textsc{Star}\!\left[#1\right]_{#2}}

\newcommand{\estar}[3]{\textsc{ElimStar}\!\left[#1\right]_{#2,#3}}

\newcommand{\ulap}{\LL}

\newcommand{\schurto}[2]{\ensuremath{\textsc{Sc}\!\left[#1\right]_{#2}}}

\newcommand{\pinv}{\dagger}

\newcommand{\proj}{\PPi}

\DeclareMathOperator{\nnz}{nnz}

%% file: intro.tex
\section{Introduction}

Linear equations in symmetric diagonally dominant matrices with non-positive off-diagonal matrices 
  (SDDM) matrices appear in many
  applications,
  including solvers for discretized scalar elliptic
  partial differential equations \cite{vdVM77,K78,CB01,BHV08,SBBRGS12}, many
  problems in computer vision and computer graphics \cite{KMT11}, and in
  machine learning \cite{ZGL03}.
They are essentially equivalent to the Laplacian matrices of graphs.
  
In this paper, we introduce a new algorithm for solving SDDM linear
equations, as well as an single-threaded implementation of the
algorithm and several variations in Julia. 
We assemble a large suite of SDDM matrices on which to to test solvers,
  and use these to evaluate the performance of our implementation 
  and to compare it with other popular solvers for SDDM
linear equations.
We include important test cases from finite-difference method
approaches to solving scalar ellipic partial differential equations, all SDDM
matrices from the \emph{SuiteSparse} collection, other classic
benchmark data, and several programmatically generated difficult
instances, including problems coming from running Interior Point
Methods for linear programming. We hope this collection of matrices will be adopted as a
benchmark in future work. 
One can use a reduction of  
  Gremban~\cite{G96} to convert the problem of solving an $n \times n$ SDD
linear equation into one of solving a $2n+1 \times 2n+1$ Laplacian linear equation,
  but in this work we do not test the performance of our
algorithm on instances coming from general SDD matrices.

Current popular implementations of SDDM solvers \cite{HY02, BGHSTMAT19, LB12} rely on variants of
multigrid methods \cite{B77,B02}.
There has been
tremendous theoretical success in developing SDDM linear equation solvers with
running times asymptotically nearly-linear in the number of non-zeros
of the matrix starting with the work of Spielman and Teng~ \cite{ST04,KMP10,KMP11,KOSZ13,PS14,KLPSS16,KS16,JS21,HJPW23},
but so far, these algorithms have not lead directly to 
practical solver implementations.
This line of work on asymptotically fast solvers builds on support theory,
an approach to preconditioning Laplacian linear equations that expands on a breakthrough of Vaidya \cite{V90,G96,BH03,BCHT04,BGHNT06}.
Ideas from support theory have also had significant impact on the Combinatorial Multigrid Solver by Koutis, Miller, and Tolliver \cite{KMT11}, which combined these ideas with a multi-grid approach.

\subsection{Related Work}

\paragraph{Code and experiments.}
Incomplete Cholesky factorizations were introduced by 
  Meijerink and van der Vorst~\cite{MvdV77} to accelerate the solution 
  of systems of linear equations in symmetric M-matrices.
These are obtained by omitting most entries that appear during Cholesky factorization,
  and the factorization they produce can be shown to crudely approximate the original matrix.
 Meijerink and van der Vorst suggested using these factorizations as a preconditioners for the Conjugate Gradient.
When they were introduced, they provided the fastest solutions to the discrete approximations of two-dimensional elliptic PDEs.
A popular implementation
  of this method is MATLAB's \texttt{ichol}.
Two key differences between these factorizations and our new Approximate Choleksy factorizations are that the Approximate Cholesky factorizations increase the magnitudes of the entries they keep and that they select entries to keep at random.

The Multigrid algorithm~\cite{BD79} provided a significant improvement in the solution
  of systems of equations arising from discrete approximations of elliptic PDEs,
  and its generalizations provide the fastest solvers for many families of linear equations.
Algebraic Multigrid (AMG) builds upon the principles of Multigrid
methods but has a wider range of applications as it does not require
the problem to have a simple underlying geometry, making it applicable
to general symmetric M-matrices~\cite{RS87}, a family that includes SDDM matrices.
Brandt designed an AMG solver applicable to Dirac
equations~\cite{B00}. For linear equations in graph Laplacians, Livne and Brandt presented a variant
of AMG solver, called the Lean Algebraic Multigrid (LAMG), and
provided a MATLAB implementation~\cite{LB12}.
Another multi-grid variant specialized for graph Laplacians was developed by Napov and Notay \cite{NN17}.

BoomerAMG, a parallel AMG implementation, was introduced in \cite{HY02}, and can be accessed via the popular software library
HyPre~\cite{FY02}.
The Portable, Extensible Toolkit for Scientific Computation (PETSc)
provides another interface for accessing the HyPre-BoomerAMG implementation and its own implementations of conjugate 
gradients~\cite{BAABBBDDEG19}.
Combinatorial Multigrid (i.e. CMG) \cite{KMT11} is a solver that combines ideas 
  from support theory \cite{V90,BH03}  with ideas from (Algebraic) Multigrid
methods.
The MATLAB implementation of CMG was applied to problems arising in computer vision and
image processing in~\cite{KMT11}.
\cite{LHLLE19} presented another implementation of CMG in C using PETSc, however, to the best
of our knowledge, this implementation is not publicly available.
As part of the Trilinos software package, MueLu \cite{BGHSTMAT19}
provides a multigrid algorithm designed for solving sparse linear systems of
equations arising from PDE discretizations using massively parallel
computing environments.

In \cite{CT03}, Chen and Toledo evaluated an SDD solver based on Vaidya's ideas on tree-based preconditioning combined with PCG. 
They compared the solver with PCG preconditioned using Incomplete Cholesky factorization and that Vaidya's approach was often much slower, but sometimes lead to convergence in cases where Incomplete Cholesky with PCG did not.
In \cite{DGMPXX16}, Deweese et al.\ studied the performance of variants of the cycle toggling SDD solver from \cite{KOSZ13}.
This algorithm works in with a `dual space' representation (sometimes known as the `flow space').
The paper focused on contrasting these implementations without a direct comparison to other solvers.
In \cite{BDG16}, Boman, Deweese, and Gilbert evaluated a cycle toggling algorithm in `primal space' (`voltage space') known as Primal Randomized Kaczmarz (PRK). They compared this with an implementation of the \cite{KOSZ13} solver and with a PCG implementation using Jacobi diagonal scaling.
The authors find that their results ``do not at present support the practical utility of'' cycle toggling solvers in primal or dual space.
Hoske et al.\ reached similar conclusions in another experimental evaluation of an implementation of the \cite{KOSZ13} cycle toggling solver \cite{HLMW16}.

A recent paper by Chen, Liang, and Biros has provided
another implementation of our sampling procedure ~\cite{CLB21} though
only the single sample-version (corresponding to our {\AC} implementation).
Unlike our implementation, the authors used the standard Cholesky
lower-triangular format for the output, whereas our implementation uses
a row-operation representation of the output.
Chen et al. use a METIS-based elimination ordering while ours is adaptive.
Chen et al. have developed a parallel version of the algorithm and
code using nested dissection ordering.
Their implementation, known as RCHOL, is available at
\url{https://github.com/ut-padas/rchol} via C++/MATLAB/Python
interfaces.
Chen, Liang, and Biros ~\cite{CLB21} experimental evaluation focused on comparing RCHOL with three other solvers: (1) the MATLAB'S preconditioned conjugate gradient combined with the ichol implementation of Incomplete Cholesky factorization, (2) Ruge--Stuben AMG (RS-AMG from pyamg), (3) the smoothed aggregation AMG (SA-AMG also from pyamg), and (4) Combinatorial Multigrid \cite{KMT11}.
Their experiments were based on three matrices arising from 3D grid Poisson problems (with uniform, variable, and anisotropic coefficients) and four matrices from the SuiteSparse matrix collection: one SDD, one SDDM and two which are neither.
They found that RCHOL outforms ichol, and found each of the four solvers RCHOL, CMG, RS-AMG, and SA-AMG performs best on some of the tested problems.

\paragraph{Theory.}
There are some theories of convergence for (Algebraic and Geometric)
Multigrid methods \cite{BH83, Y93}  and preconditioned iterative methods using
Incomplete Cholesky factorization \cite{MvdV77,gustafsson,beauwens1994approximate,notay1990solving,notay1992conditioning,SupportGraph}, but these do not establish fast rates for general SDDM matrices.

Spielman and Teng \cite{ST04} built on the work of Vaidya \cite{V90} to develop
the first provably correct
nearly-linear time algorithm for solving SDD linear equations
\cite{ST04}.
Koutis, Miller and Peng
  substantially simplified the algorithm of Spielman and
  Teng and greatly reduced its running time~\cite{KMP10,KMP11}.
Building on this line of work, \cite{CKPPR14} reduced the running time
below that of comparison sorting the input, and recently \cite{JS21}
reduced the running time to linear in the number of matrix non-zeros,
up to poly-loglog factors and a $\log(1/\epsilon)$ dependence on the error, $\epsilon$.
The solver of Kelner et al. introduced a dual-solution based approach,
resulting a simple algorithm, once a low-stretch spanning tree for the matrix has been computed.
Another line of work introduced a parallel algorithm~\cite{PS14} and
removed the reliance on combinatorial sparsification~\cite{KLPSS16}.
Finally, \cite{KS16} introduced and analyzed a very simple algorithm that is the inspiration for the algorithm in the present paper.

\section{Results and Conclusion}
We introduce a variant of the approximate Cholesky factorization
algorithm developed in \cite{KS16} that greatly improves its practical
performance while retaining its simplicitly.
Our algorithm solves linear equations in
SDDM matrices, including Laplacian matrices.
We conduct extensive experiments, testing the performance of two
variants of our algorithm across a wide range of SDDM matrices and
comparing our solver with a suite of state-of-the-art solvers.
Our experiments measure the single-threaded performance of these
solvers.
The experiments cover a broad range of equations in SDDM matrices,
including Laplacian matrices.

\subsection{Algorithms and Implementation}
\label{sec: algandimpl}
Our algorithm solves an SDDM linear
equation by reducing it to a Laplacian linear equation,
  computing an approximate Cholesky factorization of that Laplacian matrix,
  and using it as a preconditioner for the Laplacian inside the Conjugate Gradient.
The basic version of our algorithm eliminates rows (and the
corresponding column) of a Laplacian matrix one at a time.
For each off-diagonal entry in the eliminated row, a new entry is
sampled and inserted into the remaining matrix, while guaranteeing two
important properties: first, the expected value of the output
samples agree exactly with the output of Cholesky factorization, and second,
the sampled non-zeros preserve the connectivity structure
of the graph. Concretely, if the off-diagonal entries are viewed as edges
in a graph with vertices corresponding to the rows/columns of the
matrix, then the algorithm can be seen as eliminating the edges
incident to a vertex and replacing them with a tree on the neighbors
of the vertex.
In our implementations, we
 eliminate vertices with whose degree is approximately minimum at the
 time they are eliminated.
As the edges introduced during elimination are random, the resulting 
  elimination order appears unrelated to the approximate minimum degree
ordering~\cite{AMD}.
In our approximate factorization approaches, the order of elimination
of edges also changes the algorithm, and our implementation eliminates
edges in the order of lowest to heighest weight.

We refer to our implementation of the basic version of our algorithm as
Approximate Cholesky ({\AC}).
Experimentally, we find that {\AC} performs well on a
large range of SDDM matrices. 
However it performs poorly on matrices that were specially designed
 by Sushant Sachdeva to make it fail.
We call these \emph{Sachdeva stars}.
To address the problem caused by these examples, we develop a more robust, but slightly slower,
version of the algorithm: this version divides each matrix
entry into $k$ separate ``multi-entries'' (analogous to
multi-edges in a graph). After this, we run the original algorithm
on this structure -- but now the output of each sampled elimination
has the non-zero structure of a union of trees with at most as many multi-entries as
the eliminated row/column.
At $k = 2$, i.e. with just two multi-entries per original entry, this
robust version of the algorithm performs well across all our
tests.
We refer to this version of our algorithm as {\ACSM}.
While our algorithm has not been rigorously analyzed, the design
of the sampling procedure is motivated in part by
\cite{KS16}, as well as theoretical work which suggests graphs
may be well-approximated by a tree plus a few additional edges:
tree-based ultra-sparse sampling can yield convergent iterative
  methods for solving Laplacians~\cite{ST04}, even when the sampling is too sparse to yield
  concentration of measure w.r.t. spectral norms \cite{CKPPR14}; and in bounded degree graphs, a
  union of two random spanning trees forms a good cut
  sparsifier~\cite{GRV09}, while a few random spanning trees form a
  good spectral sparsifier~\cite{KS18,KKS22}.

The output of our algorithm can be represented in two different ways:
  either as a standard lower-triangular Cholesky factorization or as a 
  product of row operation matrices.
Regarded as linear operators, the two are identical
(i.e. when using the same source of randomness they will output
identical linear operators), but it is unclear which performs better
in practice, and which version has better numerical stability.
Our implementation uses the product of row operations form.

Our solvers are implemented in Julia and are available at:
.
 \begin{mdframed}[backgroundcolor=lightgray]
   \centering
   \url{https://github.com/danspielman/Laplacians.jl/}
 \end{mdframed}

\subsection{Experimental Evaluation and Comparison}
\label{sec:introexperiments}

\paragraph{Solvers.} We run a broad suite of experiments,
measuring the single-threaded performance of our {\AC} and {\ACSM}
solvers and comparing them with the following suite of popular SDDM linear
equation solvers:
\begin{itemize}
\item Combinatorial Multigrid (CMG),
\item HyPre's Krylov Solver with
  BoomerAMG preconditioning (HyPre),
\item PETSc's Krylov Solver with
  BoomerAMG preconditioning (PETSc),
\item MATLAB preconditioned conjugate gradient using the
  ichol implementation of Incomplete Cholesky factorization (ICC),
\item Lean Algebraic Multigrid (LAMG)\footnote{We omit LAMG from our result tables, as we
    found it unable to convergence to our target tolerance across all
    matrix families we tested.}.
\end{itemize}
Of the third-party solvers we test, i.e. CMG, HyPre, PETSc, ICC, and LAMG,
the latter three  seem to be dominated in performance by CMG or
HyPre, with few exceptions.
Consequently, our discussion will focus on comparing {\AC} and {\ACSM} with CMG and HyPre.

\paragraph{Summary of experimental results.}
We test three main categories of SDDM matrices.
\begin{enumerate}
    \item All large SDDM matrices in SuiteSparse.
    \item Programmatically generated matrices:
    \begin{enumerate}
        \item ``Chimeras'': a broad class of matrices that we generate programmatically to test solvers on a wide array of graph geometries.
        \item ``Sachdeva stars'': A special family of matrices designed to challenge our algorithm\footnote{We thank Sushant Sachdeva for suggesting this graph construction.}.
        \item ``Flow Interior Point Method Matrices'': We create variants of our ``Chimeras'' with different entries by running an Interior Point Method for a maximum flow problem on these graphs. This creates a sequence of linear equations in SDDM matrices with the same non-zero structure as the input graph, but with a challenging weight pattern.
        We also create sequences of SDDM matrices by running a maximum flow interior point method on \emph{Spielman graphs}, a family of graphs that are thought to be particularly challenging for interior point methods.
    \end{enumerate}
    \item Discretizations of partial differential equations on 3D grids:
    \begin{enumerate}
        \item Matrices based on fluid simulations from a Society of Petroleum Engineering (SPE) benchmark.
        \item Poisson problems on 3D grids with uniform, high-contrast coefficients, or anisotropic coefficients.
    \end{enumerate}
\end{enumerate}
The matrices we use have up to 200 million non-zeros.

{\ACSM} is able to convergence in a reasonable time on every
single test, while {\AC} is typically faster but fails to converge
on large \emph{Sachdeva stars}.
On general graphs, our solvers achieve a worst-case running time of about 7.2 $\mu s$ per non-zero entry, while on grid graphs we get a worst-case running time of about 3.6 $\mu s$ per non-zero.
This should be contrasted with all other solvers we test, which fail entirely fail to converge both on some matrices from SuiteSparse and some Chimeras, but achieve a worst case running time on grid graphs of about 0.67 $\mu s$ per non-zero (CMG) and $0.77$ $\mu s$ per non-zero (HyPre).
Our experiments show that AC and AC2 obtain good practical
performance across many different types of SDDM-matrices,
and significantly greater reliability than existing solvers.

Our test results are summarized in Table~\ref{tab:introsummarytable}.
\begin{table}[H]
   \centering
   \input{summarytable-times}
   \caption{
    Summary of worst case total solve time per non-zero matrix
    entry across all experiments.}
  \label{tab:introsummarytable}
\end{table}

\paragraph{Experiments on SuiteSparse SDDM matrices.}
We test the solvers on every SDDM matrix currently
hosted on the SuiteSparse Matrix Collection\footnote{https://sparse.tamu.edu}, 
previously known as the University of Florida Sparse Matrix Collection
\cite{DH11}, after excluding very small matrices (with less than 1000 non-zero entries) to limit the impact of start-up overhead.
{\ACSM} and {\AC} converge on all these problems in time at most 1.4 $\mu s$ per non-zero. Meanwhile HyPre and CMG fail on some instances, and PETSc and ICC respectively take 27.1  $\mu s$ and  153  $\mu s$ per non-zero in the worst case.
{\ACSM} and {\AC} are also faster than CMG, PETSc, and ICC in median times, but slower than HyPre on median times.
See Table~\ref{tab:suitesparse} for these results.

\paragraph{Experiments on programmatically generated SDDM matrices.}
We created a generator for weighted Laplacian matrices that we call \emph{Laplacian Chimeras}
  because they combine matrices of many types.
The results of these experiments on these are shown in Tables~\ref{tab:unilapchimera} and~\ref{tab:wtedlapchimera}.
We test the suite of solvers on these, as
well as a set of strictly diagonally dominant matrices derived from
the Laplacian Chimeras. Results are shown in Tables~\ref{tab:unisddchimera} and \ref{tab:wtedsddmlapchimera}.
Across the Laplacian and SDDM Chimeras we test, both CMG and HyPre fail
on some instances. On these Chimeras, the median running time of {\AC} and {\ACSM} is
lower than for CMG and HyPre.
Thus both the worst case performance and the typical performance of
{\AC} and {\ACSM} beats that of CMG and HyPre on this class of
problems.

We test our solvers on an additional generated family of graphs, which
we refer to as \emph{Sachdeva stars}, designed specifically to challenge our
family of algorithms.
On Sachdeva stars, {\ACSM} converges and {\AC} fails to converge. 
Results for Sachdeva stars are shown in Table~\ref{tab:startable}.

In addition to our tests on unweighted and randomly weighted Chimeras,
we run experiments based on using our Laplacian matrix linear equation
solvers inside an Interior Point Method (IPM) that solves an Undirected Maximum Flow
Problem to high accuracy.
Laplacian linear equations arise in IPM algorithms for
single-commodity flow problems because the algorithms proceed by
taking Newton steps that minimize a sequence of barrier function
problems, and each Newton step corresponds to a Laplacian linear
equation.
Laplacian linear equations arising from IPMs tend to have very large
ranges of weights, which leads to ill-conditioned problems, and
currently, to the best of our knowledge, no Laplacian solver has been
reliable enough to use in IPMs.
First, we consider Maximum Flow Problems on unweigted Chimera graphs.
We use a short-step IPM and solve the problems to high accuracy, see
Table~\ref{tab:chimeraIPM} for results.
Second, we consider a special family of graphs, known as \emph{Spielman
graphs}, which are conjectured to be a particularly hard instance for
short-step IPMs.
The results of this experiment are shown in Table~\ref{tab:ipmrloweps}.
Our IPM is conservative and uses far too many Newton steps to be
practical, but the experiments demonstrate that, at least for
short-step IPMs, both {\AC} and {\ACSM} work reliably across many
different instances.
We believe this suggests that {\AC} and {\ACSM} may be good candidates
for further research into developing Interior Point Methods for
single-commodity flow problems using fast linear equation solvers.
In contrast, all the other tested solvers fail on some IPM instances
-- which is unsurprising as the solvers even fail on unweighted
Chimera instances, and the IPM should create similar instances but
with more extreme weights.

\paragraph{Experiments on 3D grid discretizations.}
We find that the HyPre and CMG solvers are the fastest on grid problems,
up to 4.1 times faster than {\AC} and 6.3 times faster than {\ACSM}, see Tables
\ref{tab:uniformgridtable}, \ref{tab:checkeredboard}, 
\ref{tab:anisogridtable}, \ref{tab:wgridtable}, and \ref{tab:spetable}.
Our experiments on grids include both SDDM matrices generated from fluid flow problems using data from the Tenth SPE Comparative Solution Project \cite{CB01}, from the Society of Petroleum Engineering (SPE),
as well as Poisson 3D grids with uniform coefficients, high contrast coefficients, or anisotropic coefficients. 
These experiments suggest that our current implementation of {\AC} and {\ACSM} is slower than state-of-the-art solvers on 3D grid problems -- but the gap might be addressable by further optimization of our code.

\paragraph{Performance variability experiments.}
We run tests to demonstrate the low
variability in the performance of {\ACSM} and {\AC}, showing that their running times are reliable and consistent across many runs, see Table~\ref{tab:variance}.
We also run some experiments to shed light on why the more reliable {\ACSM} algorithm outperforms {\AC} in some settings.
In particular, we show that {\ACSM} builds a preconditioner with much smaller relative condition
number on difficult instances, see Table~\ref{tab:variants}.

\paragraph{Running time scaling with problem size.}
In Section~\ref{sec:scaling}, we show plots that display the running time of our {\AC} and {\ACSM} solvers plotted against non-zero count, demonstrating that the observed scaling of running time in problem size is nearly-linear.
Our plots for AC appear in Figures~\ref{fig:ac-nnz-vs-time} and \ref{fig:ac-nnz-vs-timebynnz} and for AC2 in Figures~\ref{fig:ac2-nnz-vs-time} and \ref{fig:ac2-nnz-vs-timebynnz}.

\subsection{Conclusion}
Our experiments suggest that {\ACSM} and {\AC} attain a
level of robustness and reliability not seen before in solvers for
SDDM linear equations, while retaining good performance across all instances.
For most problems,  {\AC} appears to be a good choice, while some very
challenging instances require using {\ACSM}.

We hope that our benchmarking data will be used in future works to provide a rigorous basis for comparison of new SDD solvers. 
This benchmarking data and an accompanying tutorial are available at:
 \begin{mdframed}[backgroundcolor=lightgray]
   \centering
    \url{https://rjkyng.github.io/SDDM2023/}
 \end{mdframed}

\paragraph{Future work.}
We plan to release a version of our code designed to handle symmetric diagonally dominant (SDD) matrices,
  which can have postitive off-diagonal entries, and symmetric block-diagonally
dominant matrices (SBDD).

\paragraph{Open problems.}
We hope that our results will motivate others to study ways to
incorporate insights from our solvers into production level codes for
solving SDDM systems.
Developing a highly-optimized version of our algorithm, even for
single-threaded systems, will be an important step in this discussion rection.
Creating an efficient parallel implementation is even more important,
and presents interesting questions about which parallelization
approach is likely to be most successful.
We believe a more in-depth study of the impact of elimination
orderings is also warranted.
For the particular case of 3D grid problems, \cite{CLB21} studied the
performance of many different elimination orders, however, it is not
clear their observations how this translates beyond this case.
They also developed a parallel implentation for 3D grid problems.

We believe the crucial feature of our algorithm is that it produces an
unbiased approximate factorization while always maintaining
connectivity in every elimination, and using slightly more samples
than just a tree (in the case of {\ACSM}).
However, many different elimination rules fullfill these criteria, and
we suspect that rules with even better practical performance can be
found.

We believe our solvers may finally be robust enough to use in Interior
Point Methods, and we hope this will lead to further applied research
on IPMs for single-commodity flow problems.

To further improve SDD solvers and to bring them into widespread practical use for applications such as IPMs, we believe more testing is required.
Adding more benchmark data sets could serve as an important guide to which algorithms perform well in practice in these settings.

\subsection*{Organization of the paper}
In Section~\ref{sec:prelim}, we introduce some standard facts
about SDDM matrices, Laplacian matrices, and linear algebra relating to
Cholesky factorization of these.

In Section~\ref{sec:elimtreesamp}, we introduce the pseudo-code for
our approximate Cholesky factorization.
We first describe a version which outputs an approximate factorization
in standard lower triangular form.
In Section~\ref{sec:variants}, we show how to introduce multi-entry
(a.k.a. multi-edge) splitting and merging into the approximate
factorization algorithm.
Combined with preconditioned conjugate gradient, this yields our
{\ACSM} implementation.
In Section~\ref{sec:rowform}, we describe an alternative
representation of the output factorization using row operations.
Our {\AC} implementation combines this row operation format with the
preconditioned conjugate gradient to obtain a linear equation solver.

In Section~\ref{sec:experiments}, we  describe our 
implementations and compare them to other solvers.

\subsection*{Acknowledgements}
The authors want to thank Sushant Sachdeva, Richard Peng, Joel Tropp,
George Biros, Houman Owhadi, and John Gilbert for many helpful conversation about
practical solvers for SDDM linear equations.
We thank Sushant Sachdeva for suggesting the Sachdeva Star to break the {\AC} solver.
We also thank Kaan Oktay for assistance with setting up PETSc and
running experiments.

%% file: summarytable-times.tex
\resizebox{\textwidth}{!}{%
\begin{threeparttable}
\begin{tabular}{
      |c
    | S[table-format=2.3] S[table-format=2.2]  
    | S[table-format=2.3] S[table-format=2.3] S[table-format=2.2]  S[table-format=2.2]
    |
}
\hline
    \multirow{3}{*}{Instance family}
&  {AC}   & {{\ACSM}}     
& {CMG}  &  {HyPre}   &  {PETSc}        & {ICC}    
    \\
\cline{2-7} 
    &
    \multicolumn{2}{c|}{$t_{\text{total}} / \nnz $ }  
    &
    \multicolumn{4}{c|}{$t_{\text{total}} / \nnz $ }  
    \\
\cline{2-7} 
& {$(\mu s)$} & {$(\mu s)$} 
& {$(\mu s)$}  &  {$(\mu s)$}  & {$(\mu s)$}  & {$(\mu s)$}\\ \hline %
{SuiteSparse Matrix Collection}&0.994&1.4&Inf&Inf&27.1&153\\
{Unweighted Chimeras}&3.95&7.19&Inf&Inf&Inf&Inf\\
{Weighted Chimeras}&4.39&4.79&Inf&Inf&N/A&Inf\\
{Unweighted SDDM Chimeras}&3.4&5.44&Inf&Inf&Inf&Inf\\
{Weighted SDDM Chimeras}&4.34&4.36&Inf&Inf&N/A&Inf\\
{Maximum flow IPMs}&1.2&2.39&1.9\tnote{**}&5.04\tnote{*}&N/A&Inf\\
{Sachdeva stars}&4.76\tnote{**}&1.04&0.868&3.52\tnote{**}&7.22\tnote{**}&Inf\\
{SPE benchmark}&1.62&1.91&0.511&0.648&2.58&4.29\\
{Uniform coefficient Poisson grid}&1.46&2.49&0.624&0.587&2.62&2.22\\
{High contrast coefficient Poisson grid}&2.33&3.57&0.572&0.752&3.85&4.76\\
{Anisotropic coef. Poisson grid, variable discretization}&1.48&2.5&0.659&0.614&2.6&1.92\\
{Anisotropic coef. Poisson grid, variable weight}&1.67&2.5&0.674&0.772&4.92&2.44\\
    \hline
  \end{tabular}
\begin{tablenotes}
\footnotesize
\item[*] relative residual error exceeded tolerance $10^{-8}$ by a factor $(1,10^{4}]$.
\item[**] relative residual error exceeded tolerance $10^{-8}$ by a factor $(10^{4},10^{8})$.
\item[] \hspace{-1em} Inf: Solver crashed or returned a solution with relative residual error 1 (error of the trivial all-zero solution).
\item[] \hspace{-1em} N/A: Experiment omitted as the solver crashed too often with (only occurred for PETSc).
\end{tablenotes}
\end{threeparttable}
} %

%% file: prelims.tex
\section{Preliminaries}
\label{sec:prelim}

In this section, we introduce several families of matrices for which 
we are able to construct fast linear system solvers.

\paragraph{Symmetric Diagonally Dominant (SDD) matrices.} A matrix $\MM \in
\R^{n \times n}$ is said to be Symmetric Diagonally Dominant (SDD),
if it is symmetric and for each row $i$,
\begin{align}
\MM(i,i) \geq \sum_{j \neq i} \abs{\MM(i,j)}
.
\end{align}
It can be shown that every SDD matrix is positive semi-definite.

\paragraph{Symmetric Diagonally Dominant (SDDM) matrices.} A matrix $\MM \in
\R^{n \times n}$ is said to be Symmetric Diagonally Dominant (SDDM) if
it is SDD and has non-positive off-diagonal entries.

\paragraph{Graph Laplacians a.k.a. Laplacians.}%

We consider an undirected graph $G =(V,E)$, with
positive edges weights $w : E \to \R_{+}$.  Let $n = \abs{V}$ be the number of vertices and ,
$m = \abs{E}$ the number of edges, and
$V = \setof{1, \ldots, n}$. 
Let $\vecind_{i}$ denote the $i^{\textrm{th}}$
standard basis vector.  Given an ordered pair of vertices $(u,v)$, we
define the pair-vector $b_{u,v} \in \R^{n}$ as
$\bb_{u,v} = \vecind_{v} - \vecind_{u}.$ 
For a edge $e$ with endpoints $u$ and $v$ (arbitrarily ordered), we define $\bb_{e} = \bb_{u,v}.$
By assigning an arbitrary direction to each edge of $G$ we define
the Laplacian of $G$ as $\ulap = \sum_{e \in E} w(e) \bb_{e} \bb_{e}^{\top}.$
Note that the Laplacian does not depend on the choice of direction for
each edge.
Given a single edge $e$, we refer to $w(e) \bb_{e} \bb_{e}^{\top}$
as the Laplacian of $e$.
The Laplacian of a weighted, undirected graph is unique given the graph and vice versa, and we will treat these as interchangeable given this equivalence.

\paragraph{Laplacians and multi-graphs.}
For some variants of our algorithms, we will need to work with multi-graphs.
We can consider an undirected multi-graph $G =(V,E)$ with
positive multi-edges weights $w : E \to \R_{+}$.
Note that different multi-edges between the same vertices may have different weights.
In this setting, we again let $n$ denote the number of vertices and we let $m = \abs{E}$ denote the number of multi-edges.
Similar to the case of graphs without multi-edges, we can assign arbitrary directions to multi-edges to define the Laplacian of a multi-edge $e$ as  $w(e) \bb_{e} \bb_{e}^{\top}$ and define the Laplacian of $G$ as 
$\ulap = \sum_{e \in E} w(e) \bb_{e} \bb_{e}^{\top}.$
Note that different multi-graphs may have the same Laplacian. This occurs because the coefficient of the term $\bb_{(u,v)} \bb_{(u,v)}^{\top}$ only depends on the sum $\sum_{e \in E_{uv}} w(e)$ where $E_{uv}$ is the set of multi-edges between vertices $u$ and $v$.
Thus, for example, we can take a multi-graph $G$ and produce another multi-graph $G'$ with the same Laplacian by splitting each multi-edge $e$ of $G$ into $k$ copies with weight $w(e)/k$ each.

In Section~\ref{sec:variants}, we introduce algorithms that use multi-graphs when computing an approximate Cholesky factorization of a Laplacian.

\begin{fact}
If $G$ is connected, then the kernel of the corresponding Laplacian $\ulap$ is the span of 
the vector $\vecone.$
\end{fact}

\paragraph{Upper and lower triangular matrices.}
We say a square matrix $\matup$
  is upper triangular if it has non-zero entries $\matup(i,j) \neq 0$
  only for $i \leq j$ (i.e. above the diagonal). 
  Similarly, we say a square matrix $\matlow$
  is lower triangular if it has non-zero entries $\matlow(i,j) \neq 0$
  only for $i \geq j$ (i.e. below the diagonal). 
  Often, we will work with matrices that are not upper or lower
  triangular, but which for we know a permutation matrix $\PP$ s.t. 
  $\PP \matup \PP^{\trp}$ is upper (respectively lower) triangular.
  For computational purposes, this is essentially equivalent to having
  a upper or lower triangular matrix, and \emph{we will refer to such
    matrices as upper (or lower) triangular}. 
  The algorithms we develop for factorization will always compute the
  necessary permutation.

%% file: algo.tex
\section{A Linear Equation Solver for Laplacian Matrices}

Our linear equation solvers are based on Preconditioned Conjugate
Gradient (PCG), largely following \cite{BBCDDDEPRVdV94}. Our approach to testing for stagnation is largely based on the MATLAB PCG implementation.
We also refer the reader to the textbook by Golub and Van Loan
\cite[Chapter 11]{GVL13} for additional background.
We assume the reader is familiar with basics of PCG and preconditioners.

Our goal in this section is to describe our algorithm(s) for computing approximate Cholesky factorizations of  Laplacian matries, which we then use as preconditioners inside the PCG.
We use the overall algorithmic framework of \cite{KS16}, but we introduce a new approximate elimination rule that seems to perform much better in practice. In particular, it yields a good preconditioner with a much lower sampling rate than that of \cite{KS16}.

We reduce the problem of solving linear equations in SDDM matrices to that of solving linear equations in Laplacian matrices by using a technique introduced by 
Gremban \cite[Lemma 4.2]{G96}:
An SDDM Matrix $M$ can be written in the form $D + L$ where $L$ is a Laplacian and $D$ is a non-negative diagonal matrix.
Let $d$ be the vector of entries on the diagonal of $D$.
Gremban constructs a Laplacian matrix $\widehat{L}$ that is identical to $L$ except that it has one extra vertex that is connected to vertex $i$ by an edge of weight $d[i]$.
To solve the system $Mx = b$, we first solve $\widehat{L} y = \hat{b}$ where $\hat{b}$ is identical to $b$ except at the extra vertex, where its value is set to make the sum of the entries of $\hat{b}$ zero.
The vector $x$ is then obtained by subtracting the value of $y$ at the extra vertex from the values of $y$ at the original vertices.

\subsection{Background: Approximate Cholesky Factorization of Laplacian Matrices}
\label{sec:apxChol}

In this section, we briefly review the how to compute a Cholesky factorization of a Laplacian matrix, using the additive formulation of \cite{KS16}.
We then briefly discuss how the \emph{elimination clique} that arises from each row and column elimination can be sampled and how this yields a fast algorithm for approximate Cholesky decomposition \cite{KS16}.
This approximate Cholesky decomposition may then be used as a preconditioner in PCG to build a linear equation solver.

\paragraph{Producing a Cholesky factorization.}
Cholesky factorization expresses a matrix in the form $\matlow \matlow^{T}$.
A Cholesky factorization can be found for every symmetric positive definite matrix,
  and for every positive semidefinite matrix after a possible permutation of its rows and columns.
One perspective on Cholesky factorization is that it proceeds by
writing a symmetric positive semi-definite matrix as a sum of rank 1
matrices. We will outline this view here.

Let $\MM$ be a symmetric positive semi-definite matrix.
As a convenient notational convention, we define
$\SS^{(0)} \defeq \MM$, the ``0th'' Schur complement. 
We then define 
\begin{align}
\nonumber
  \ll_{i} &=\frac{1}{\sqrt{\SS^{(i-1)}(i,i)}} \SS^{(i-1)}(:,i) \\
\nonumber
  \SS^{(i)} &= \SS^{(i-1)} - \ll_{i}\ll_{i}^{\trp}
,
\label{eq:schurIterDef}
\end{align}
unless $\SS^{(i-1)}(i,i) = 0$, $\SS^{(i-1)}(:,i) = \veczero$, and $\SS^{(i-1)}(i,:) =
\veczero$, in which case we define $\SS^{(i)} =
\SS^{(i-1)}$. 
We do not define the Schur complement when the diagonal is zero but
off-diagonals are not.
By setting
\[
\matlow =
\begin{pmatrix}
  \ll_{1} & \ll_{2} & \cdots & \ll_{n}
\end{pmatrix}
\]
we get a Cholesky factorization $\MM = \matlow\matlow^{\trp}$.

It can be shown that the algorithm described above always produces a
Cholesky factorization when applied to a Laplacian.
For some other classes of matrices, it does not always
succeed and pivoting may be required (e.g. if the first column has a zero diagonal but other non-zero entries).

\paragraph{The clique structure of Schur complements}
In this section, we recall a standard proof that the Schur complement
of a Laplacian is another Laplacian, while also observing that the Schur
complement of a Laplacian onto $n-1$ indices has additional structure
that will help us develop algorithms for approximating these Schur
complements.

Given a Laplacian $\ulap$,
let $\vstar{\ulap}{v} \in \R^{n \times n}$ denote the Laplacian corresponding to the 
edges incident on vertex $v$ (the \emph{star} on $v$), i.e.
\begin{equation}
\label{eq:starDef}
\vstar{\ulap}{v} \defeq \sum_{e \in E : e \ni v} w(e) \bb_{e} \bb_{e}^{\top}
.
\end{equation}
For example, if the first column of $\ulap$ is
$
\begin{pmatrix}
  d \\
  -\aa 
\end{pmatrix}
,$
then
$
\vstar{\ulap}{1} =
\begin{bmatrix}
d & -\aa^{\top} \\
-\aa & \diag(\aa)
\end{bmatrix}
.
$
We can write the Schur complement 
$\schurto{\ulap}{[n]\setminus\setof{v}}$ as
\[
\schurto{\ulap}{[n]\setminus\setof{v}} = \ulap-\vstar{\ulap}{v} +\vstar{\ulap}{v}  -
\frac{1}{\ulap(v,v)}\ulap(:, v) \ulap(v, :)
.
\]
It is immediate that $\ulap-\vstar{\ulap}{v} $ is a Laplacian matrix,
since $\ulap-\vstar{\ulap}{v}  = \sum_{e \in E : e \not\ni v} w(e) \bb_{e}
\bb_{e}^{\top}$.
A more surprising (but well-known) fact is that 
\begin{align}
\label{eq:Cv-def}
\cliq{\ulap}{v} \defeq \vstar{\ulap}{v}  - \frac{1}{\ulap(v, v)}\ulap(:, v) \ulap(v, :)
\end{align}
is also a Laplacian, and its edges form a clique on the neighbors of $v$.
It suffices to show it for $v = 1.$ We write $i \sim j$ to indicate $(i,j) \in E.$
Then
\begin{equation*}
\label{eq:cliquestructure}
\cliq{\ulap}{1}
=
\begin{bmatrix}
\matzero & \veczero^{\top} \\
\veczero & \diag(\aa)- \frac{\aa\aa^{\top} }{d}
\end{bmatrix}
= 
\sum_{i \sim 1} \sum_{j \sim 1, i < j} 
\frac{w(1,i) w(1,j)} 
{d}
\bb_{(i,j)} \bb_{(i,j)}^{\top}
.
\end{equation*}
Thus
$\schurto{\ulap}{[n]\setminus\setof{v}}$ is a Laplacian since it is a sum of two Laplacians.
By induction, for all $C \subseteq [n] $, $\schurto{\ulap}{C}$ is a Laplacian.

\paragraph{Laplacian Cholesky factorization with cliques.}
We can use the clique structure described above to write Cholesky factorization of a
Laplacian matrix in a convenient form. We state the pseudo-code below
in Algorithm~\ref{alg:LapCholeskyFactorization}.

\begin{algorithm}
  \caption{Algorithm
    \textsc{Cholesky}{($\LL$)} outputs a lower
    triangular Cholesky factor $\matlow \in \R^{n \times n}$ of the input Laplacian
    ${\LL \in \R^{n \times n}}$.}
  \label{alg:LapCholeskyFactorization}
  \begin{algorithmic}[1]
    \Procedure{Cholesky
    }{$\LL$}
    \State $\SS \gets \LL$
     \For{$v=1$ to $n-1$}  \Comment the order can be chosen adaptively
    \State $\ll_{v} \gets \frac{1}{\sqrt{\SS(v,v)}} \SS(:,v)$ \Comment
    $\SS(:,v)$ is the $v$th column of $\SS$
    \State $\SS \gets \SS  - \vstar{\SS}{v} + \textsc{Clique}(v,\SS)$ \label{lne:cholCliqueAdd}
    \EndFor 
    \State \Return
   $
    \matlow =
    \begin{pmatrix}
      \ll_{1} & \ll_{2} & \cdots & \ll_{n-1}  & \veczero
    \end{pmatrix}
    $
    \Comment
    The final column should be the all-zero vector.
    \EndProcedure
  \end{algorithmic}
\end{algorithm}

\paragraph{A Framework for Approximate Cholesky Factorization via Clique Sampling.}
\cite{KS16} introduced the idea of constructing a preconditioner by 
  replacing the elimination clique in
  Cholesky factorization with an unbiased sampled approximation of it.
The major factors governing the performance of such a preconditioner are the sampling procedure and the order in which vertices are eliminated.
For now, we state the generic algorithm by 
 using $\textsf{CliqueSample}$ as a placeholder for the sampler and an arbitrary order.
This gives
  the generic algorithm for approximate
  Cholesky factorization
  stated in Algorithm~\ref{alg:GeneralApproximateCholeskyFactorizationColwise}

\begin{algorithm}
  \caption{Algorithm
    \textsc{ApproximateCholesky}{($\LL$)} outputs a lower
    triangular matrix $\matlow \in \R^{n \times n}$ that gives an approximate Cholesky
    factorization of the input Laplacian
    ${\LL \in \R^{n \times n}}$.}
  \label{alg:GeneralApproximateCholeskyFactorizationColwise}
  \begin{algorithmic}[1]
    \Procedure{ApproximateCholesky
    }{$\LL$}
    \State $\SS \gets \LL$
     \For{$v=1$ to $n-1$} \Comment the order can be chosen adaptively
    \State $\ll_{v} \gets \frac{1}{\sqrt{\SS(v,v)}} \SS(:,v)$ \Comment
    $\SS(:,v)$ is the $v$th column of $\SS$
    \State $\SS \gets \SS  - \vstar{\SS}{v} + \textsf{CliqueSample}(v,\SS)$
    \EndFor 
    \State \Return
   $
    \matlow =
    \begin{pmatrix}
      \ll_{1} & \ll_{2} & \cdots & \ll_{n-1}  & \veczero
    \end{pmatrix}
    $
    \Comment
    The final column should be the all-zero vector.
    \EndProcedure
  \end{algorithmic}
\end{algorithm}

Formally, the requirement that the sampler be unbiased means that 
 \[
    \expct{}{\textsf{CliqueSample}(v,\SS)} = \cliq{\SS}{v}.
 \]
\cite{KS16} proved that unbiased clique sampling produces 
  a factorization that in expectation equals the original input matrix
 \footnote{See \cite{KS16} Section 4.1.
They considered a specific sampling rule, but their proof uses no
properties of the sampling rule except that it is unbiased.}.

They then showed that
\begin{claim}
\label{clm:unbiasedchol}
  When $\textsc{ApproximateCholesky}$ is instantiated with an unbiased clique
  sampling routine $\textsf{CliqueSample}$ whose output is positive semidefinite, 
  its output 
  $\matlow = \textsc{ApproximateCholesky}(\LL)$ satisfies
  \[
    \expct{}{\matlow \matlow^{\trp}} = \LL
    .
  \]
\end{claim}
For completeness, we include a proof in Appendix~\ref{sec:unbiased}.

We will introduce a sampling procedure $\textsc{CliqueTreeSample}$ that outputs at most as many multi-edges as the multi-degree of the vertex being eliminated.
The analysis of \cite{KS16} may be used to bound the running time of approximate elimination with a clique sampler that satisfies such a guarantee.
Suppose the clique sampler has the property that it always outputs at most as many (multi-)edge samples as the (multi-edge) degree of the vertex $v$ in $\SS$.
If $\textsc{ApproximateCholesky}$ is then run with an elimination ordering that always picks a vertex whose degree in the remaining graph is upper bounded by a constant multiple of the average degree 
  of vertices in the remaining graph, then the number of non-zeros in the output $\matlow$ will be at most
$O(|E|\log|E|)$.
This bound always holds, and does not depend on the random choices made by the clique sampler.
If the degree bound holds in expectation for each elimination, the non-zero bound will hold in expectation as well.
We state this claim formally below.

\begin{claim}
  \label{clm:fillin}
  Suppose $\textsf{CliqueSample}(v,\SS)$ outputs at most as many multi-edge samples
  as the multi-edge degree of the vertex $v$ in $\SS$.
  Then if $\textsc{ApproximateCholesky}$ is run with an elimination ordering that ensures that when a vertex $v$ is eliminated, we have
  \[
  \text{degree of $v$ in current graph} \leq O(\text{average degree in current graph}).
  \]
  then $\textsc{ApproximateCholesky}$ always outputs a lower triangular matrix $\matlow$ with $O(|E|\log|E|)$ non-zeros.
In fact, it suffices for the degree of $v$ in the current graph to be bounded by a constant times the number of multi-edges in the original graph divided by the number of vertices in the current graph.
  
  If the bound on the degree of $v$ holds at each elimination step in expectation over a random choice of $v$, then $\textsc{ApproximateCholesky}$ outputs a lower triangular matrix $\matlow$ with $O(|E|\log|E|)$ non-zeros in expectation.
\end{claim}
For completeness, we include a proof of the claim.
\begin{proof}
  Observe that the total number of edges in $\SS$ is non-increasing
  over time, as each elimination removes the edges incident on the
  vertex being eliminated and the $\textsc{CliqueSample}$ call adds back at most
  the same number.

  Hence $\SS$ contains at most $|E|$ edges.
  The number of non-zeros in column $\ll_i$ of $\matlow$ is
  proportional to the degree of the vertex being eliminated in the current graph of
  $\SS$.
  This is hence on average equal to the average degree in $\SS$, which
  is at most $\frac{2|E|}{n-i}$.
  Thus the total number of non-zeros in $\matlow$ is bounded
  by
  \[
    \sum_{i = 1}^{n-1} \frac{2|E|}{n-i} = O(|E|\log|E|)
    .
  \]
  If the bound per elimination holds in expectation, the final bound $O(|E|\log|E|)$ will again hold, but only in expectation. 
\end{proof}

\paragraph{Formal analysis of approximate elimination.} 
In \cite{KS16}, the authors described a clique sampling procedure with the property that if $\LL$ is the Laplacian of a graph where every (multi-)edge has leverage score at most $O(1/\log^2|V|)$, and a random elimination order is used, then with high probability, the output decomposition $\matlow \matlow^{\trp}$ has relative condition number with $\LL$ of at most 4. This in turn will ensure that PCG converges to a highly accurate solution in few steps when used to solve a linear equation in $\LL$.
By sub-dividing edges into multi-edges, every Laplacian with $m$ off-diagonal non-zeros can be associated with a graph with $O(m\log^2|V|)$ multi-edges that each have leverage score $O(1/\log^2|V|)$.
Thus, the clique sampling procedure of \cite{KS16} can produce a good approximate Cholesky decomposition in time $O(m\log^3 m)$.

\subsection{Edgewise Elimination and A New Clique Sampling Rule}
\label{sec:edgewise}
In this section, we describe how the elimination clique created in each step of Cholesky factorization can in fact be further decomposed into a sum of \emph{elimination stars} on the neighbors of the vertex being eliminated.
We then show how subsampling each elimination star gives rise to a clique sampler with a desirably property: The output is always a connected graph.
We describe two variants of our sampling rule. The simplest approximates the elimination clique with a tree, and hence we call it $\textsc{CliqueTreeSample}$.
Our second variant approximates the elimination clique with a tree and potentially a few extra edges, roughly as many as appear in the tree.
This variant arises naturally from combining our $\textsc{CliqueTreeSample}$ approach with the multi-edge splitting that is used in \cite{KS16} to improve the approximation quality by doing more fine-grained sampling. 
We call this $\textsc{CliqueTreeSampleMultiEdge}$, and a slight but important tweak of it is called  $\textsc{CliqueTreeSampleMultiEdgeMerge}$.
This tweak limits the maximum number of multi-edges between two vertices, which gives large performance gains in practice.

\paragraph{Decomposing the elimination clique into elimination stars.}
When we eliminate a vertex, we can further decompose the resulting
clique $\cliq{\ulap}{v}$ as a sum of smaller graphs.
This decomposition may appear a little mysterious, but it arises
naturally if we consider the process of eliminating a column as a
sequence of eliminations that each remove one entry of the column.
To arrive at the decomposition, we pick an ordering on the
neighbors of the vertex being eliminated.
To simplify the discussion, we consider eliminating the first vertex
$v = 1$, and the order we pick for eliminating the neighbors is by
increasing vertex number. However, any vertex can be eliminated this
way and any ordering on the neighbors can be used.
As before, we let the first column of the Laplacian $\LL \in \R^{n \times
  n}$ be by
$
\begin{pmatrix}
  d \\
  -\aa (2:n)
\end{pmatrix}$
where $\aa
\in
\R^{n}$ has first entry 0.
We also define $\aa_i \in \R^{n}$ given by
\[
  \aa_i(j)=
  \begin{cases}
\aa(l) & \text{for } j>i \\
0 &  \text{ for } j \leq i
\end{cases}
\]
Thus $\aa_i$ is the off-diagonal entries of the first column, after
removing the first $i$ entries.

We can then break $\cliq{\ulap}{1}$ into a sum of smaller terms, where
each term is a graph Laplacian consisting of edges from neighbor $i$
of vertex $1$ to neighbors $l \geq i$.

This smaller graph Laplacian is given by
\begin{equation}
\label{eq:schurfromsingleedge}
\estar{\LL}{1}{i}
\defeq
\begin{pmatrix}
 \diag\left(\frac{\aa(i) \vecone^\trp\aa_i}{d} \ee_i + \frac{\aa(i)}{d} \aa_i \right)
  - \frac{\aa(i)}{d}\left( \ee_i\aa_i^\trp+\aa_i\ee_i^\trp \right)
\end{pmatrix}
\end{equation}
We can also write this as 
\begin{equation}
\label{eq:schurfromsingleedgesimple}
\estar{\LL}{1}{i}
=
\begin{pmatrix}
  \sum_{j> i}
  \frac{\aa(j) \aa(i)}{d} (\ee_i - \ee_j) (\ee_i - \ee_j)^{\trp}
\end{pmatrix}
\end{equation}

We call $\estar{\LL}{1}{i}$ the \emph{elimination star} of neighbor $i$.
It is straightforward to verify that $\estar{\LL}{1}{i}$ is a graph Laplacian and
it is non-zero only when $\aa(i) \neq 0$ and consists of edges from
the $i$th neigbor neighbors $l \geq i$ where $\aa(l) \neq 0$.
Note that $\estar{\LL}{1}{i}$ depends on the ordering chosen on the neighbors, and
again, any ordering can be used.

Furthermore, these Laplacians together add up to give the clique created
by eliminating vertex $1$.
This is immediate from the fact that each $S_i$ contains exactly the
off-diagonal entries of $\cliq{\ulap}{1}$ for $i$ and for each $j >
i$, and thus summing up across all $i$, we get $\cliq{\ulap}{1}$.
Because $\SS_i$ is nonzero only when $i$ is a neighbor of vertex 1 in
the current graph, we can write
\begin{equation}
  \label{eq:edgestars}
  \cliq{\ulap}{1}  = \sum_{i\sim 1} \estar{\LL}{1}{i}
\end{equation}

\paragraph{Approximating the elimination clique by sampling a tree.}
\label{sec:elimtreesamp}
Like \cite{KS16}, we will approximate the elimination clique by
sampling.
However, we engineer our sampler so that each elimination star
$\estar{\LL}{1}{i}$ is approximated by a single randomly chosen reweighted edge,
chosen so that in expectation, this edge Laplacian yields the
elimination star $\estar{\LL}{1}{i}$.

We choose a random index 
$\gamma_i$ according to a probability distribution given by
\[
\Pr[ \gamma_i = j ] =
p_i(j) \text{ where } p_i(j)= \frac{\aa_i(j)}{\vecone^\trp\aa_i}
.
\]
We can also write this as
\[
p_i(j)  = 
\begin{cases}
\frac{\aa(j)}{\sum_{l > i} \aa(l)}
& \text{for } j>i \\
0
&  \text{ for } j \leq i
\end{cases}
.
\]
To make the expectations work out,
if the outcome is $\gamma_i = j$, 
we then choose a weight of 
\begin{equation}
\tilde{w}(i,j) = \frac{1}{p_i(j)} \frac{1}{d} \aa(j) \aa(i)
= \frac{\aa(i)\vecone^\trp\aa_i}{d}
.\label{eq:edgesampleweight}
\end{equation}
and we output the single edge Laplacian

\begin{equation}
  \SStil_i =
   \tilde{w}(i,j)  (\ee_i - \ee_j) (\ee_i - \ee_j)^{\trp}\label{eq:edgesample}
 \end{equation}
 
We can show that $\expct{\gamma_i}{\SStil_i} = \estar{\LL}{1}{i}$, i.e. the
sampling of the elimination star is correct in expectation.

\begin{claim}
  \label{clm:expctestarsample}
  $\expct{\gamma_i}{\SStil_i} = \estar{\LL}{1}{i}$
\end{claim}
\begin{proof}
By Equation~\eqref{eq:edgesample} and
Equation~\eqref{eq:schurfromsingleedgesimple}, we have that
\[
  \expct{\gamma_i}{\SStil_i} =
  \sum_{j > i} p_i(j)  \tilde{w}(i,j)  (\ee_i - \ee_j) (\ee_i - \ee_j)^{\trp}
  =
  \sum_{j > i} \frac{1}{d} \aa(j) \aa(i) (\ee_i - \ee_j) (\ee_i - \ee_j)^{\trp}
  =
  \estar{\LL}{1}{i}.
  \]
\end{proof}

We now approximate the whole elimination clique
$\cliq{\LL}{1}$ by the union of one sample for each of the elimination
stars, leading to an approximation of the clique by 
$\sum_{i\sim 1} \SStil_i \approx \cliq{\LL}{1}$.
We summarize the pseudo-code for this sampling procedure
  in Algorithm~\ref{alg:CliqueTreeSample}, which we
 call $\textsc{CliqueTreeSample}$ for reasons that will be
apparent in a moment.

\begin{algorithm}[H]
  \caption{Algorithm 
    \textsc{CliqueTreeSample}{$(v,\SS)$} returns $\CCtil$
    which approximates the elimination clique $\cliq{\SS}{v}$}
  \label{alg:CliqueTreeSample}
  \begin{algorithmic}[1]
    \Procedure{CliqueTreeSample}{$v,\SS$}
    \State $\CCtil \gets \matzero \in \R^{n \times n}$
    \State $d \gets \SS(v,v)$; $\aa \gets -\SS(v,:)$; $\aa(v) \gets 0$
    \ForAll{$i$ s.t. $\aa(i) \neq 0$}      \Comment pick any ordering on
    the neighbors
    \State $c \gets \aa(i)$; $\aa(i) \gets 0$
    \State Sample index $j$ with probability $p(j) = \frac{\aa(j)}{\vecone^\trp\aa_{i}}$
    \State $\CCtil \gets \CCtil + c \cdot\frac{ \vecone^\trp\aa}{d} (\ee_i - \ee_j) (\ee_i - \ee_j)^{\trp}$
    \EndFor 
    \State \Return $\CCtil$
    \EndProcedure
  \end{algorithmic}
\end{algorithm}

We can directly conclude from
Equation~\eqref{eq:schurfromsingleedgesimple} and
Claim~\ref{clm:expctestarsample} that, in expectation, the output of \textsc{Clique\-Tree\-Sample}{(v,\SS)} 
equals $\cliq{\SS}{v}$.

\begin{corollary}
  \label{cor:expctcliquesample}
    \[
    \expct{}{\textsc{CliqueTreeSample}(v,\SS)} = \cliq{\SS}{v}
  \]
\end{corollary}

It turns out that $\textsc{CliqueTreeSample}$ always outputs the
Laplacian of a tree of the neighbors of the vertex being
eliminated. 

\begin{claim}
  The random matrix returned by $\textsc{CliqueTreeSample}{(v,\SS)}$ is always
  the Laplacian of a tree on the neighbors of $v$ in the graph
  associated with $\SS$.
\end{claim}

\begin{proof}
  We can trivially ignore all the vertices that are not neigbors of
  vertex $v$, and w.l.o.g. assume the neighbors of $v$ are numbered $1,\ldots,k$.
  The result follows by a simple induction, starting from $i = k-1$
  and down to $i = 1$.
  Note that the $i$th sample connects vertex $i$ to a vertex $j >i$.
  Our induction hypothesis is that samples $i,\ldots,k$ form a tree
  on vertices $i,\ldots,k$.
  The base case $i = k-1$ is trivially true as the last sample connects vertices
  $k-1$ and $k$ by a single edge (deterministically).
  For the inductive step, observe that the $i$th sample connects to one of the vertices
  $i+1,\ldots,k$, where by induction, we already have a tree.
  Thus $i$ is now connected by an edge to the tree on vertices
  $i+1,\ldots,k$. And the new graph must again form a tree as it
  connects $k-i+1$ vertices using $k-i$ edges.
\end{proof}

\paragraph{Approximate Cholesky factorization.}
We plug in $\textsc{CliqueTreeSample}$ 
  for \textsf{CliqueSample} as our choice of clique sampling routing in Algorithm~\ref{alg:GeneralApproximateCholeskyFactorizationColwise} to obtain a new algorithm for computing an approximate Cholesky factorization.
By Corollary~\ref{cor:expctcliquesample}, the sampling rule is unbiased, and hence by Claim~\ref{clm:unbiasedchol}, the overall output satisfies $\expct{}{\matlow \matlow^{\trp}} = \LL$.

We eliminate the vertices in an order that guarantees that the degree of the vertex being eliminated is always at most twice the number of multi-edges in the original graph divided by the number of vertices in the current graph.
Thus, Claim~\ref{clm:fillin} applies and the number of non-zeros in the output decomposition is bounded by $O(m\log m)$ when the input Laplacian has $O(m)$ non-zeros.

While we expect that one could obtain better preconditioners by always eliminating a vertex of lowest degree in the current graph, the priority queue required to implement this would slow down the algorithm too much. We use a lazy approximate version of this rule because it is faster.

One should also expect the order of elimination of stars within a clique to impact the quality 
  of the preconditioner.
We eliminate starts in order
 of lowest to highest edge weight.

\subsection{A More Accurate Sampling Rule Using Multi-Edge Sampling}
\label{sec:variants}
In this section, we introduce two variants of the $\textsc{CliqueTreeSample}$ sampling rule.
The first variant is conceptually simplest but performs poorly in practice, necessitating our introduction of the second variant, which is slightly more complex but performs much better in practice.

\paragraph{Approximate Cholesky factorization with edge splits and multi-edges.}
Inspired by \cite{KS16}, our first variant splits the edges of the graph into $k$ copies of multi-edge, each with
$\frac{1}{k}$ of the original edge weights.
Note that the resulting multi-graph has the same Laplacian as the original graph, the associated Laplacian only depends on the sum of multi-edge weights between each pair of vertices.

While Algorithm~\ref{alg:GeneralApproximateCholeskyFactorizationColwise} works with a Laplacian and its associated graph, our next algorithms will explicitly maintain a multi-graph and its associated Laplacian.
Because of the initial splitting of the edges and sampling based on a larger number of multi-edges, 
the variants we introduce in this section is more robust and is observed to produce better approximations of the Cholesky factorization. 
However, this comes at the cost of a higher number of non-zeroes in the output factorization and longer running time.

We start by introducing a new overall factorization algorithm framework, Algorithm~\ref{alg:GeneralApproximateCholeskyFactorizationColwiseMultiEdge}, which computes a factorization of the input Laplacian by calling a sampling routine $\textsc{CliqueSampleMultiEdge}$ that uses multi-edges when sampling.
We then provide two different of variants of $\textsc{CliqueSampleMultiEdge}$, stated in  Algorithms~\ref{alg:CliqueTreeSampleMultiEdge} and~\ref{alg:CliqueTreeSampleMultiEdgeMerge}.

\begin{algorithm}
  \caption{Algorithm
    \textsc{ApproximateCholeskyMultiEdge}{($\LL$)} outputs a lower
    triangular matrix $\matlow \in \R^{n \times n}$ that gives an approximate Cholesky
    factorization of the input Laplacian
    ${\LL \in \R^{n \times n}}$.}
  \label{alg:GeneralApproximateCholeskyFactorizationColwiseMultiEdge}
  \begin{algorithmic}[1]
    \Procedure{ApproximateCholeskyMultiEdge
    }{$\LL$}
    \State $\SS \gets \LL$ and $G \gets$ a multi-graph with $\LL$ as its Laplacian.
    \State 
    \Comment{Different algorithm variants may use different rules for computing the multi-graph $G$ given $\LL$}
     \For{$v=1$ to $n-1$} \Comment the order can be chosen adaptively
    \State $\ll_{v} \gets \frac{1}{\sqrt{\SS(v,v)}} \SS(:,v)$ \Comment
    $\SS(:,v)$ is the $v$th column of $\SS$
    \State $\SS \gets \SS  - \vstar{\SS}{v} + \textsc{CliqueSampleMultiEdge}(v,\SS,G)$
    \State
        \Comment{$\textsc{CliqueSampleMultiEdge}(v,\SS,G)$ updates the multi-graph $G$ and the updated $\SS$ is the associated Laplacian.}
    \EndFor 
    \State \Return
   $
    \matlow =
    \begin{pmatrix}
      \ll_{1} & \ll_{2} & \cdots & \ll_{n-1}  & \veczero
    \end{pmatrix}
    $
    \Comment
    The final column should be the all-zero vector.
    \EndProcedure
  \end{algorithmic}
\end{algorithm}

Our first multi-edge-based clique sampling algorithm is Algorithm~\ref{alg:CliqueTreeSampleMultiEdge}, which we refer to as \textsc{Clique\-Tree\-Sample\-Multi\-Edge}.

\begin{algorithm}[H]
    \caption{Algorithm 
      \textsc{CliqueTreeSampleMultiEdge}{$(v,\SS,G)$} updates the multi-graph $G$ to eliminate vertex $v$ and add an approximate elimination clique with multi-edges in its place
      \emph{ and } 
      returns $\CCtil$
      which approximates the elimination clique $\cliq{\SS}{v}$.
      }
    \label{alg:CliqueTreeSampleMultiEdge}
    \begin{algorithmic}[1]
      \Procedure{CliqueTreeSampleMultiEdge}{$v,\SS,G$}
      \State $\CCtil \gets \matzero \in \R^{n \times n}$
      \State 
      \label{lne:nbrsdef}
      Let $N(G,v)$ be the set of vertices that neighbor $v$ in $G$
      \State Let $E(G,v,i)$ be the set of multi-edges between $v$ and $i$ in $G$
      \State Let $w(G,v,i) = \sum_{e \in E(G,v,i)} w(e)$ be the sum of weights of multi-edges between $v$ and $i$ in $G$
      \State
      \label{lne:wdef}
      Let $w(G,v) = \sum_{i \in N(G,v)} w(G,v,i)$  be the sum of weights of multi-edges incident to $v$ in $G$
      \State \Comment{The quantities above are updated as $G$ changes}
      \State $d \gets w(G,v)$ it \emph{initial} weight of multi-edges incident to $v$
      \State \Comment{$d$ is \emph{not} updated as $G$ changes}
      \ForAll{vertices $i \in N(G,v)$}      \Comment{pick any ordering on
      the neighbors}
      \State $t \gets \abs{E(G,v,i)}$
      \State $\bar{w}_{vi} \gets w(G,v,i)/t$
      \Comment{Average of the multi-edge weights between $v$ and $i$ in $G$}
      \State Delete multi-edges $E(G,v,i)$ from $G$.
      \State 
      \Comment{Update values defined in Lines~\ref{lne:nbrsdef}-\ref{lne:wdef} correspondingly}
      \State $w_{\text{new}} \gets \bar{w}_{vi} \cdot \frac{w(G,v)}{d}$
      \Comment{Weight used for samples}
      \For{$h = 1$ to $t$}
      \State Sample index $j$ with probability $p(j) = \frac{w(G,v,j)}{w(G,v)}$
      \State $\CCtil \gets \CCtil + w_{\text{new}} \cdot (\ee_i - \ee_j) (\ee_i - \ee_j)^{\trp}$
      \State Add a multi-edge between $i$ and $j$ of weight $w_{\text{new}}$ in to $G$
      \EndFor 
      \EndFor
      \State \Return $\CCtil$
      \EndProcedure
    \end{algorithmic}
\end{algorithm}

Algorithm~\ref{alg:CliqueTreeSampleMultiEdge} can be instantiated with different orderings on the neighbors of $v$.
We use an ordering by increasing weight $w(G,v,i)$.
This is motivated by numerical stability of the sampling and heuristic arguments about global variance. It is likely other orderings may work well too.

Using an ordering by increasing weight also means the algorithm can be implemented in expected linear time in the degree of $v$ by using the alias method for sampling, and combining with the rejection sampling to avoid updating the sampling weights until a constant fraction of entries have been processed.
However, in practice we opt for maintaining probabilities using a binary search tree, with $O(k \log k)$ overall complexity processing a vertex of degree $k$.

For each edge being eliminated, we always average the weights of its multi-edges. This step can be omitted, but in practice 
we observe that this step makes the code numerically more stable.
We also expect it improves spectral approximation slightly as the averaged samples have lower maximum leverage score.

Like our other clique sampling routines, Algorithm~\ref{alg:CliqueTreeSampleMultiEdge} is correct in expectation, as stated in the following claim.
\begin{claim}
The output of Algorithm~\ref{alg:CliqueTreeSampleMultiEdge} satisfies
    \label{cor:expctcliquesamplemultiedge}
      \[
      \expct{}{\textsc{CliqueTreeSampleMultiEdge}(v,\SS,G)} = \cliq{\SS}{v}
    \]
\end{claim}
The proof is simple and similar to the proof in the simple graph case, and hence we omit it. 

\paragraph{Approximate Cholesky factorization with edge splits and multi-edge merging.}
We provide a second \textsc{CliqueSampleMultiEdge} routine which 
is a compromise between the previous variant and our basic version of \textsc{CliqueTreeSample} (Algorithm~\ref{alg:CliqueTreeSample}).
  This variant also splits the edges into $k$ copies of multi-edge but only keeps track of the multi-edges up to some limit $l$. In other words, if there are more than $l$ multi-edges for the same edge, this variant merges them down to $l$ copies.
    This compromise works well in practice, as shown by our experiments in Section~\ref{sec:experiments}.
    We believe this is 
    because it allows fine-grained sampling than Algorithm~\ref{alg:CliqueTreeSample} without leading to a build-up of large numbers of multi-edges late in the elimination as the graph gets denser as can happen in Algorithm~\ref{alg:CliqueTreeSampleMultiEdge}.
  Algorithm~\ref{alg:CliqueTreeSample} is equivalent to a special case of this variant where we set both $k$ and $l$ to be $1$, while
  Algorithm~\ref{alg:CliqueTreeSampleMultiEdge} can be recovered by setting $l$ to be infinity.

  We summarize the \textsc{CliqueSample} procedure of this 
  variant as Algorithm~\ref{alg:CliqueTreeSampleMultiEdgeMerge}.
  
  \begin{algorithm}[H]
    \caption{Algorithm 
      \textsc{CliqueTreeSampleMultiEdgeMerge}{$(l,v,\SS,G)$} updates the multi-graph $G$ to eliminate vertex $v$ and add an approximate elimination clique with multi-edges in its place
      \emph{ and } 
      returns $\CCtil$
      which approximates the elimination clique $\cliq{\SS}{v}$.
      }
    \label{alg:CliqueTreeSampleMultiEdgeMerge}
    \begin{algorithmic}[1]
      \Procedure{CliqueTreeSampleMultiEdgeMerge}{$l,v,\SS,G$}
      \State $\CCtil \gets \matzero \in \R^{n \times n}$
      \State Let $N(G,v)$ be the set of vertices that neighbor $v$ in $G$ \label{lne:nbrsdef2}
      \State Let $E(G,v,i)$ be the set of multi-edges between $v$ and $i$ in $G$
      \State Let $w(G,v,i) = \sum_{e \in E(G,v,i)} w(e)$ be the sum of weights of multi-edges between $v$ and $i$ in $G$
      \State Let $w(G,v) = \sum_{i \in N(G,v)} w(G,v,i)$  be the sum of weights of multi-edges incident to $v$ in $G$ \label{lne:wdef2}
      \State \Comment{The quantities above are updated as $G$ changes}
      \State $d \gets w(G,v)$ it \emph{initial} weight of multi-edges incident to $v$
      \State \Comment{$d$ is \emph{not} updated as $G$ changes}
      \ForAll{vertices $i \in N(G,v)$}      \Comment{pick any ordering on
      the neighbors}
      \State $t \gets \min(\abs{E(G,v,i)},l)$ \label{lne:multEdgeCopiesMerge}
      \State $\bar{w}_{vi} \gets w(G,v,i)/t$
      \Comment{Average of the multi-edge weights between $v$ and $i$ in $G$.}
      \State Delete multi-edges $E(G,v,i)$ from $G$
      \State 
      \Comment{Update values defined in Lines~\ref{lne:nbrsdef2}-\ref{lne:wdef2} correspondingly}
      \State $w_{\text{new}} \gets \bar{w}_{vi} \cdot \frac{w(G,v)}{d}$
      \Comment{Weight used for samples}
      \For{$h = 1$ to $t$}
      \State Sample index $j$ with probability $p(j) = \frac{w(G,v,j)}{w(G,v)}$
      \State $\CCtil \gets \CCtil + w_{\text{new}} \cdot (\ee_i - \ee_j) (\ee_i - \ee_j)^{\trp}$
      \State Add a multi-edge between $i$ and $j$ of weight $w_{\text{new}}$ in to $G$
      \EndFor 
      \EndFor
      \State \Return $\CCtil$
      \EndProcedure
    \end{algorithmic}
\end{algorithm}
Note that Algorithm~\ref{alg:CliqueTreeSampleMultiEdgeMerge} differs only from Algorithm~\ref{alg:CliqueTreeSampleMultiEdge} by the choice of $t \gets \min(\abs{E(G,v,i)},l)$ in Line~\ref{lne:multEdgeCopiesMerge}.

Again, the procedure is correct in expectation.
  \begin{claim}
    \label{cor:expctcliquesamplemultiedgelimit}
      \[
      \expct{}{\textsc{CliqueTreeSampleMultiEdgeMerge}(l, v,\SS,G)} = \cliq{\SS}{v}
    \]
\end{claim}

Again the proof is simple and we omit it.

This variant has the advantage of producing a relatively more robust approximate Cholesky factorization, 
but at the same time not creating too many non-zeroes in the output factorization. Again, with a higher 
number of initial splitting and limit for merging, this variant produces a more reliable factorization, 
but also has a longer runtime.

\input{rowform}

\subsection{Our Full Algorithms}
We provide two SDDM solver implementations based on the algorithms described in the previous sections. 
Both solvers convert an $n \times n$ SDDM linear equation to an $(n+1) \times (n+1)$ Laplacian linear equation
using the standard approach of adding an additional row and column with entries corresponding to the diagonal excess.
We then solve this linear equation using PCG with a preconditioner given by an approximate Cholesky factorization of the Laplacian.

Our first implementation, AC, uses the $\textsc{CliqueTreeSample}$ sampling rule (Algorithm~\ref{alg:CliqueTreeSample}) and formats the output as row operations (Algorithm~\ref{alg:ApproximateEdgewiseCholesky}).
Our second implementation, AC-s$x$m$y$, uses the \textsc{Clique\-Tree\-Sample\-MultiEdge\-Merge} sampling rule (Algorithm~\ref{alg:CliqueTreeSampleMultiEdgeMerge}) 
with split parameter $x$ and merge parameter $y$ and formats the output as row operations (Algorithm~\ref{alg:ApproximateEdgewiseCholesky}).

Because we work with Laplacian matrices (which have kernel of dimension one or higher), we need to account for the kernel when using an (approximate) Cholesky decomposition as a preconditioner.
This is standard, but for completeness we describe this step in Appendix~\ref{sec:cholBackground}.

%% file: rowform.tex
\subsection{A Row-Operation Format for Cholesky Factorization Output}
\label{sec:rowform}
In this section, we describe an alternative approach to representing
a Cholesky factorization. In this form, it becomes a product of row operation matrices.
The two representations, either as a single lower triangular
or as a product of row operation matrices is, still result in the
exact same matrix, and the representations only take constant factor
difference in space consumption.
Our implemented code uses this alternate row operation representation, which we have found to be more efficient as it allows certain optimizations of space usage.
However, we believe other approaches can make the standard Cholesky factorization comparably efficient.

We can use the row-operations form in conjunction with the
clique sampling introduced in the previous section, and again we
obtain a fast algorithm for approximate Cholesky factorization.
The choice of row-operations form or lower-triangular form has no
impact on the sampling algorithm and does not change the output of the
algorithm, it only computes a different representation of the output\footnote{This is true in the RealRAM model, but the two
  forms may have different numerical stability properties in finite
  precision arithmetic.}. 

It is well-known that Cholesky factorization
can be interpreted as a sequence of row and column operations.
We can use this interpretation to write the lower-triangular matrix
$\matlow$ of a Cholesky factorization $\matlow \matlow^{\trp}$ as a
product of matrices that perform a row operation.

A less appreciated fact, however, is that we have some flexibility
when choosing the scaling of the row operation.
We introduce a particular choice of scaling with an interesting property:
When we apply a row operation to a Laplacian matrix, as an intermediate step in
eliminating column $v$, we then remove the entry of row $i$
corresponding to the edge from $v$ to $i$ and in the process we create
new entries that correspond exactly to the
\emph{elimination star} defined in
Equation~\eqref{eq:schurfromsingleedge}.

This motivates our decomposition of the elimination clique from
Equation~\eqref{eq:Cv-def} into elimination stars.

We assume we are dealing with a connected
graph, and describe a row-operation representation of partial Cholesky
factorization when it is used to eliminate a single vertex.
We can eliminate multiple vertices by repeated application.

We let
\[
\LL = 
\left(
\begin{array}{ccc}
d & -w  & -\aa^\trp \\
  -w & \multicolumn{2}{l}{
       \multirow{2}{*}{
             \hspace{-0.9em}
      $\begin{pmatrix}
         w &  \veczero^\trp \\
        \veczero & \diag(\aa)
       \end{pmatrix} + \LL_{-1}$
        }} \\ %
  -\aa & 
\end{array} \right)
\]
\paragraph{Degree-1 elimination.}
When eliminating a vertex of degree 1, i.e. when $\aa = \veczero$,
we use the following factorization.
\begin{align*}
\left(
\begin{array}{ccc}
w & 0 & \veczero^\trp \\
 0 & \multicolumn{2}{c}{
     \multirow{2}{*}{
      $\LL_{-1}$
        }} \\ %
  \veczero & 
\end{array}
             \right)
=
\left(
\begin{array}{ccc}
  1 & 0 &\veczero^{\trp} \\
  1  & \multicolumn{2}{c}{
       \multirow{2}{*}{
      $\II$
        }} \\
   \veczero
\end{array}
\right)
\left(
\begin{array}{ccc}
w & -w  & \veczero^\trp \\
  -w & \multicolumn{2}{l}{
       \multirow{2}{*}{
             \hspace{-0.9em}
      $\begin{pmatrix}
         w &  \veczero^\trp \\
        \veczero & \matzero^{\phantom{\trp}}
       \end{pmatrix} + \LL_{-1}$
        }} \\ %
  \veczero & 
\end{array} \right)
\left(
\begin{array}{ccc}
  1 & 1 &\veczero^{\trp} \\
  0  & \multicolumn{2}{c}{
       \multirow{2}{*}{
      $\II$
        }} \\
   \veczero
\end{array}
\right)
\end{align*}
And hence, by applying inverses of the outer matrices in the
factorization,
\begin{align*}
  \left(
\begin{array}{ccc}
w & -w  & \veczero^\trp \\
  -w & \multicolumn{2}{l}{
       \multirow{2}{*}{
             \hspace{-0.9em}
      $\begin{pmatrix}
         w &  \veczero^\trp \\
        \veczero & \matzero^{\phantom{\trp}}
       \end{pmatrix} + \LL_{-1}$
        }} \\ %
  \veczero & 
\end{array} \right)
=
\left(
\begin{array}{ccc}
  1 & 0 &\veczero^{\trp} \\
  -1  & \multicolumn{2}{c}{
       \multirow{2}{*}{
      $\II$
        }} \\
   \veczero
\end{array}
\right)
\left(
\begin{array}{ccc}
w & 0 & \veczero^\trp \\
 0 & \multicolumn{2}{c}{
     \multirow{2}{*}{
      $\LL_{-1}$
        }} \\ %
  \veczero & 
\end{array}
\right)
\left(
\begin{array}{ccc}
  1 & -1 &\veczero^{\trp} \\
  0  & \multicolumn{2}{c}{
       \multirow{2}{*}{
      $\II$
        }} \\
   \veczero
\end{array}
\right)
\end{align*}
\paragraph{Edge elimination.}
When the vertex we're in the process of eliminating has degree more
than 1, we instead apply the following factorization,
where $d = \vecone^{\trp} \aa + w$ and $\theta = w/d$.
\begin{align}
\label{eq:elimsingleedge}
  &\left(
\begin{array}{ccc}
d (1-\theta)^2 & 0  & -(1-\theta) \aa^\trp \\
  0 & \multicolumn{2}{l}{
       \multirow{2}{*}{
             \hspace{-0.9em}
      $\begin{pmatrix}
         w - w^2/d &  -\frac{w}{d} \aa^\trp \\
        -\frac{w}{d} \aa & \diag(\aa)
       \end{pmatrix} + \LL_{-1}$
        }} \\ %
 -(1-\theta) \aa& 
\end{array} \right)
\\ \nonumber
&=
\left(
\begin{array}{ccc}
  1-\theta & 0 &\veczero^{\trp} \\
  \theta   & \multicolumn{2}{c}{
       \multirow{2}{*}{
      $\II$
        }} \\
   \veczero
\end{array}
\right)
\left(
\begin{array}{ccc}
d & -w  & -\aa^\trp \\
  -w & \multicolumn{2}{l}{
       \multirow{2}{*}{
             \hspace{-0.9em}
      $\begin{pmatrix}
         w &  \veczero^\trp \\
        \veczero & \diag(\aa)
       \end{pmatrix} + \LL_{-1}$
        }} \\ %
  -\aa & 
\end{array} \right)
\left(
\begin{array}{ccc}
  1-\theta & \theta &\veczero^{\trp} \\
  0  & \multicolumn{2}{c}{
       \multirow{2}{*}{
      $\II$
        }} \\
   \veczero
\end{array}
  \right)
\end{align}
Note that the matrix on the LHS is a Laplacian because $d
(1-\theta)^2=  (1-\theta) \vecone^\trp\aa$.
And again, by applying inverses of the outer matrices in the
factorization, we have
\begin{align*}
&
\left(
\begin{array}{ccc}
d & -w  & -\aa^\trp \\
  -w & \multicolumn{2}{l}{
       \multirow{2}{*}{
             \hspace{-0.9em}
      $\begin{pmatrix}
         w &  \veczero^\trp \\
        \veczero & \diag(\aa)
       \end{pmatrix} + \LL_{-1}$
        }} \\ %
  -\aa & 
\end{array} \right)
\\
&=
\left(
\begin{array}{ccc}
  \frac{1}{1-\theta} & 0 &\veczero^{\trp} \\
  \frac{-\theta}{1-\theta}  & \multicolumn{2}{c}{
       \multirow{2}{*}{
      $\II$
        }} \\
   \veczero
\end{array}
\right)
\left(
\begin{array}{ccc}
d (1-\theta)^2 & 0  & -(1-\theta) \aa^\trp \\
  0 & \multicolumn{2}{l}{
       \multirow{2}{*}{
             \hspace{-0.9em}
      $\begin{pmatrix}
         w - w^2/d &  -\frac{w}{d} \aa^\trp \\
        -\frac{w}{d} \aa & \diag(\aa)
       \end{pmatrix} + \LL_{-1}$
        }} \\ %
 -(1-\theta) \aa& 
\end{array} \right)
\left(
\begin{array}{ccc}
  \frac{1}{1-\theta} & \frac{-\theta}{1-\theta} &\veczero^{\trp} \\
  0  & \multicolumn{2}{c}{
       \multirow{2}{*}{
      $\II$
        }} \\
   \veczero
\end{array}
  \right)
\end{align*}
\paragraph{Schur complement invariance.}
We can also see that the Schur complement onto the remaining vertices
is 
\begin{align*}
\SS = 
  \schurto{
\left(
\begin{array}{ccc}
d & -w  & -\aa^\trp \\
  -w & \multicolumn{2}{l}{
       \multirow{2}{*}{
             \hspace{-0.9em}
      $\begin{pmatrix}
         w &  \veczero^\trp \\
        \veczero & \diag(\aa)
       \end{pmatrix} + \LL_{-1}$
        }} \\ %
  -\aa & 
\end{array}
\right)
}{[n]\setminus\setof{1}}
&=
\begin{pmatrix}
         w - w^2/d &  -\frac{w}{d} \aa^\trp \\
        -\frac{w}{d} \aa & \diag(\aa) - \frac{1}{d}\aa \aa^{\trp}
\end{pmatrix} + \LL_{-1}
  \\
  &=
\schurto{
\left(
\begin{array}{ccc}
d (1-\theta)^2 & 0  & -(1-\theta) \aa^\trp \\
  0 & \multicolumn{2}{l}{
       \multirow{2}{*}{
             \hspace{-0.9em}
      $\begin{pmatrix}
         w - w^2/d &  -\frac{w}{d} \aa^\trp \\
        -\frac{w}{d} \aa & \diag(\aa) 
       \end{pmatrix} + \LL_{-1}$
        }} \\ %
 -(1-\theta) \aa& 
\end{array} \right)
}{[n]\setminus\setof{1}}
\end{align*}

\paragraph{Eliminating a vertex, one edge at a time.}
Let us summarize these observations into a statement about how to
write a factorization of $\LL$ that eliminates the first vertex,
\emph{which we will now denote by vertex 0}.
Assume vertex 0 has degree $k$, and that its neighbors are vertices
$1,2,\ldots,k$.
\begin{equation}
  \label{eq:lapuniformlayout}
  \LL = 
\left(
\begin{array}{ccc}
d & -\aa^\trp \\
-\aa& \diag(\aa) + \LL_{-1}
\end{array} \right)
\end{equation}

We denote the weight on the edges from vertex 0 to its neighbors by
$\aa(1), \aa(2), \ldots, \aa(k)$.
We then write
\begin{equation}
  \label{eq:edgeelimfactor}
  \LL
  =
  \matlow_1 \matlow_2 \ldots \matlow_k
\left(
\begin{array}{ccc}
\phi & \veczero^\trp \\
  \veczero &  \SS
\end{array}
\right)
\matlow_k^{\trp}\ldots  \matlow_2^{\trp} \matlow_1^{\trp}
\end{equation}
where $\SS = \schurto{\LL}{[n]\setminus\setof{1}}$, and 
where for $i < k$, we have
\[
  \matlow_i
  =
  \left(
  \begin{array}{ccc}
  \frac{1}{1-\theta_i}  & \veczero^{\trp} \\
    \frac{-\theta_i}{1-\theta_i}
    \ee_i
    & \II
  \end{array}
\right)
=
       \II +  \frac{\theta_i}{1-\theta_i} (\ee_1 - \ee_i) \ee_1^{\trp}
\]
where $\ee_i$ is the $i$ basis vector in dimension $n -1$,
and
$
  \theta_i= \frac{\aa(i) \Pi_{j < i} (1-\theta_i) }{d \cdot \Pi_{j < i} (1-\theta_i)^2} = \frac{\aa(i)}{d \cdot \Pi_{j < i} (1-\theta_i)}
  $.
We can simplify this, using the observation that 
$\Pi_{j < i} (1-\theta_j) =
\frac{d - \sum_{j<i} \aa(j)}{d}$.
Hence
\[
  \theta_i=  \frac{\aa(i)}{d - \sum_{j<i} \aa(j)}
.
  \]

The last factor is given by,
  \[
  \matlow_k
  =
  \left(
  \begin{array}{ccc}
  1 & \veczero^{\trp} \\
  - \ee_k
    & \II
  \end{array}
\right)
=
       \II -  \ee_i \ee_1^{\trp}
\]
and $\phi = \aa(k) \Pi_{j < k} (1-\theta_j) = \frac{\aa(k)^2}{d}$.

\paragraph{Computing an edgewise Cholesky factorization.}
In this section, we briefly remark how to repeatedly eliminate vertices to obtain a full
edgewise Cholesky
factorization.

Let us slightly modify the notation for the first elimination to write
\begin{equation*}
  \LL
  =
  \matlow_1^{(1)} \matlow_2^{(1)} \cdots \matlow_{k_1}^{(1)}
\left(
\begin{array}{ccc}
\phi_1 & \veczero^\trp \\
  \veczero &  \SS
\end{array}
\right)
( \matlow_{k_1}^{(1)})^{\trp}\cdots  (\matlow_2^{(1)})^{\trp}
(\matlow_1^{(1)})^{\trp}
,
\end{equation*}
where $k_1$ denotes the degree of the vertex we eliminated,
and $\SS = \schurto{\LL}{[n]\setminus\setof{1}}$.

Now if we recursively factor the remaining matrix
$\SS$
we can eventually write

\begin{equation}
  \label{eq:fullfactorfactorperedge}
  \LL
  =
 \left(
   \Pi_{i = 1}^n \Pi_{j \sim i }
    \matlow_j^{(i)}
\right)
\PPhi
  \left(
   \Pi_{i = 1}^n \Pi_{j \sim i }
    \matlow_j^{(i)}
\right)^{\trp}
\end{equation}
Note that $j \sim i$ denotes the set of $j$ that are neighbors of $i$
in the Schur complement that $i$ is being eliminated from.

$\PPhi$ is a diagonal matrix with $\PPhi(i,i)$ being the
diagonal entry that results from the final single-edge elimination of vertex $i$.

Each row operation matrix $\matlow_j^{(i)}$ only has two where it
differs from the identity matrix and hence $\matlow_j^{(i)} x$ can be
computed from $x$ in $O(1)$ time by only modifying a constant number
of entries of $x$.
Similarly, the inverse of these row operation matrices can be applied
in $O(1)$.

We summarize edgewise Cholesky
factorization in the pseudo-code below. The output is a diagonal
matrix $\PPhi$ and a sequence of row-operation matrices
$\setof{\matlow^{(v)}_i}_{v,i}$, which gives an edgewise Cholesky
factorization of a Laplacian in the sense of
Equation~\eqref{eq:fullfactorfactorperedge}.

\begin{algorithm}
  \caption{Algorithm
    \textsc{EdgewiseCholesky}{($\LL$)} outputs an
    edgewise Cholesky
    factorization of the input Laplacian
    ${\LL \in \R^{n \times n}}$.}
  \label{alg:EdgewiseCholesky}
  \begin{algorithmic}[1]
    \Procedure{EdgewiseCholesky}{$\LL$}
    \State $\PPhi \gets \matzero$
    \State $\SS \gets \LL$
     \For{$v=1$ to $n-1$} \Comment Eliminate any $n-1$ vertices in any
     order.
    \State $d \gets \SS(v,v)$; $\aa \gets -\SS(v,:)$; $\aa(v) \gets 0$
     \For{$i  \in \setof{ \text{nbrs. of } v}$, excluding one nbr. denoted by
       index $k_v$}
     \State
     $\theta_i \gets  \frac{\aa(i)}{\vecone^{\trp} \aa}$
          \Comment The neighbor eliminations can
     \State
          $
       \matlow^{(v)}_i
       \gets
       \II +  \frac{\theta_i}{1-\theta_i} (\ee_v - \ee_i) \ee_v^{\trp}
       $
            \Comment be performed in any
     order.
       \State $\aa(i) \gets 0$
       \EndFor
     \State  $\matlow^{(v)}_{k_v}
     \gets     
     \II + \ee_{k_v}\ee_v^{\trp}
       $
       \State $\PPhi(v,v) \gets \frac{\aa(k_v)^2}{d}$
    \State $\SS \gets \SS  - \vstar{\SS}{v} + \cliq{\SS}{v}$
    \EndFor 
    \State \Return
   $\PPhi, \setof{\matlow^{(v)}_i}_{v,i}$
    \EndProcedure
  \end{algorithmic}
\end{algorithm}

\paragraph{Obtaining an algorithm for edgewise approximate Cholesky factorization.}
In Section~\ref{sec:apxChol},
we saw a meta-algorithm $\textsc{ApproximateCholesky}$
(Algorithm~\ref{alg:GeneralApproximateCholeskyFactorizationColwise}),
which shows how we can replace the elimination clique in standard to
Cholesky factorization
($\textsc{Cholesky}$, Algorithm~\ref{alg:LapCholeskyFactorization})
with a sampled clique approximation
to obtain an approximate Cholesky factorization.
And we introduced a new approach to clique sampling,
$\textsc{CliqueTreeSample}$ (Algorithm~\ref{alg:CliqueTreeSample}),
which we can plug into the $\textsc{ApproximateCholesky}$ meta-algorithm.

Similarly, if we replace the elimination clique in the edgewise Cholesky
factorization of the previous section, then we immediately get an
approximate edgewise Cholesky factorization.

This could be framed as a meta-algorithm with an arbitrary
clique sampling routine, but in the interest
of brevity, we state our pseudo-code directly using
$\textsc{CliqueTreeSample}$.
The algorithm appears below as Algorithm~\ref{alg:ApproximateEdgewiseCholesky}.

\begin{algorithm}
  \caption{Algorithm
    \textsc{ApproximateEdgewiseCholesky}{($\LL$)} outputs an
    approximate edgewise Cholesky
    factorization $\PPhi, \setof{\matlow^{(v)}_i}_{v,i}$ such that
    $
      \left(
   \Pi_{i = 1}^n \Pi_{j \sim i }
    \matlow_j^{(i)}
\right)
\PPhi
  \left(
   \Pi_{i = 1}^n \Pi_{j \sim i }
    \matlow_j^{(i)}
\right)^{\trp}
\approx \LL$.
  }
  \label{alg:ApproximateEdgewiseCholesky}
  \begin{algorithmic}[1]
    \Procedure{ApproximateEdgewiseCholesky}{$\LL$}
    \State $\PPhi \gets \matzero$
    \State $\SS \gets \LL$
     \For{$v=1$ to $n-1$} \Comment Eliminate any $n-1$ vertices in any
     order.
    \State $d \gets \SS(v,v)$; $\aa \gets -\SS(v,:)$; $\aa(v) \gets 0$
     \For{$i  \in \setof{ \text{nbrs. of } v}$, excluding one nbr. denoted by
       index $k_v$}
     \State
     $\theta_i \gets  \frac{\aa(i)}{\vecone^{\trp} \aa}$
          \Comment The neighbor eliminations can
     \State
          $
       \matlow^{(v)}_i
       \gets
       \II +  \frac{\theta_i}{1-\theta_i} (\ee_v - \ee_i) \ee_v^{\trp}
       $
            \Comment be performed in any
     order.
       \State $\aa(i) \gets 0$
       \EndFor
     \State  $\matlow^{(v)}_{k_v}
     \gets     
     \II + \ee_{k_v}\ee_v^{\trp}
       $
       \State $\PPhi(v,v) \gets \frac{\aa(k_v)^2}{d}$
    \State $\SS \gets \SS  - \vstar{\SS}{v} +
    \textsc{CliqueTreeSample}(v,\SS)$
    \EndFor 
    \State \Return
   $\PPhi, \setof{\matlow^{(v)}_i}_{v,i}$
    \EndProcedure
  \end{algorithmic}
\end{algorithm}

The outputs of
$\textsc{ApproximateEdgewiseCholesky}$
and $\textsc{ApproximateCholesky}$ are identical when viewed as linear
operators.
We this claim formally below.

\begin{claim}
  \label{clm:equivalence}
  \quad \\
Suppose we run both $\textsc{ApproximateCholesky}$ and 
$\textsc{ApproximateEdgewiseCholesky}$
\begin{itemize}
\item using the same elimination ordering,
\item and using \textsc{CliqueTreeSample} as the clique
  sampling routine, and using the same outputs
  of \textsc{Clique\-Tree\-Sample}, i.e. the same outcomes of the random
  samples.
\end{itemize}
  Then the output factorizations $\matlow$ from
  $\textsc{ApproximateCholesky}$ and $\PPhi,
  \setof{\matlow^{(v)}_i}_{v,i}$ from
  \textsc{Approximate\-Edgewise\-Cholesky} will satisfy
  \[
    \matlow \matlow^{\trp} =
    \left(
   \Pi_{i = 1}^n \Pi_{j \sim i }
    \matlow_j^{(i)}
\right)
\PPhi
  \left(
   \Pi_{i = 1}^n \Pi_{j \sim i }
    \matlow_j^{(i)}
  \right)^{\trp}
\]
In other words, the only difference between the two algorithms is in
the formatting of the output.
\end{claim}

We prove the claim in Appendix~\ref{sec:equivalence}.

\begin{remark}
  The claim above shows the algorithms differ only in the formatting
  of their outputs.
  In fact, we can also convert the output of either to the output of the
  other in linear time.
\end{remark}

From this, we can also deduce that
\[
\expct{}{
\left(
   \Pi_{i = 1}^n \Pi_{j \sim i }
    \matlow_j^{(i)}
\right)
\PPhi
  \left(
   \Pi_{i = 1}^n \Pi_{j \sim i }
    \matlow_j^{(i)}
\right)^{\trp}} = \LL
\]

Furthermore, we can see by inspection that Claim~\ref{clm:fillin} essentially applies: Each row operation $\matlow_j^{(i)}$ consists of an identity matrix with two entries adjusted, one on the diagonal, and one below the diagonal.
This means the matrix and its inverse can be applied in $O(1)$ time, and the overall sequence of row operations (or their inverses) can be applied in $O(m\log m)$ time when the input matrix has $O(m)$ non-zero entries.

%% file: results.tex
\section{Experimental Evaluation}
\label{sec:experiments}

\subsection{Implementation}
\label{sec:impl}
We provide two main implementations of our SDDM and Laplacian linear
equation solvers, which we refer to as {\AC} and AC-s$x$m$y$ (AC with
$x$ edge splits and merge threshold  $y$).
For our evaluation of AC-s$x$m$y$ we focus on the parameter setting
``split=2, merge=2'', which we refer to as AC-s2m2, or simply {\ACSM}.
{\AC} produces the same factorization as  AC-s1m1 would, but uses a
slightly more efficient implementation, specific to this parameter setting.

Both {\AC} and {\ACSM} use preconditioned conjugate gradient (PCG) to compute
a solution to the input system of linear equations, up to the
specified residual error tolerance.
Our implementation of PCG is based on the book \emph{Templates for the solution of linear systems: building blocks for iterative methods}~\cite{BBCDDDEPRVdV94}.
Both solvers are initially designed for Laplacian matrices, and use a standard reduction to convert linear equations in SDDM matrices into linear equations in Laplacians.

In Sections~\ref{sec:edgewise}-\ref{sec:rowform}, we describe the overall pseudo-code of
these two algorithms.
\paragraph{Elimination order.} We used a non-standard elimination ordering
in our implementations, which, however, is close in spirit to
approximate minimum degree heuristics.
This ordering is greedy on the unweighted degree of the vertices, in other words,
the vertex being eliminated is the vertex with the (approximately) least unweighted degree. We implemented 
an approximate priority queue to support dynamically updating the unweighted degrees of vertices during 
the elimination.

\paragraph{Tuning the sampling quality of  AC-s$x$m$y$: {\AC} vs. {\ACSM} and more.} 
In addition to our basic version of the approximate Cholesky
factorization implemented in {\AC}, 
we also implemented the two variants described in
Section~\ref{sec:variants}.
In our experimental evaluations, we will
refer to the basic version as {\AC}, the version that splits the edges into $x$ copies initially with no merging as 
AC-s$x$, and the version that splits the edges into $x$ copies initially and merges multi-edges up $y$ as 
AC-s$x$m$y$.
For users, we recommend  AC-s2m2, which we refer to as {\ACSM}. Our experiments suggest
that AC-s2m2 is a reliable, robust, and still fast choice.
As we have briefly explained in Section~\ref{sec:variants}, one
expects that AC-s$x$ produces the most reliable 
factorization while AC is the fastest.
AC-s$x$m$y$ is faster than AC-s$x$ at computing a factorization, often much faster, while still having significantly
improved reliability compared to AC.
In our experiments, we find that AC and AC-s2m2 provide the best
performance.
For ``easy'' problems, the factorizations given by AC are reliable
enough and hence has the best runtime.
However, for ``harder'' problems, such as the Sachdeva star, AC cannot
produce a reliable factorization and the PCG algorithm might run into stagnation. In contrast to that, AC-s2m2 remains to
produce reliable factorizations for Sachdeva stars, while maintaining a relatively good runtime. Therefore, AC should be used
as the main version for ``easier'' problems while AC-s2m2 for
``harder'' problems. In Section~\ref{sec:variantseval}, we provide a detailed experimental comparison between AC, AC-s1, AC-s2, AC-s2m2, and AC-s3m3. In other experiments,
we will focus on AC and AC-s2m2 ({\ACSM}).

\paragraph{Output format and forward/backward substitution.}
Both our implementations {\AC} and {\ACSM} use a non-standard row operation form of the algorithm and the output factorization is in a customized format. 
This customized format can be transformed into a standard Cholesky lower-triangular format. 
We implemented a routine for applying the inverse of the output factorization without 
converting to the standard Cholesky lower-triangular format.
Our non-standard format for the output factorization is based on the
edge-wise elimination format described in Section~\ref{sec:rowform}.
We believe a standard Cholesky factorization format can be
equally efficient and potentially easier to parallelize.
In the case of single-threaded performance, we believe the distinct between the edgewise-elimination
format and a standard Cholesky factorization is not crucial.

\subsection{Solver Comparisons: Summary of Experiments}

\paragraph{Experiment set-up.}
All benchmarking was ran on a $24$ Core Intel Xeon Gold $6240$R $2.4$GHz Processor with $512$ GB memory, on Ubuntu.
Experiments were run with single-threaded versions of each solver.

\paragraph{Evaluated linear equation solvers.}
We compare our solver to
\begin{itemize}
\item PETSc implementation of PCG with the PETSc version of the
  HyPre-BoomerAMG as the preconditioner. We run PETSc with a single processor.
\item HyPre implementation of PCG with HyPre-BoomerAMG
  as the solver.
\item Combinatorial Multigrid (CMG)'s implementation of PCG with its combinatorial preconditioner.
\item MATLAB implementation of PCG with MATLAB's ichol Incomplete Cholesky factorization (ICC) as the preconditioner.
\item Lean Algebraic Multigrid (LAMG)'s Laplacian linear
  solver\footnote{We omit LAMG from our result tables, as we
    found it unable to convergence to our target tolerance across all
    matrix families we tested.}
\end{itemize}

\paragraph{Evaluation statistics.}
\begin{itemize}
\item To generate a right hand side vector for a linear equation in a matrix $\MM$, we apply $\MM$ to a random Gaussian $\gg$ and normalize, yielding $\MM\gg/\norm{\MM\gg}_2$ as the right hand side.
This ensures the right hand side is in the image of $\MM$. 
\item For each linear equation, we report the time and iterations required to obtain
  a relative residual of $10^{-8}$.
  If the solver returns a result with higher residual error, we mark the time entry with $*$ if the residual error exceeds the target by a factor $(1,10^{4}]$ and by $**$ if it exceeds the target by a factor $(10^{4},10^{8})$, and finally we  report the time as $\infty$ (INF) if the residual error is too large by a factor $10^{8}$ or more.
  Note that the zero vector obtains a residual error matching this trivial bound.
\item We report wall-clock time.
\item We exclude the time required to load the input matrix into
  memory before starting the solver (because our file loading system
  for HyPre is slow).
\end{itemize}

For each solver, for each system of linear equations, we report the following statistics, or a subset thereof. 
When we summarize across many instances, we report the median, 75th percentile, and worst-case running times seen across these instances.
\begin{itemize}
\item $n$: The linear equation is in $n$ variables.
\item $\nnz$: the number of non-zeros in the original matrix.
\item$t_{\text{total}}$: the total time in seconds (wall-clock time) to build our approximate Cholesky
  factorization and solve a linear equation in the input matrix.
  $t_{\text{total}} = t_{\text{build}} + t_{\text{solve}}$.
\item $t_{\text{build}}$: the total time in seconds (wall-clock time) to build our approximate Cholesky
  factorization of the input matrix.
\item$t_{\text{solve}}$: the total time  in seconds (wall-clock time) to solve a linear equation in the
  input matrix.
\item $N_{\text{iter}}$: The number of iterations of PCG required to
  reach the desired relative residual.
\end{itemize}

\paragraph{ Linear equation systems used for experimental evaluation.}
We present evaluations of our solver on a host of different examples, falling into three categories:
\begin{enumerate}
\item Matrices from the SuiteSparse Matrix Collection which are SDDM or approximately SDDM.
\item Programmatically generated SDDM matrices.
\begin{enumerate}
\item ``Chimera'' graphs Laplacians and variants that are strictly SDDM.
\item ``Sachdeva-star'' graphs Laplacians.
\item Laplacian matrices arising from solving max flow problems using an interior point method.
\end{enumerate}
    \item Discretizations of partial differential equations on 3D grids:
    \begin{enumerate}
        \item Matrices based on fluid simulations in a Society of Petroleum Engineering (SPE) benchmark.
        \item Uniform coefficient Poisson problems on 3D grids.
\item Variable coefficient Poisson problems on 3D grids with a
  checkerboard board pattern using variable resolution and fixed weight.
\item Anisotropic coefficient Poisson problems on 3D grids with variable discretization and fixed weight.
\item Anisotropic coefficient Poisson problems on 3D grids with fixed resolution
  and variable weight.
    \end{enumerate}
\end{enumerate}
We describe our experiments and timing results on the classes listed above in detail in Section~\ref{sec:experimentdescription}.
In Section~\ref{sec:solvervariance} we run a few additional tests where we report the variation in solve time for our solvers. These tests show that the variance is solve time is very low.
In Section~\ref{sec:variantseval}, we run experiments that suggest why {\ACSM} sometimes outperforms {\AC}. These experiments show that on difficult instances, we show that the preconditioner computed by {\ACSM} has much better relative condition number than that computed by {\AC}.

\paragraph{Results and conclusions.}
Table~\ref{tab:summarytable} provides an overview the worst case
performance of all the tested solvers across all the tests we ran.
The rest of our results are shown in
Tables~\ref{tab:uniformgridtable}-\ref{tab:variants}.
For a discussion and overview of the experiments, we refer the reader to Section~\ref{sec:introexperiments}.

\begin{table}[H]
     \begin{center}
  \input{summarytable-times}
    \end{center}
  \caption{Summary of worst case total solve time per non-zero matrix
    entry across all experiments. }
  \label{tab:summarytable}
\end{table}

\subsection{Solver Comparisons: Description of Experiments and Results}
\label{sec:experimentdescription}

\input{suitesparse}

\input{chimeraAndIPM}

\input{poisson}

\subsection{Variation in the Performance of our Solvers}
\label{sec:solvervariance}
In this section, we report the statistics on the variation in the
performance of our solvers AC and {\ACSM}) (i.e. AC-s2m2). 
We evaluate the variation in the performance of these solvers with a
weighted Chimera with 10 million variables, a weighted SDDM Chimera with 
10 million variables, the Sachdeva Star with parameter $k=700$, high contrast coefficient Poisson grid with $200$ million non-zeros and axis interval $32$, and 
Anisotropic coefficient Poisson grid with fixed discretization and weight $0.001$.
For each system of linear equations, we run each solver 10 times.

Results are shown in Table~\ref{tab:variance}. Experiments showed that our solvers' solve time and total time is tightly concentrated, with
the only exception occurring with AC on Sachdeva star. This is because AC failed to produce reliable factorizations for
Sachdeva stars. In comparison, AC-s2m2 still manage to exhibit a very small variance of runtime on Sachdeva star.

\begin{threeparttable}[H]
  \centering
\resizebox{\textwidth}{!}{%
\begin{tabular}{
  |c
  |S[table-format=3.1]
  |S[table-format=1.2] S[table-format=1.2] S[table-format=1.2]
  |S[table-format=1.3] S[table-format=1.2] S[table-format=1.3]
  |S[table-format=1.2] S[table-format=1.2] S[table-format=1.2]
  |S[table-format=1.3] S[table-format=1.3] S[table-format=1.3]
  |
  }
\hline
  {Instance}
  & {non-zeros}
  & \multicolumn{3}{|c|}{AC { $t_{\text{solve}} / \nnz$ }}
  & \multicolumn{3}{c|}{{\ACSM} (ac-s2m2) { $t_{\text{solve}} / \nnz$ }}
  & \multicolumn{3}{c|}{AC { $t_{\text{total}} / \nnz$ }} 
  & \multicolumn{3}{c|}{{\ACSM} (ac-s2m2)  { $t_{\text{total}} / \nnz$ }}\\ \hline
& {nnz} &{ median   }&{ 0.75    }&{ max   }&{ median   }&{ 0.75    }&{ max   }&{ median  }&{ 0.75   }& {max} &{ median  }&{ 0.75   }& {max} \\ \hline
& {$(M)$}&{$\mu s$}  & {$\mu s$} & {$\mu s$} & {$\mu s$} &{$\mu s$}  & {$\mu s$} & {$\mu s$} &{$\mu s$}  & {$\mu s$} & {$\mu s$} &{$\mu s$}  & {$\mu s$} \\ \hline
Weighted chimera&50.2&3.78&3.87&3.95&2.59&2.6&2.68&4.65&4.74&4.82&4.09&4.11&4.18\\
Weighted SDDM chimera&49.7&3.3&3.38&3.48&2.44&2.5&2.58&4.45&4.53&4.64&4.26&4.33&4.4\\
Sachdeva Star&172&2.56\tnote{*}&3.52\tnote{**}&4.05\tnote{**}&0.208&0.212&0.214&2.96\tnote{*}&3.86\tnote{**}&4.37\tnote{**}&0.977&0.982&0.986\\
High contrast coefficient Poisson grid&200&1.71&1.73&1.77&2.14&2.15&2.25&2.6&2.62&2.66&3.34&3.35&3.45\\
Anisotropic coef. Poisson grid, variable weight&200&1.12&1.14&1.18&1.15&1.16&1.2&1.72&1.75&1.82&2.33&2.34&2.38\\

 \hline
\end{tabular}
} %
 \caption{Variation in total time and solve time of ac and ac-s2m2 computed across 10 runs per system of linear equations.}
\label{tab:variance}
\begin{tablenotes}
\footnotesize
\item[*] relative residual error exceeded tolerance $10^{-8}$ by a factor $(1,10^{4}]$.
\item[**] relative residual error exceeded tolerance $10^{-8}$ by a factor $(10^{4},10^{8})$.
\item[] \hspace{-1em} Inf: Solver crashed or returned a solution with relative residual error 1 (trivial solution).
\item[] \hspace{-1em} N/A: Experiment omitted as the solver crashed too often with (only occurred for PETSc).
\end{tablenotes}
\end{threeparttable}

\subsection{Comparison of Algorithm Variants: Splitting and Merging}
\label{sec:variantseval}
In this section, we report some further experimental statistics on the qualities of all variants of the approximate Cholesky 
factorization on the difficult problems described earlier. The factorizations reported here are

\begin{itemize}
  \item Approximate Cholesky without splitting and with merging of multi-edges down to 1 copy (AC).
  \item Approximate Cholesky without splitting and without merging (AC-s1).
  \item Approximate Cholesky with the original edges split into $2$ and without merging (AC-s2).
  \item Approximate Cholesky with the original edges split into $2$ and merging multi-edges down to $2$ copies (AC-s2m2).
  \item Approximate Cholesky with the original edges split into $3$ and merging multi-edges down to $3$ copies (AC-s3m3).
\end{itemize}

For these factorizations, we report the following statistics:

\begin{itemize}
  \item Total time to solve each linear equation, normalized by non-zero count.
  \item Condition number of the SDDM matrix and the approximate factorizations. This captures the quality of 
  approximate factorizations.
  \item Ratio between the size of the output factorization and the input Laplacian (adjacency) matrix.
\end{itemize}

Results are shown in Table~\ref{tab:variants}.
Here we remark that AC and AC-s1 did not achieve the target $10^{-8}$ tolerance on the Sachdeva star when $k=700$.
For Sachdeva stars, most sampled edges of our solvers
occur within the complete graphs on the leaves, hence causing the ratio between the sizes of output factorizations and input Sachdeva star Laplacians to be close to $1$. We see that the initial splitting plays a very significant 
role in the quality of the output factorization. AC-s1 produced factorizations that have very large condition numbers relative to the input matrix on Sachdeva star, even though it does not compress any multi-edges created during 
the elimination. On the other hand, the condition number and iteration counts for AC-s2m2 are still very much comparable to that of AC-s2, despite that it only keeps up to 2 multi-edges.
 This suggests that allowing no limits on multi-edges does not increase quality much compared to versions that merge, and versions without merging are noticeably slower.
AC-s3m3 produced the most reliable
factorizations and beats AC-s2 in terms of the runtime, but is not as fast as AC-s2m2.

\begin{landscape}
  \newcommand{\lstablelen}{24cm}
  \begin{threeparttable}[H]
    \centering
  \resizebox{\lstablelen}{!}{%
  \begin{tabular}{
    |c
    |S[table-format=3.1]
    |S[table-format=1.2] S[table-format=3.0] S[table-format=1.2] S[table-format=5.1]
    |S[table-format=1.2] S[table-format=3.0] S[table-format=1.2] S[table-format=4.1]
    |S[table-format=1.2] S[table-format=2.0] S[table-format=1.2] S[table-format=3.1]
    |
    }
  \hline
    {Instance}
    & {non-zeros}
    & \multicolumn{4}{|c|}{AC}
    & \multicolumn{4}{c|}{AC-s1}
    & \multicolumn{4}{c|}{AC-s2} \\ \hline
  & {nnz} &{ $t_{\text{total}} / \nnz $   }&{$N_{\text{iter}}$}&{ Ratio of sizes }&{ Condition num.   }&{ $t_{\text{total}} / \nnz $   }&{$N_{\text{iter}}$}&{ Ratio of sizes}&{ Condition num.   }&{ $t_{\text{total}} / \nnz $   }&{$N_{\text{iter}}$}&{ Ratio of sizes   }& {Condition num.}  \\ \hline
  & {$(M)$} &{$\mu s$}&  &  & &{$\mu s$}&  &  & &{$\mu s$}&  &  &  \\ \hline
  Weighted chimera&50.2&4.42&55&3.56&87.0&6.66&49&4.53&75.9&9.58&35&6.67&18.5\\
Weighted SDDM chimera&49.7&4.3&51&3.29&59.8&5.63&48&3.99&75.4&7.11&31&5.81&28.1\\
Sachdeva Star&172&2.2\tnote{*}&408&1.0&9660.0&3.13\tnote{*}&189&1.0&3090.0&5.97&31&1.0&102.0\\
High contrast coefficient Poisson grid&200&2.24&57&2.78&60.2&3.82&55&3.44&66.8&6.32&44&4.86&35.8\\
Anisotropic coef. Poisson grid, variable weight&200&1.78&41&2.66&24.7&3.24&32&3.2&16.3&5.73&23&4.41&8.78\\
  \end{tabular}
  } %
  \resizebox{\lstablelen}{!}{%
  \begin{tabular}{
    |c
    |S[table-format=3.1]
    |S[table-format=1.3] S[table-format=2.0] S[table-format=1.2] S[table-format=2.1]
    |S[table-format=1.2] S[table-format=2.0] S[table-format=1.2] S[table-format=2.2]
    |
    }
  \hline
    {Instance}
    & {non-zeros}
    & \multicolumn{4}{|c|}{AC-s2m2}
    & \multicolumn{4}{c|}{AC-s3m3}\\ \hline
  & {nnz} &{ $t_{\text{total}} / \nnz $   }&{$N_{\text{iter}}$}&{ Ratio of sizes }&{ Condition num.   }&{ $t_{\text{total}} / \nnz $   }&{$N_{\text{iter}}$}&{ Ratio of sizes}&{ Condition num.   }  \\ \hline
  & {$(M)$} &{$\mu s$}&  &  & &{$(\mu s)$}&  &  &  \\ \hline
 Weighted chimera&50.2&4.07&33&5.32&29.7&4.9&26&6.73&12.7\\
Weighted SDDM chimera&49.7&3.99&30&4.8&19.0&4.49&26&6.0&14.7\\
Sachdeva Star&172&0.84&44&1.0&64.6&1.15&31&1.0&21.4\\
High contrast coefficient Poisson grid&200&3.44&44&3.98&37.3&4.4&38&4.92&30.4\\
Anisotropic coef. Poisson grid, variable weight&200&2.33&26&3.69&11.6&3.12&21&4.46&7.51\\
    \hline
  \end{tabular}
  } %
   \caption{Comparison between variants of the approximate Cholesky factorization using different multi-edge splitting and merging approaches.}
  \label{tab:variants}
  
  \begin{tablenotes}
\footnotesize
\item[*] relative residual error exceeded tolerance $10^{-8}$ by a factor $(1,10^{4}]$.
\item[**] relative residual error exceeded tolerance $10^{-8}$ by a factor $(10^{4},10^{8})$.
\item[] \hspace{-1em} Inf: Solver crashed or returned a solution with relative residual error 1 (trivial solution).
\item[] \hspace{-1em} N/A: Experiment omitted as the solver crashed too often with (only occurred for PETSc).
\end{tablenotes}
  \end{threeparttable}
\end{landscape}

\subsection{Comparison of Algorithm Variants: Elimination Order}
\label{sec:ordervariantseval}

In this section, we investigate the impact of our approximate greedy elimination ordering (see Section~\ref{sec:impl}) compared with a randomized ordering.
The use of a completely random ordering was suggested in \cite{KS16} which introduced randomized approximate Gaussian elimination.
In \cite{KS16}, the random ordering plays a role in the proof that the algorithm produces a good factorization.
This paper also showed that one can randomize among the vertices with at most twice the average unweighted degree and still get a provably correct algorithm.
In our AC and other solvers, we instead eliminate a vertex with approximately minimum unweighted degree in each round.
To understand the impact of this choice, we compare AC with a variant that uses uniformly random elimination order, which we call AC-random.

In our first experiment, we reuse the setup from Section~\ref{sec:solvervariance}, i.e. we run the AC-random solver 10 times on each of 5 different systems of linear equations and we report the median, 75th percentile, and maximum solve time, normalized by non-zero count, to compare it with our results for AC.
The results of this experiment are shown in Table~\ref{tab:variance-random}.

In our second experiment, we reuse the setup from Section~\ref{sec:variantseval}, i.e. we run the AC-random solver once on each of 5 different systems of linear equations and compute the total solve time per non-zero, the ratio of sizes between the input and the output factorization, and the condition number.
The results of this experiment are shown in Table~\ref{tab:cn-random}.

In both experiments, we see that AC outperforms AC-random except on Sachdeva stars, where AC does worse on all parameters. This is not so surprising, as the Sachdeva star is designed to challenge the greedy ordering approach in particular. 
Concretely, eliminating the central vertex of the Sachdeva star early tends to give worse performance, and the parameters we use for the star are such that the greedy ordering eliminates the central vertex first.
The results suggest that our greedy ordering approach generally works well, but perhaps could be improved slightly by adding a little more randomness to the elimination ordering.

\begin{threeparttable}[H]
  \centering
\resizebox{\textwidth}{!}{%
\begin{tabular}{
  |c
  |S[table-format=3.1]
  |S[table-format=1.2] S[table-format=1.2] S[table-format=1.2]
  |S[table-format=1.2] S[table-format=1.2] S[table-format=2.2]
  |S[table-format=1.2] S[table-format=1.2] S[table-format=1.2]
  |S[table-format=1.2] S[table-format=2.2] S[table-format=2.2]
  |
  }
\hline
  {Instance}
  & {non-zeros}
  & \multicolumn{3}{|c|}{AC { $t_{\text{solve}} / \nnz$ }}
  & \multicolumn{3}{c|}{AC-random { $t_{\text{solve}} / \nnz$ }}
  & \multicolumn{3}{c|}{AC { $t_{\text{total}} / \nnz$ }} 
  & \multicolumn{3}{c|}{AC-random  { $t_{\text{total}} / \nnz$ }}\\ \hline
& {nnz} &{ median   }&{ 0.75    }&{ max   }&{ median   }&{ 0.75    }&{ max   }&{ median  }&{ 0.75   }& {max} &{ median  }&{ 0.75   }& {max} \\ \hline
& {$(M)$}&{$\mu s$}  & {$\mu s$} & {$\mu s$} & {$\mu s$} &{$\mu s$}  & {$\mu s$} & {$\mu s$} &{$\mu s$}  & {$\mu s$} & {$\mu s$} &{$\mu s$}  & {$\mu s$} \\ \hline
Weighted chimera&50.2&3.78&3.87&3.95&2.59&2.6&2.68&4.65&4.74&4.82&4.09&4.11&4.18\\
Weighted SDDM chimera&49.7&3.3&3.38&3.48&2.44&2.5&2.58&4.45&4.53&4.64&4.26&4.33&4.4\\
Sachdeva Star&172&2.56\tnote{*}&3.52\tnote{**}&4.05\tnote{**}&0.208&0.212&0.214&2.96\tnote{*}&3.86\tnote{**}&4.37\tnote{**}&0.977&0.982&0.986\\
High contrast coefficient Poisson grid&200&1.71&1.73&1.77&2.14&2.15&2.25&2.6&2.62&2.66&3.34&3.35&3.45\\
Anisotropic coef. Poisson grid, variable weight&200&1.12&1.14&1.18&1.15&1.16&1.2&1.72&1.75&1.82&2.33&2.34&2.38\\
 \hline
\end{tabular}
} %
 \caption{Variation in total time of AC and AC-random computed across 10 runs per system of linear equations.}
\label{tab:variance-random}
\begin{tablenotes}
\footnotesize
\item[*] relative residual error exceeded tolerance $10^{-8}$ by a factor $(1,10^{4}]$.
\item[**] relative residual error exceeded tolerance $10^{-8}$ by a factor $(10^{4},10^{8})$.
\item[] \hspace{-1em} Inf: Solver crashed or returned a solution with relative residual error 1 (trivial solution).
\item[] \hspace{-1em} N/A: Experiment omitted as the solver crashed too often with (only occurred for PETSc).
\end{tablenotes}
\end{threeparttable}

\begin{threeparttable}[H]
    \centering
  \resizebox{\textwidth}{!}{%
  \begin{tabular}{
    |c
    |S[table-format=3.1]
    |S[table-format=1.3] S[table-format=3.0] S[table-format=1.2] S[table-format=4.1]
    |S[table-format=1.2] S[table-format=3.0] S[table-format=1.2] S[table-format=4.1]
    |
    }
  \hline
    {Instance}
    & {non-zeros}
    & \multicolumn{4}{|c|}{AC}
    & \multicolumn{4}{c|}{AC-random}\\ \hline
  & {nnz} &{ $t_{\text{total}} / \nnz $   }&{$N_{\text{iter}}$}&{ Ratio of sizes }&{ Condition num.   }&{ $t_{\text{total}} / \nnz $   }&{$N_{\text{iter}}$}&{ Ratio of sizes}&{ Condition num.   }  \\ \hline
  & {$(M)$} &{$\mu s$}&  &  & &{$(\mu s)$}&  &  &  \\ \hline
Weighted chimera&50.2&4.42&55&3.56&87.0&9.64&112&5.03&193.0\\
Weighted SDDM chimera&49.7&4.3&51&3.29&59.8&8.47&90&5.19&487.0\\
Sachdeva Star&172&2.2\tnote{*}&408&1.0&9660.0&1.18&210&1.0&4290.0\\
High contrast coefficient Poisson grid&200&2.24&57&2.78&60.2&9.94&145&3.56&436.0\\
Anisotropic coef. Poisson grid, variable weight&200&1.78&41&2.66&24.7&4.97&65&3.61&176.0\\
    \hline
  \end{tabular}
  } %
   \caption{Comparison between approximate Cholesky factorization with different elimination orders.}
  \label{tab:cn-random}
  
  \begin{tablenotes}
\footnotesize
\item[*] relative residual error exceeded tolerance $10^{-8}$ by a factor $(1,10^{4}]$.
\item[**] relative residual error exceeded tolerance $10^{-8}$ by a factor $(10^{4},10^{8})$.
\item[] \hspace{-1em} Inf: Solver crashed or returned a solution with relative residual error 1 (trivial solution).
\item[] \hspace{-1em} N/A: Experiment omitted as the solver crashed too often with (only occurred for PETSc).
\end{tablenotes}
  \end{threeparttable}

\input{scaling}

%% file: suitesparse.tex
\subsubsection{SuiteSparse Matrix Collection}
In this section, we report the performances of the solvers on matrices from the SuiteSparse Matrix Collection \footnote{https://sparse.tamu.edu}, 
previously known as the University of Florida Sparse Matrix Collection \cite{DH11}. In addition to the 
SDDM matrices in SuiteSparse Matrix Collection, we also included matrices that are ``approximately'' SDD, in the 
sense that they have positive diagonals, non-positive off diagonals and is symmetric, while not diagonally dominated but close. 
In particular, for such an $n\times n$ matrix $M$, if $\forall i\in[n],$ such that $(M\vecone)_i<0, \abs{(M\vecone)_i / M_{ii}}\leq \epsilon$ for some small
value $\epsilon$, then we say $M$ is ``approximately'' SDDM. In other words, for these non-SDDM matrices, we also included those satisfying 
$\abs{\min_{i\in [n]}(M\vecone)_i / M_{ii}}\leq \epsilon$ (and we call this value the SDDM-nearness. 
In our experiments, we set this $\epsilon$ to be ten times the machine epsilon. In particular, we included the following non-SDDM matrices:

\begin{itemize}
  \item McRae/ecology1, with SDDM-nearness $1.6\times 10^{-16}$
  \item McRae/ecology2, with SDDM-nearness $1.6\times 10^{-16}$
  \item HB/nos7, with SDDM-nearness $9.09\times 10^{-17}$
\end{itemize}

Similarly, in our experiments, we treat some of the SDDM matrices as approximately Laplacian. In particular, if $\max_{i\in [n} (M\vecone) / M_{ii} \leq \epsilon $, then we say $M$ is approximately Laplacian. If the tolerance $\epsilon$ is small enough, then solving the system while discarding the excess on the diagonals is sufficient for solving the system in $M$.

We excluded matrices with less than $1000$ non-zeros so that the overheads are negligible. 
We also excluded the matrix Cunningham/m3plates, which is a diagonal matrix with zero diagonals. This matrix significantly slowed down some of the solvers, but not AC and AC2.
In total, we tested the solvers on $28$ matrices 
from the SuiteSparse Matrix Collection. We observe that our solvers work well on the ``approximately'' SDDM matrices as well.

Results are shown in Table~\ref{tab:suitesparse}. 
Only AC and {\ACSM} performed well across all the tested matrices from the SuiteSparse Matrix Collection.

\begin{threeparttable}[H]
  \centering
\resizebox{\textwidth}{!}{%
\begin{tabular}{
  |S[table-format=1.3] S[table-format=1.3] S[table-format=1.3]
  |S[table-format=1.3] S[table-format=1.3] S[table-format=1.1]
  |S[table-format=2.1] S[table-format=3.0] S[table-format=3.0]
  |
  }
\hline
    \multicolumn{3}{|c|}{AC $t_{\text{total}} / \nnz $}
  & \multicolumn{3}{c|}{{\ACSM} $t_{\text{total}} / \nnz $}
  & \multicolumn{3}{c|}{CMG $t_{\text{total}} / \nnz $} 
  \\ 
  \cline{1-9}
{ median   }&{ 0.75    }&{ max   }&{ median   }&{ 0.75    }&{ max   }&{ median  }&{ 0.75   }& {max} 
\\ \cline{1-9}
{$(\mu s)$}  & {$(\mu s)$} & {$(\mu s)$} & {$(\mu s)$} &{$(\mu s)$}  & {$(\mu s)$} & {$(\mu s)$} &{$(\mu s)$}  & {$(\mu s)$} 
\\ \hline
0.434&0.503&0.994&0.553&0.698&1.4&6.29&Inf&Inf\\
  \hline
\end{tabular}
} %
\\ \vspace{0.5em}
\resizebox{\textwidth}{!}{%
\begin{tabular}{
  |S[table-format=1.4] S[table-format=1.2] S[table-format=3.0]
  |S[table-format=1.1] S[table-format=2.1] S[table-format=3.0]
  |S[table-format=2.1] S[table-format=2.1] S[table-format=3.0]
  |
  }
\hline
    \multicolumn{3}{|c|}{HyPre $t_{\text{total}} / \nnz $}
  & \multicolumn{3}{c|}{PETSc $t_{\text{total}} / \nnz $}
  & \multicolumn{3}{c|}{ICC $t_{\text{total}} / \nnz $} 
  \\ 
  \cline{1-9}
{ median   }&{ 0.75    }&{ max   }&{ median   }&{ 0.75    }&{ max   }&{ median  }&{ 0.75   }& {max} 
\\ \cline{1-9}
{$(\mu s)$}  & {$(\mu s)$} & {$(\mu s)$} & {$(\mu s)$} &{$(\mu s)$}  & {$(\mu s)$} & {$(\mu s)$} &{$(\mu s)$}  & {$(\mu s)$} 
\\ \hline
0.0586&0.345&Inf&1&1.83&27.1&10.1&48.8&153\\
  \hline
\end{tabular}
} %
 \caption{SuiteSparse Matrix Collection Time / nnz}
\label{tab:suitesparse}
\begin{tablenotes}
\footnotesize
\item[*] relative residual error exceeded tolerance $10^{-8}$ by a factor $(1,10^{4}]$.
\item[**] relative residual error exceeded tolerance $10^{-8}$ by a factor $(10^{4},10^{8})$.
\item[] \hspace{-1em} Inf: Solver crashed or returned a solution with relative residual error 1 (error of the trivial all-zero solution).
\item[] \hspace{-1em} N/A: Experiment omitted as the solver crashed too often with (only occurred for PETSc).
\end{tablenotes}
\end{threeparttable}

%% file: chimeraAndIPM.tex
\newpage
\subsubsection{``Chimera'' Graph Laplacian and SDDM Problems.}
\label{sec:Chimera}

Both our solver and Algebraic Multigrid methods are designed to handle all
symmetric diagonally dominant systems of linear equations, not only
those arising from discretizations of 3D cube problems.

In this section, we study the performance of solvers on SDDM and
Laplacian systems of linear equations with a wide range of non-zero
structures.
We first focus on graph Laplacians arising from a class of recursively generated graphs which we call \emph{Chimeras}.

Chimeras are formed by randomly combining standard base graphs.
We generate Chimera graphs with a given vertex count by a recursive process which uses a seeded pseudo-random generator to pseudo-randomly choose among the following options. 
\begin{enumerate}[label=(\alph*)]
    \item Form a graph from one of several classes: paths, trees, 2D grids, rings and generalized rings, the largest connected component of an Erd\H{o}s-Renyi graph, random regular graphs, and random preferential-attachment-like graphs.
    \item Recursively form two smaller Chimeras and combine them in one of several ways:
    by adding random edges between them, or forming their Cartesian graph product, or forming a `generalized necklace' from the two graphs.
    We define a generalized necklace of two graphs $G$ and $H$ as the result of expanding each vertex in $G$ to an
instance of $H$, and for each edge in $G$ adding a number of random edges between the corresponding copies of $H$ in the new graph.
    \item Select one of two operations that create a larger graph from a single smaller graph.
    The first operation is a pseudo-random two-lift, which doubles the vertex set and replace each edge with either a pair of internal edges internal to the original set and the copy or a pair crossing edges between the original and the copy.
    The second operation is `thickening,' which adds a subset of the current 2-hop paths as new edges in the graph.
\end{enumerate}
We consider both unweighted Chimeras and weighted variants.
Weights are obtained by a pseudo-random process that first either chooses unform edge weights in $\left[0,1\right]$ or chooses random potentials to assign vertices and then sets the weights of edges to the differences between the potentials at the endpoints.
These potentials can either be uniform in $\left[0,1\right]$, or the result of multiplying such vertex potentials by a small power of the  normalized Laplacian matrix or the corresponding walk matrix.
Finally, with probability one half, all the edge weights are replaced by their reciprocals.

For each vertex count $n$, when we report statistics for $C$ different instances, these are always the Chimeras with $n$ vertices and seed indices $i = 1, \ldots, C$. 
We always choose the Chimeras generated by the first seed indices, and we do not exclude any Chimeras. 
Our Chimera generator is deterministic given the seed.
When using Chimeras for benchmarking, we strongly encourage this approach, as arbitrary exclusions could give misleading statistics.

We also consider a variant of the graph Laplacian problem, where we
modify the linear system by %
setting a boundary condition of $u(x) = 0$ for $\sizeof{V}^{2/3}$ of the vertices $x \in V$.
The boundary size of $\sizeof{V}^{2/3}$ was chosen to get a similar
boundary size as in the 3D cube examples.
We enforce the boundary condition at vertices with index divisible by $\sizeof{V}^{1/3}$, which spreads the boundary vertices across the graph.
Expanding these tests to include `adversarially' chosen boundary conditions remains an interesting open question.
This modification results in a strictly SDDM linear equation.

Results for unweighted and weighted Laplacian Chimeras are shown in Table~\ref{tab:unilapchimera} and Table~\ref{tab:wtedlapchimera} respectively.
Results for unweighted and weighted SDDM Chimeras are shown in Table~\ref{tab:unisddchimera} and Table~\ref{tab:wtedsddmlapchimera} respectively. Our experiments showed that 
Chimeras are a particularly challenging problem because only AC and {\ACSM}
successfully reached the desired tolerance $10^{-8}$ for all problem
instances.
Beyond this worst case behavior, our solver AC also achieved the best median runtimes, and 75th
percentile runtime, beating HyPre and CMG 
which were best solvers for the Poisson grid problems. 
\newpage

\begin{threeparttable}[h]
    \centering
  \resizebox{\textwidth}{!}{%
  \begin{tabular}{
    |S[table-format=5.0]
    |S[table-format=3.0]
    |S[table-format=1.3] S[table-format=1.3] S[table-format=1.3]
    |S[table-format=1.3] S[table-format=1.3] S[table-format=1.2]
    |S[table-format=1.2] S[table-format=1.2] S[table-format=2.1]
    |
    }
\hline
    { variables } & {\# instances}
    & \multicolumn{3}{c|}{AC $t_{\text{total}} / \nnz $}
    & \multicolumn{3}{c|}{{\ACSM} $t_{\text{total}} / \nnz $}
    & \multicolumn{3}{c|}{CMG $t_{\text{total}} / \nnz $} 
    \\ \hline
  {  $n$ }             &                            
   &{ median   }&{ 0.75    }&{ max   }&{ median   }&{ 0.75    }&{ max   }&{ median  }&{ 0.75   }& {max} 
   \\ \hline
  { $(K)$}  &  &{$(\mu s)$}  & {$(\mu s)$} & {$(\mu s)$} & {$(\mu s)$} &{$(\mu s)$}  & {$(\mu s)$} & {$(\mu s)$} &{$(\mu s)$}  & {$(\mu s)$} 
  \\ \hline
  10&103&0.421&0.467&0.576&0.687&0.804&1.24&6.5&7.71&16.1\\
100&105&0.676&0.87&1.56&1.08&1.44&3.18&1.54&2.04&Inf\\
1000&23&1.5&1.73&2.39&1.72&2.83&3.58&1.75&2.44&Inf\\
10000&8&2.79&3.4&3.95&3.57&5.06&7.19&3.99&5.45&20.8\\
    \hline
  \end{tabular}
} %
\\ \vspace{0.5em}
  \resizebox{\textwidth}{!}{%
  \begin{tabular}{
    |S[table-format=5.0]
    |S[table-format=3.0]
    |S[table-format=1.3] S[table-format=1.2] S[table-format=2.2]
    |S[table-format=2.2] S[table-format=3.1] S[table-format=3.1]
    |S[table-format=1.2] S[table-format=1.2] S[table-format=1.2]
    |
    }
\hline
    { variables } & {\# instances}
    & \multicolumn{3}{c|}{HyPre $t_{\text{total}} / \nnz $}
                        & \multicolumn{3}{c|}{PETSc $t_{\text{total}} / \nnz $}
    & \multicolumn{3}{c|}{ICC $t_{\text{total}} / \nnz $} \\ \hline
  {  $n$ }             &                             &{ median   }&{ 0.75    }&{ max   }&{ median   }&{ 0.75    }&{ max   }&{ median  }&{ 0.75   }& {max}  \\ \hline
  { $(K)$}  &  &{$(\mu s)$}  & {$(\mu s)$} & {$(\mu s)$} & {$(\mu s)$} &{$(\mu s)$}  & {$(\mu s)$} & {$(\mu s)$} &{$(\mu s)$}  & {$(\mu s)$} \\ \hline
  10&103&0.919&1.31&2.86\tnote{**}&5.86&11.9&38.8\tnote{*}&5.09&6.57&Inf\\
100&105&1.38&2.29&6.23\tnote{**}&18.9&48.9&182\tnote{*}&1.68&2.7&Inf\\
1000&23&3.19&5.73&15.4&97.6&370\tnote{*}&Inf&1.59&3.16&Inf\\
10000&8&6.7&7.8&Inf&Inf&Inf&Inf&3.04&4.82&5.63\\
    \hline
  \end{tabular}
} %
   \caption{
     Unweighted Laplacian Chimera Time / nnz}
  \label{tab:unilapchimera}
  \begin{tablenotes}
\footnotesize
\item[*] relative residual error exceeded tolerance $10^{-8}$ by a factor $(1,10^{4}]$.
\item[**] relative residual error exceeded tolerance $10^{-8}$ by a factor $(10^{4},10^{8})$.
\item[] \hspace{-1em} Inf: Solver crashed or returned a solution with relative residual error 1 (error of the trivial all-zero solution).
\item[] \hspace{-1em} N/A: Experiment omitted as the solver crashed too often with (only occurred for PETSc).
\end{tablenotes}
\end{threeparttable}

\vspace{1em}

\begin{threeparttable}[H]
  \centering
\resizebox{\textwidth}{!}{%
\begin{tabular}{
  |S[table-format=5.0]
  |S[table-format=3.0]
  |S[table-format=1.3] S[table-format=1.3] S[table-format=1.3]
  |S[table-format=1.3] S[table-format=1.3] S[table-format=1.3]
  |S[table-format=1.2] S[table-format=1.2] S[table-format=2.1]
  |
  }
\hline
  { variables } & {\# instances}
  & \multicolumn{3}{c|}{AC $t_{\text{total}} / \nnz $}
  & \multicolumn{3}{c|}{{\ACSM} $t_{\text{total}} / \nnz $}
  & \multicolumn{3}{c|}{CMG $t_{\text{total}} / \nnz $} 
  \\ \hline
{  $n$ }             &                             &{ median   }&{ 0.75    }&{ max   }&{ median   }&{ 0.75    }&{ max   }&{ median  }&{ 0.75   }& {max} 
\\ \hline
{ $(K)$}  &  &{$(\mu s)$}  & {$(\mu s)$} & {$(\mu s)$} & {$(\mu s)$} &{$(\mu s)$}  & {$(\mu s)$} & {$(\mu s)$} &{$(\mu s)$}  & {$(\mu s)$} 
\\ \hline
10&103&0.365&0.405&0.455&0.572&0.729&0.939&7.72&9.44&16.4\\
100&105&0.605&0.745&1.31&0.908&1.26&2.25&1.87&2.54&Inf\\
1000&23&1.5&1.8&2.19&1.67&2.53&3.37&1.94&2.64&Inf\\
10000&8&3.16&4.21&4.39&3.83&4.19&4.79&2.89&5.04&11.9\\
  \hline
\end{tabular}
} %
\\ \vspace{0.5em}
\resizebox{\textwidth}{!}{%
\begin{tabular}{
  |S[table-format=5.0]
  |S[table-format=3.0]
  |S[table-format=1.3] S[table-format=2.2] S[table-format=2.2]
  |S[table-format=2.2] S[table-format=3.1] S[table-format=3.1]
  |S[table-format=1.2] S[table-format=2.2] S[table-format=2.1]
  |
  }
\hline
  { variables } & {\# instances}
  & \multicolumn{3}{c|}{HyPre $t_{\text{total}} / \nnz $}
                      & \multicolumn{3}{c|}{PETSc $t_{\text{total}} / \nnz $}
  & \multicolumn{3}{c|}{ICC $t_{\text{total}} / \nnz $} \\ \hline
{  $n$ }             &                             &{ median   }&{ 0.75    }&{ max   }&{ median   }&{ 0.75    }&{ max   }&{ median  }&{ 0.75   }& {max}  \\ \hline
{ $(K)$}  &  &{$(\mu s)$}  & {$(\mu s)$} & {$(\mu s)$} & {$(\mu s)$} &{$(\mu s)$}  & {$(\mu s)$} & {$(\mu s)$} &{$(\mu s)$}  & {$(\mu s)$} \\ \hline
10&103&0.968&1.36&5.18\tnote{**}&4.58&11.2&32.9\tnote{*}&5.7&7.01&Inf\\
100&105&1.25&2.3&7.7\tnote{**}&12.9&33.3&151\tnote{*}&2.59&5.07&Inf\\
1000&23&3.24&5.8&15.8\tnote{*}&52.5&218&861\tnote{*}&3.54&7.67&Inf\\
10000&8&4.88&10.7\tnote{**}&Inf&N/A&N/A&N/A&7.27&11.6\tnote{*}&15.2\tnote{*}\\
  \hline
\end{tabular}
} %
 \caption{
   Weighted Laplacian Chimera Time / nnz}
\label{tab:wtedlapchimera}
\begin{tablenotes}
\footnotesize
\item[*] relative residual error exceeded tolerance $10^{-8}$ by a factor $(1,10^{4}]$.
\item[**] relative residual error exceeded tolerance $10^{-8}$ by a factor $(10^{4},10^{8})$.
\item[] \hspace{-1em} Inf: Solver crashed or returned a solution with relative residual error 1 (error of the trivial all-zero solution).
\item[] \hspace{-1em} N/A: Experiment omitted as the solver crashed too often with (only occurred for PETSc).
\end{tablenotes}
\end{threeparttable}

\newpage

\begin{threeparttable}[H]
  \centering
\resizebox{\textwidth}{!}{%
\begin{tabular}{
  |S[table-format=4.2]
  |S[table-format=3.0]
  |S[table-format=1.3] S[table-format=1.3] S[table-format=1.2]
  |S[table-format=1.3] S[table-format=1.2] S[table-format=1.3]
  |S[table-format=1.2] S[table-format=1.2] S[table-format=2.2]
  |
  }
\hline
  { variables } & {\# instances}
                      & \multicolumn{3}{c|}{AC $t_{\text{total}} / \nnz $}
  & \multicolumn{3}{c|}{{\ACSM} $t_{\text{total}} / \nnz $}
  & \multicolumn{3}{c|}{CMG $t_{\text{total}} / \nnz $} 
  \\ \hline
{  $n$ }             &                             &{ median   }&{ 0.75    }&{ max   }&{ median   }&{ 0.75    }&{ max   }&{ median  }&{ 0.75   }& {max}
\\ \hline
{ $(K)$}  &  &{$(\mu s)$}  & {$(\mu s)$} & {$(\mu s)$} & {$(\mu s)$} &{$(\mu s)$}  & {$(\mu s)$} & {$(\mu s)$} &{$(\mu s)$}  & {$(\mu s)$}
\\ \hline
9.54&103&0.379&0.391&0.42&0.543&0.6&0.739&7.18&9.66&16.2\\
97.9&116&0.585&0.684&1.5&0.878&1.19&1.9&1.43&1.85&Inf\\
990&29&1.43&1.62&1.8&1.58&2.01&2.83&1.58&2.29&Inf\\
9950&12&2.7&2.95&3.4&2.84&3.47&5.44&3.76&4.91&Inf\\
  \hline
\end{tabular}
} %
\\ \vspace{0.5em}
\resizebox{\textwidth}{!}{%
\begin{tabular}{
  |S[table-format=4.2]
  |S[table-format=3.0]
  |S[table-format=1.3] S[table-format=1.2] S[table-format=2.2]
  |S[table-format=2.2] S[table-format=2.2] S[table-format=3.1]
  |S[table-format=1.3] S[table-format=1.2] S[table-format=1.2]
  |
  }
\hline
  { variables } & {\# instances}
  & \multicolumn{3}{c|}{HyPre $t_{\text{total}} / \nnz $}
                      & \multicolumn{3}{c|}{PETSc $t_{\text{total}} / \nnz $}
  & \multicolumn{3}{c|}{ICC $t_{\text{total}} / \nnz $} \\ \hline
{  $n$ }             &                             &{ median   }&{ 0.75    }&{ max   }&{ median   }&{ 0.75    }&{ max   }&{ median  }&{ 0.75   }& {max}  \\ \hline
{ $(K)$}  &  &{$(\mu s)$}  & {$(\mu s)$} & {$(\mu s)$} & {$(\mu s)$} &{$(\mu s)$}  & {$(\mu s)$} & {$(\mu s)$} &{$(\mu s)$}  & {$(\mu s)$} \\ \hline
9.54&103&0.904&1.09&2.18&3.56&7.71&21.4&5.12&6.68&9.26\\
97.9&116&1.25&1.75&4.89&11.2&35&167&1.15&1.51&3.39\\
990&29&2.8&4.14&11.4&32.7&86.5&Inf&0.901&1.2&Inf\\
9950&12&5.92&7.73&Inf&Inf&Inf&Inf&1.28&1.44&3.61\\
  \hline
\end{tabular}
} %
 \caption{Unweighted SDDM Chimera Time / nnz}
\label{tab:unisddchimera}
\begin{tablenotes}
\footnotesize
\item[*] relative residual error exceeded tolerance $10^{-8}$ by a factor $(1,10^{4}]$.
\item[**] relative residual error exceeded tolerance $10^{-8}$ by a factor $(10^{4},10^{8})$.
\item[] \hspace{-1em} Inf: Solver crashed or returned a solution with relative residual error 1 (error of the trivial all-zero solution).
\item[] \hspace{-1em} N/A: Experiment omitted as the solver crashed too often with (only occurred for PETSc).
\end{tablenotes}
\end{threeparttable}

\vspace{1em}

\begin{threeparttable}[H]
  \centering
\resizebox{\textwidth}{!}{%
\begin{tabular}{
  |S[table-format=4.2]
  |S[table-format=3.0]
  |S[table-format=1.3] S[table-format=1.3] S[table-format=1.3]
  |S[table-format=1.3] S[table-format=1.3] S[table-format=1.3]
  |S[table-format=1.2] S[table-format=2.2] S[table-format=2.2]
  |
  }
\hline
  { variables } & {\# instances}
  & \multicolumn{3}{c|}{AC $t_{\text{total}} / \nnz $}
  & \multicolumn{3}{c|}{{\ACSM} $t_{\text{total}} / \nnz $}
  & \multicolumn{3}{c|}{CMG $t_{\text{total}} / \nnz $} 
  \\ \hline
{  $n$ }             &                             &{ median   }&{ 0.75    }&{ max   }&{ median   }&{ 0.75    }&{ max   }&{ median  }&{ 0.75   }& {max}
 \\ \hline
{ $(K)$}  &  &{$(\mu s)$}  & {$(\mu s)$} & {$(\mu s)$} & {$(\mu s)$} &{$(\mu s)$}  & {$(\mu s)$} & {$(\mu s)$} &{$(\mu s)$}  & {$(\mu s)$} 
\\ \hline
9.54&103&0.375&0.385&0.431&0.537&0.592&0.714&7.6&10.6&17.1\\
97.9&105&0.606&0.701&1.43&0.909&1.11&1.82&1.44&1.86&2.54\tnote{**}\\
990&23&1.4&1.59&2.16&1.66&2.14&2.42&1.55&1.95&Inf\\
9950&8&3.07&3.99&4.34&3.34&4.18&4.36&2.65&4.19&12.6\\
  \hline
\end{tabular}
} %
\\ \vspace{0.5em}
\resizebox{\textwidth}{!}{%
\begin{tabular}{
  |S[table-format=4.2]
  |S[table-format=3.0]
  |S[table-format=1.3] S[table-format=1.2] S[table-format=2.2]
  |S[table-format=3.0] S[table-format=3.0] S[table-format=3.0]
  |S[table-format=1.2] S[table-format=1.2] S[table-format=2.1]
  |
  }
\hline
  { variables } & {\# instances}
  & \multicolumn{3}{c|}{HyPre $t_{\text{total}} / \nnz $}
                      & \multicolumn{3}{c|}{PETSc $t_{\text{total}} / \nnz $}
  & \multicolumn{3}{c|}{ICC $t_{\text{total}} / \nnz $} \\ \hline
{  $n$ }             &                             &{ median   }&{ 0.75    }&{ max   }&{ median   }&{ 0.75    }&{ max   }&{ median  }&{ 0.75   }& {max}  \\ \hline
{ $(K)$}  &  &{$(\mu s)$}  & {$(\mu s)$} & {$(\mu s)$} & {$(\mu s)$} &{$(\mu s)$}  & {$(\mu s)$} & {$(\mu s)$} &{$(\mu s)$}  & {$(\mu s)$} \\ \hline
9.54&103&0.873&1.26&2.4&N/A&N/A&N/A&4.45&5.65&10.5\\
97.9&105&0.946&1.66&4.73&N/A&N/A&N/A&1.04&1.39&10.3\tnote{*}\\
990&23&2.61&4.29&11.4\tnote{*}&N/A&N/A&N/A&1.06&2.63&Inf\\
9950&8&4.39&6.4&Inf&N/A&N/A&N/A&2&4.06&14.5\tnote{*}\\
  \hline
\end{tabular}
} %
 \caption{
   Weighted SDDM Chimera Time / nnz}
\label{tab:wtedsddmlapchimera}
\begin{tablenotes}
\footnotesize
\item[*] relative residual error exceeded tolerance $10^{-8}$ by a factor $(1,10^{4}]$.
\item[**] relative residual error exceeded tolerance $10^{-8}$ by a factor $(10^{4},10^{8})$.
\item[] \hspace{-1em} Inf: Solver crashed or returned a solution with relative residual error 1 (error of the trivial all-zero solution).
\item[] \hspace{-1em} N/A: Experiment omitted as the solver crashed too often with (only occurred for PETSc).
\end{tablenotes}
\end{threeparttable}

\newpage
\subsubsection{Sachdeva Star Graph Laplacians.}
\begin{figure}[H]
  \centering
    \includegraphics[width=0.3\textwidth]{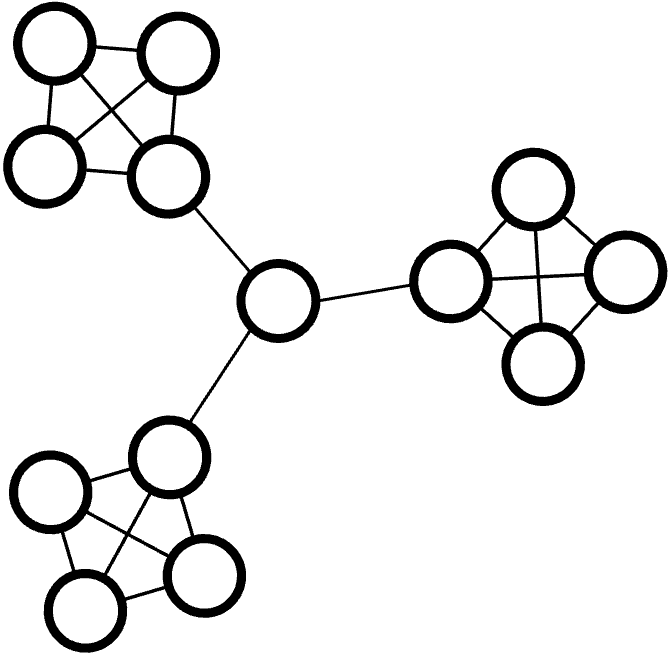}
  \caption{Sachdeva star with $l = 3$ and $k = 4$.}
  \label{fig:sachdevastar}
\end{figure}

In this section, we report the performance of solvers on a graph
Laplacian specifically designed to make our
$\textsc{ApproximateCholesky}$/$\textsc{ApproximateEdgewiseCholesky}$ algorithm fail.
The construction was suggested by Sushant Sachdeva.
We construct a star graph with $l$ leaves,
and then replace each leaf with a complete graph on $k$ vertices. 
In the experiments, we set $l$ to be $\frac{k}{2}$.
We then consider a linear equation in the associated graph Laplacian.
An example is shown in Figure~\ref{fig:sachdevastar}.

Results are shown in Table~\ref{tab:startable}. This is the most challenging 
problem for AC. The iteration counts of AC spiked when compared
to other problems. This suggests that AC failed to produce reliable 
factorizations for the Sachdeva star Laplacians. In contrast, {\ACSM}
produced much more reliable factorizations and the iteration counts are 
still comparable to other problems of similar sizes. In Section~\ref{sec:variantseval}
we have a more detailed comparison between variants of our approximate 
Cholesky solvers on Sachdeva star. We observe that CMG performs the best
for this problem, while HyPre, PETSc and ICC failed to achieve the desired
tolerance on most instances.

\begin{landscape}
  \newcommand{\lstablelen}{24cm}
  \begin{threeparttable}[H]
    \centering
  \resizebox{\lstablelen}{!}{%
  \begin{tabular}{
    |S[table-format=3.0]
    |S[table-format=6.0]
    |S[table-format=3.0]
    |S[table-format=3.0] S[table-format=1.3] S[table-format=1.3] S[table-format=3.3e1]
    |S[table-format=2.0] S[table-format=1.3] S[table-format=1.3] S[table-format=1.3]
    |S[table-format=1.0] S[table-format=1.4] S[table-format=1.3] S[table-format=1.2e-1]
    |
    }
  \hline
    { $k$ }& { nonzeros }&{ variables }&
                                                    \multicolumn{4}{c|}{AC}     & \multicolumn{4}{c|}{{\ACSM}}     & \multicolumn{4}{c|}{CMG} 
                                                     \\
  \hline
    { } & {$\nnz$}  &  {$n$} 
        &{ $N_{\text{iter}}$ }&{ $t_{\text{solve}} / \nnz$ }&{ $t_{\text{total}} / \nnz$ }&{ $\text{res. err.}$ }
        &{ $N_{\text{iter}}$ }&{ $t_{\text{solve}} / \nnz$ }&{ $t_{\text{total}} / \nnz$ }&{ $\text{res. err.}$ }
        &{ $N_{\text{iter}}$ }&{ $t_{\text{solve}} / \nnz$ }&{ $t_{\text{total}} / \nnz$ }&{ $\text{res. err.}$ }
    \\ \hline %
        & {$(K)$} & {$(K)$}  
        & & {$(\mu s)$} & {$(\mu s)$} & {$(10^{-8})$}
        & & {$(\mu s)$} & {$(\mu s)$} & {$(10^{-8})$}
        & & {$(\mu s)$} & {$(\mu s)$} & {$(10^{-8})$}
    \\
    \hline
    100&500&5&83&0.226&0.398&0.9&28&0.0793&0.478&0.64&7&0.431&0.868&3.85e-06\\
150&1690&11&149&0.502&0.748&0.96&34&0.117&0.581&0.632&7&0.181&0.486&2.46e-06\\
200&4000&20&167&0.682&0.931&0.955&37&0.155&0.716&0.877&7&0.124&0.405&1.56e-06\\
250&7810&31&241&1.07&1.43&0.998&38&0.169&0.703&0.98&7&0.107&0.442&1.09e-06\\
300&13500&45&289&1.33&1.59&0.88&39&0.179&0.719&0.854&7&0.0962&0.436&8.87e-07\\
350&21400&61&355&1.65&1.97&0.99&40&0.188&0.733&0.996&7&0.0905&0.426&8.18e-07\\
400&32000&80&459&2.12&2.44&0.967&40&0.186&0.718&0.985&7&0.0884&0.423&8.17e-07\\
450&45600&101&486&2.27&2.58&0.999&43&0.202&0.745&0.902&7&0.0848&0.418&6.6e-07\\
500&62500&125&452&2.08&2.4&0.944&43&0.198&0.747&0.808&7&0.0833&0.415&6.01e-07\\
550&83200&151&149&0.698\tnote{**}&1.02\tnote{**}&1.88e+04&41&0.192&0.753&0.792&7&0.0826&0.417&1.83e-06\\
600&108000&180&773&3.55&3.93&0.997&45&0.208&0.849&0.741&7&0.0821&0.415&5.94e-07\\
650&137000&211&953&4.43\tnote{*}&4.76\tnote{*}&2.4&44&0.201&1.04&0.687&7&0.0823&0.41&1.3e-06\\
700&172000&245&478&2.19\tnote{*}&2.62\tnote{*}&126&45&0.207&0.859&0.882&7&0.0823&0.412&8.1e-07\\
750&211000&281&905&4.15\tnote{*}&4.47\tnote{*}&2.33&43&0.198&0.855&0.916&7&0.0813&0.409&1.41e-06\\
800&256000&320&669&3.07\tnote{*}&3.39\tnote{*}&32.9&45&0.207&0.883&0.956&7&0.0812&0.408&1.35e-06\\
    \hline
  \end{tabular}
  } %
  \resizebox{\lstablelen}{!}{%
  \begin{tabular}{
    |S[table-format=3.0]
    |S[table-format=6.0]
    |S[table-format=3.0]
    |S[table-format=3.0] S[table-format=1.3] S[table-format=1.3] S[table-format=3.3e1]
    |S[table-format=1.0] S[table-format=1.4] S[table-format=1.2] S[table-format=1.2e1]
    |S[table-format=3.0] S[table-format=1.6] S[table-format=1.3] S[table-format=1.2e-1]
    |
    }
  \hline
    { $k$ }& { nonzeros }&{ variables }
    & \multicolumn{4}{c|}{HyPre}
                                                    &\multicolumn{4}{c|}{PETSc}     %
                                                     & \multicolumn{4}{c|}{ICC}        
                                                     \\
  \hline
    { } & {$\nnz$}  &  {$n$} 
        &{ $N_{\text{iter}}$ }&{ $t_{\text{solve}} / \nnz$ }&{ $t_{\text{total}} / \nnz$ }&{ $\text{res. err.}$ }
        &{ $N_{\text{iter}}$ }&{ $t_{\text{solve}} / \nnz$ }&{ $t_{\text{total}} / \nnz$ }&{ $\text{res. err.}$ }
        &{ $N_{\text{iter}}$ }&{ $t_{\text{solve}} / \nnz$ }&{ $t_{\text{total}} / \nnz$ }&{ $\text{res. err.}$ }
    \\ \hline %
        & {$(K)$} & {$(K)$}  
        & & {$(\mu s)$} & {$(\mu s)$} & {$(10^{-8})$}
        & & {$(\mu s)$} & {$(\mu s)$} & {$(10^{-8})$}
        & & {$(\mu s)$} & {$(\mu s)$} & {$(10^{-8})$}
    \\
    \hline
    100&500&5&100&0.44\tnote{*}&0.84\tnote{*}&225&2&0.032\tnote{**}&1.08\tnote{**}&2.76e+04&Inf&Inf&Inf&1e+08\\
150&1690&11&100&0.45\tnote{*}&0.995\tnote{*}&493&2&0.0341\tnote{**}&1.47\tnote{**}&2.63e+04&Inf&Inf&Inf&1e+08\\
200&4000&20&100&0.692\tnote{*}&1.4\tnote{*}&2.87e+03&2&0.0376\tnote{**}&1.93\tnote{**}&2.09e+04&1&0.0789&0.168&8.83e-08\\
250&7810&31&100&0.852\tnote{*}&1.73\tnote{*}&15.2&1&0.041\tnote{**}&2.37\tnote{**}&1.68e+04&1&0.0475&0.146&9.52e-08\\
300&13500&45&100&0.849\tnote{**}&1.89\tnote{**}&3.91e+04&2&0.0427\tnote{**}&2.81\tnote{**}&1.56e+04&1&0.0316&0.14&1.02e-07\\
350&21400&61&100&0.855\tnote{**}&2.07\tnote{**}&1.14e+06&5&0.106\tnote{*}&3.31\tnote{*}&4.21e+03&Inf&Inf&Inf&1e+08\\
400&32000&80&100&0.842\tnote{*}&2.22\tnote{*}&1.54e+03&1&0.0394\tnote{**}&3.69\tnote{**}&1.2e+04&1&0.022&0.151&1.3e-07\\
450&45600&101&100&0.852\tnote{**}&2.4\tnote{**}&2.71e+06&2&0.0388\tnote{**}&4.13\tnote{**}&1.21e+04&1&0.0202&0.163&1.32e-07\\
500&62500&125&100&0.848\tnote{**}&2.55\tnote{**}&1.47e+05&1&0.0374\tnote{**}&4.57\tnote{**}&1.03e+04&Inf&Inf&Inf&1e+08\\
550&83200&151&100&0.834\tnote{*}&2.7\tnote{*}&3.36e+03&2&0.0371\tnote{*}&5.05\tnote{*}&9.35e+03&1&0.0183&0.188&1.46e-07\\
600&108000&180&47&0.396&2.42&0.562&2&0.0572\tnote{*}&5.47\tnote{*}&9.17e+03&1&0.0177&0.2&1.51e-07\\
650&137000&211&100&0.831\tnote{*}&3.02\tnote{*}&4.37e+03&2&0.0568\tnote{*}&5.93\tnote{*}&8.53e+03&1&0.0169&0.212&1.62e-07\\
700&172000&245&100&0.831\tnote{*}&3.19\tnote{*}&236&3&0.0567\tnote{*}&6.35\tnote{*}&8.93e+03&Inf&Inf&Inf&1e+08\\
750&211000&281&100&0.828\tnote{**}&3.35\tnote{**}&5.89e+04&2&0.0371\tnote{*}&6.77\tnote{*}&8.59e+03&1&0.0161&0.235&1.71e-07\\
800&256000&320&100&0.821\tnote{*}&3.52\tnote{*}&2.44e+03&2&0.0388\tnote{*}&7.22\tnote{*}&8.43e+03&Inf&Inf&Inf&1e+08\\
    \hline
  \end{tabular}
  } %
  \caption{Sachdeva Star with ($k/2$) Leaves Replaced by Complete Graph on $k$ Vertices}
  \label{tab:startable}
  \begin{tablenotes}
\footnotesize
\item[*] relative residual error exceeded tolerance $10^{-8}$ by a factor $(1,10^{4}]$.
\item[**] relative residual error exceeded tolerance $10^{-8}$ by a factor $(10^{4},10^{8})$.
\item[] \hspace{-1em} Inf: Solver crashed or returned a solution with relative residual error 1 (trivial solution).
\item[] \hspace{-1em} N/A: Experiment omitted as the solver crashed too often with (only occurred for PETSc).
\end{tablenotes}
  \end{threeparttable}
\end{landscape}

\subsubsection{Interior Point Method Matrices}
In this section we report the performances of solvers on matrices
generated from running an Interior Point Method on 
undirected Maximum Flow problems on two classes of graphs.
Our Maximum Flow interior point solver is based on the~\cite{KLS20} framework, and uses too many iterations to be practical.
Hence, our main interest is to examine if the resulting linear equations are tractable using various solvers.

\paragraph{IPM on Chimeras.} First, we tested the solvers using Laplacians generated from running a
Maximum Flow IPM unweighted on Chimeras (see section~\ref{sec:Chimera} for
details on Chimera graphs).
We ran the IPM on five Chimeras with $10^5$ vertices and with the
generating seed index $i$ ranging from $1$ to $5$.
We compute and store a number of the Laplacian linear equations
encountered during a run of the IPM:
Concretely, the IPM has a homotopy parameter, which measure roughly
how close we are to having an optimal maximum flow solution.
As the algorithm progresses, the homotopy parameter changes and the
currently maintained approximate maximum flow solution improves in quality. 
We collect matrices at five different stages, corresponding to
different values of the homotopy parameter, and hence different
current solution qualites. 
First, we collect matrices when the maximum flow solution is guaranteed error at
most $10^{-1}$ by the homotopy parameter; then we collect when the
guaranteed error is $10^{-2}, 10^{-3}, 10^{-4}$ and $10^{-5}$.
At each of these five collection stages, we store roughly $5$
Laplacians corresponding to consequtive Newton steps in each precision
range, for a total of 25 Laplacians per Maximum Flow problem.

In general, the weights in the Laplacian matrices will get increasingly
extreme as the precision parameter decreases (meaning the error is smaller).
Thus we might expect the problems to become increasingly difficult as
the precision parameter decreases.
Table~\ref{tab:chimeraIPM} reports detailed performances of the
solvers in each precision range, showing the median, 75th percentile,
and worst case running times across the five Chimera instances we tested.
``AC'', ``AC2'', ``CMG'' and ``ICC'' reached the desired tolerance on
all instances.
``AC'' is clearly the fastest solver across precision ranges, while
the other solvers still managed a comparable run time.
We omitted PETSc because we found it to be very unstable on these
IPM-weighted Chimera Laplacians and frequently crashed.

\paragraph{IPM on Spielman graphs.} After our Chimera IPM experiments, we tested a family of graphs called
Spielman graphs.
These are conjectured to be a hard instance for
short-step Maximum Flow IPMs -- i.e. short-step IPMs need many Newton steps to
converge on graphs.
Spielman graphs are built recursively across a
number of levels. There is unique Spielman graph with $k$ levels,
and it has $O(k^3)$ vertices and edges.
This family has a high degree of symmetry, which allows an algorithm
based on implicit representations (which we call R-space representation)
of the flow solution to compute intermediate states of the IPM very
quickly.
We use this to speed up
our experiments, by only employing a Laplacian solver at a few points
spaced evenly throughout a run of the IPM, and using R-space
representation for majority of the Newton steps.
We run the IPM up to precision $10^{-6}$. For each number of levels
$k = 100,200,300,400,500,600$ (or in terms of number of edges of
the Spielman graphs: from roughly a million up to $200$ millions), we take roughly $10$ Laplacians that occurs when running 
Maximum Flow IPM.
The Laplacians are distributed evenly across precision ranges include
the first and the final Laplacians during each run of the IPM. 
Table~\ref{tab:ipmrloweps} reports the details of the performances of
the solvers on Laplcians from running our Maximum Flow IPM on
Spielman graphs.
Only ``AC'', ``AC2'' and ``HyPre'' reached the desired tolerence on all instances. The running time 
of these three solvers are roughly comparable.
The results of this experiment are shown in
Table~\ref{tab:ipmrloweps}.
\newpage 
\begin{threeparttable}[H]
  \centering
\resizebox{\textwidth}{!}{%
\begin{tabular}{
  |S[table-format=3.2]
  |S[table-format=2.0]
  |S[table-format=1.3] S[table-format=1.3] S[table-format=1.2]
  |S[table-format=1.2] S[table-format=1.2] S[table-format=1.2]
  |S[table-format=1.2] S[table-format=1.2] S[table-format=1.2]
  |
  }
\hline
  { IPM param. } & {\# instances}
  & \multicolumn{3}{c|}{AC $t_{\text{total}} / \nnz $}
  & \multicolumn{3}{c|}{{\ACSM} $t_{\text{total}} / \nnz $}
  & \multicolumn{3}{c|}{CMG $t_{\text{total}} / \nnz $} 
  \\ \hline
{  }             &                             &{ median   }&{ 0.75    }&{ max   }&{ median   }&{ 0.75    }&{ max   }&{ median  }&{ 0.75   }& {max} 
\\ \hline
{ $(10^{-3})$}  &  &{$(\mu s)$}  & {$(\mu s)$} & {$(\mu s)$} & {$(\mu s)$} &{$(\mu s)$}  & {$(\mu s)$} & {$(\mu s)$} &{$(\mu s)$}  & {$(\mu s)$} 
\\ \hline
100&28&0.68&0.774&1.2&1.27&1.79&2.39&1.2&1.61&1.8\\
10&27&0.678&0.777&1.15&1.27&1.8&2.36&1.2&1.56&1.69\\
1&27&0.68&0.763&1.14&1.25&1.79&2.38&1.21&1.66&1.78\\
0.1&27&0.683&0.784&1.15&1.27&1.79&2.36&1.19&1.58&1.78\\
0.01&19&0.681&0.777&1.15&1.28&1.78&2.37&1.22&1.59&1.9\\
  \hline
\end{tabular}
} %
\\ \vspace{0.5em}
\resizebox{\textwidth}{!}{%
\begin{tabular}{
  |S[table-format=3.2]
  |S[table-format=2.0]
  |S[table-format=1.2] S[table-format=1.2] S[table-format=1.2]
  |S[table-format=3.0] S[table-format=3.0] S[table-format=3.0]
  |S[table-format=1.3] S[table-format=1.2] S[table-format=1.2]
  |
  }
\hline
  { IPM param. } & {\# instances}
  & \multicolumn{3}{c|}{HyPre $t_{\text{total}} / \nnz $}
                      & \multicolumn{3}{c|}{PETSc $t_{\text{total}} / \nnz $}
  & \multicolumn{3}{c|}{ICC $t_{\text{total}} / \nnz $} \\ \hline
{   }             &                            &{ median   }&{ 0.75    }&{ max   }&{ median  }&{ 0.75   }& {max}
&{ median  }&{ 0.75   }& {max}  \\ \hline
{ $(10^{-3})$}  &   & {$(\mu s)$} &{$(\mu s)$}  & {$(\mu s)$} & {$(\mu s)$} &{$(\mu s)$}  & {$(\mu s)$}
& {$(\mu s)$} &{$(\mu s)$}  & {$(\mu s)$} \\ \hline
100&28&1.54&2.27&5.04&N/A&N/A&N/A&0.857&1.19&2.84\\
10&27&1.52&2.29&4.82&N/A&N/A&N/A&0.873&1.18&2.82\\
1&27&1.52&2.32&4.94\tnote{*}&N/A&N/A&N/A&0.855&1.19&2.92\\
0.1&27&1.53&2.29&4.86\tnote{*}&N/A&N/A&N/A&0.858&1.2&2.88\\
0.01&19&1.52&2.3&4.58\tnote{*}&N/A&N/A&N/A&0.851&1.19&2.91\\
  \hline
\end{tabular}
} %
 \caption{
   Chimera IPM}
\label{tab:chimeraIPM}
\begin{tablenotes}
\footnotesize
\item[*] relative residual error exceeded tolerance $10^{-8}$ by a factor $(1,10^{4}]$.
\item[**] relative residual error exceeded tolerance $10^{-8}$ by a factor $(10^{4},10^{8})$.
\item[] \hspace{-1em} Inf: Solver crashed or returned a solution with relative residual error 1 (trivial solution).
\item[] \hspace{-1em} N/A: Experiment omitted as the solver crashed too often with (only occurred for PETSc).
\end{tablenotes}
\end{threeparttable}

\vspace{1em}

\begin{threeparttable}[H]
  \centering
\resizebox{\textwidth}{!}{%
\begin{tabular}{
  |S[table-format=4.2e1]
  |S[table-format=2.0]
  |S[table-format=1.3] S[table-format=1.3] S[table-format=1.3]
  |S[table-format=1.3] S[table-format=1.3] S[table-format=1.3]
  |S[table-format=1.2] S[table-format=1.2] S[table-format=1.2]
  |
  }
\hline
  { edges } & {\# instances}
  & \multicolumn{3}{c|}{AC $t_{\text{total}} / \nnz $}
  & \multicolumn{3}{c|}{{\ACSM} $t_{\text{total}} / \nnz $}
  & \multicolumn{3}{c|}{CMG $t_{\text{total}} / \nnz $} 
  \\ \hline
{  $n$ }             &                             &{ median   }&{ 0.75    }&{ max   }&{ median   }&{ 0.75    }&{ max   }&{ median  }&{ 0.75   }& {max} 
\\ \hline
{ $(K)$}  &  &{$(\mu s)$}  & {$(\mu s)$} & {$(\mu s)$} & {$(\mu s)$} &{$(\mu s)$}  & {$(\mu s)$} & {$(\mu s)$} &{$(\mu s)$}  & {$(\mu s)$} 
\\ \hline
1030&10&0.162&0.172&0.247&0.28&0.287&0.313&1.18&1.2&1.35\\
8100&10&0.194&0.195&0.223&0.296&0.305&0.369&1.03&1.03\tnote{*}&1.07\tnote{**}\\
2.72e+04&10&0.229&0.25&0.286&0.346&0.354&0.389&1.26\tnote{*}&1.37\tnote{*}&1.57\tnote{*}\\
6.44e+04&10&0.247&0.268&0.291&0.357&0.367&0.396&1.13\tnote{*}&1.35\tnote{*}&1.62\tnote{*}\\
1.26e+05&10&0.264&0.279&0.362&0.422&0.428&0.428&1.43&1.44\tnote{*}&1.69\tnote{*}\\
2.17e+05&11&0.418&0.437&0.735&0.502&0.51&0.809&1.5&1.57\tnote{*}&1.74\tnote{*}\\
  \hline
\end{tabular}
} %
\\ \vspace{0.5em}
\resizebox{\textwidth}{!}{%
\begin{tabular}{
  |S[table-format=5.2e1]
  |S[table-format=2.0]
  |S[table-format=1.3] S[table-format=1.3] S[table-format=1.3]
  |S[table-format=1.3] S[table-format=1.3] S[table-format=1.3]
  |S[table-format=1.2] S[table-format=1.2] S[table-format=1.2]
  |
  }
\hline
  { edges } & {\# instances}
  & \multicolumn{3}{c|}{HyPre $t_{\text{total}} / \nnz $}
                      & \multicolumn{3}{c|}{PETSc $t_{\text{total}} / \nnz $}
  & \multicolumn{3}{c|}{ICC $t_{\text{total}} / \nnz $} \\ \hline
{  $n$ }             &                             &{ median   }&{ 0.75    }&{ max   }&{ median   }&{ 0.75    }&{ max   }&{ median  }&{ 0.75   }& {max}  \\ \hline
{ $(K)$}  &  &{$(\mu s)$}  & {$(\mu s)$} & {$(\mu s)$} & {$(\mu s)$} &{$(\mu s)$}  & {$(\mu s)$} & {$(\mu s)$} &{$(\mu s)$}  & {$(\mu s)$} \\ \hline
1030&10&0.346&0.361&0.449&0.76&1.06&1.81\tnote{*}&3.56&3.64&3.99\\
8100&10&0.454&0.48&1.22&1.22&1.85&3.32\tnote{*}&Inf&Inf&Inf\\
2.72e+04&10&0.533&0.551&0.621&1.84&2.47&2.61&Inf&Inf&Inf\\
6.44e+04&10&0.583&0.62&0.725&1.68&2.05&3.27&Inf&Inf&Inf\\
1.26e+05&10&0.633&0.665&0.768&2.38&2.6&3.24&Inf&Inf&Inf\\
2.17e+05&11&0.643&0.693&0.811&2.31&2.82&3.98&Inf&Inf&Inf\\
  \hline
\end{tabular}
} %
 \caption{
   Maximum Flow IPM Laplacians from Spielman graphs.}
\label{tab:ipmrloweps}
\begin{tablenotes}
\footnotesize
\item[*] relative residual error exceeded tolerance $10^{-8}$ by a factor $(1,10^{4}]$.
\item[**] relative residual error exceeded tolerance $10^{-8}$ by a factor $(10^{4},10^{8})$.
\item[] \hspace{-1em} Inf: Solver crashed or returned a solution with relative residual error 1 (trivial solution).
\item[] \hspace{-1em} N/A: Experiment omitted as the solver crashed too often with (only occurred for PETSc).
\end{tablenotes}
\end{threeparttable}

%% file: poisson.tex
\newpage
\subsubsection{SPE Benchmark}
In this section, we report the performance of solvers on SDDM matrices from the Society of Petroleum Engineering benchmark~\cite{CB01}. 
Here we made use of the matrices from~\cite{CCBRTD20}. These matrices are generated from a $3$D waterflood of a geostatistical model. 
The problem size ranges from $0.5$M to $16$M.

Results are shown in Table~\ref{tab:spetable}. On this problem, CMG and HyPre have particularly good performances and beats
other solvers. On the other hand, ICC performed poorly with very large iteration counts for the two largest problem
instances. The runtimes of {\ACSM} are very close to that of AC, but again, SPE benchmark is not hard enough
for {\ACSM} to out-perform AC.

\begin{landscape}
  \newcommand{\lstablelen}{24cm}
  \begin{table}[H]	
    \centering
    \resizebox{\lstablelen}{!}{%
      \begin{tabular}{
        |S[table-format=6.0]
        |S[table-format=5.0]
        |S[table-format=2.0] S[table-format=1.3] S[table-format=1.3] S[table-format=1.3]
        |S[table-format=2.0] S[table-format=1.3] S[table-format=1.2] S[table-format=1.3]
        |S[table-format=2.0] S[table-format=1.3] S[table-format=1.3] S[table-format=1.3]
        |
        }
        \hline
        { nonzeros }&{ variables }
        &  \multicolumn{4}{c|}{AC}
        &  \multicolumn{4}{c|}{{\ACSM}}
        &  \multicolumn{4}{c|}{CMG}
        \\
        \hline
        {$\nnz$}  &  {$n$}
        &{ $N_{\text{iter}}$ }&{ $t_{\text{solve}} / \nnz$ }&{ $t_{\text{total}} / \nnz$ }&{ $\text{res. err.}$ }
        &{ $N_{\text{iter}}$ }&{ $t_{\text{solve}} / \nnz$ }&{ $t_{\text{total}} / \nnz$ }&{ $\text{res. err.}$ }
        &{ $N_{\text{iter}}$ }&{ $t_{\text{solve}} / \nnz$ }&{ $t_{\text{total}} / \nnz$ }&{ $\text{res. err.}$ }
        \\
        {$(K)$} & {$(K)$}
        &  & {$(\mu s)$} & {$(\mu s)$} & {$(10^{-8})$}
        &  & {$(\mu s)$} & {$(\mu s)$} & {$(10^{-8})$}
        &  & {$(\mu s)$} & {$(\mu s)$} & {$(10^{-8})$}
        \\  \hline %
        2910&422&17&0.24&0.718&0.433&13&0.23&1.07&0.288&13&0.222&0.436&0.504\\
14600&2100&22&0.359&0.967&0.337&14&0.291&1.31&0.315&14&0.184&0.363&0.679\\
28500&4100&22&0.404&0.928&0.474&14&0.337&1.34&0.416&15&0.194&0.364&0.884\\
55800&8000&39&0.893&1.49&0.495&27&0.796&1.79&0.256&25&0.336&0.51&0.744\\
112000&16000&41&1.03&1.62&0.331&26&0.868&1.91&0.28&25&0.338&0.511&0.861\\
        \hline
      \end{tabular}
    } %
    \resizebox{\lstablelen}{!}{%
      \begin{tabular}{
        |S[table-format=6.0]
        |S[table-format=5.0]
        |S[table-format=2.0] S[table-format=1.3] S[table-format=1.3] S[table-format=1.4]
        |S[table-format=1.0] S[table-format=1.2] S[table-format=1.2] S[table-format=1.5e-1]
        |S[table-format=3.0] S[table-format=1.3] S[table-format=1.3]S[table-format=1.3]
        |
        }
        \hline
        { nonzeros }&{ variables }
        &  \multicolumn{4}{c|}{HyPre} 
        &  \multicolumn{4}{c|}{PETSc}
        &  \multicolumn{4}{c|}{ICC}    \\
        \hline
        {$\nnz$}  &  {$n$}
        &{ $N_{\text{iter}}$ }&{ $t_{\text{solve}} / \nnz$ }&{ $t_{\text{total}} / \nnz$ }&{ $\text{res. err.}$ }
        &{ $N_{\text{iter}}$ }&{ $t_{\text{solve}} / \nnz$ }&{ $t_{\text{total}} / \nnz$ }&{ $\text{res. err.}$ }
        &{ $N_{\text{iter}}$ }&{ $t_{\text{solve}} / \nnz$ }&{ $t_{\text{total}} / \nnz$ }&{ $\text{res. err.}$ }
        \\
        {$(K)$} & {$(K)$}
        &  & {$(\mu s)$} & {$(\mu s)$} & {$(10^{-8})$}
        &  & {$(\mu s)$} & {$(\mu s)$} & {$(10^{-8})$}
        &  & {$(\mu s)$} & {$(\mu s)$} & {$(10^{-8})$}
        \\  \hline %
        2910&422&6&0.199&0.508&0.0937&9&1.63&2.38&1e-06&36&0.466&0.527&0.99\\
14600&2100&6&0.213&0.518&0.378&9&1.71&2.46&4.33e-06&57&0.588&0.683&0.999\\
28500&4100&6&0.215&0.519&0.384&9&1.73&2.49&5.78e-06&59&0.596&0.699&0.935\\
55800&8000&10&0.34&0.646&0.287&9&1.77&2.57&0.00928&380&4.18&4.29&0.999\\
112000&16000&10&0.342&0.648&0.301&9&1.78&2.58&0.00939&377&4.08&4.2&0.998\\
        \hline
      \end{tabular}
    } %
    \caption{SPE Benchmark}
    \label{tab:spetable}
  \end{table}
\end{landscape}

\subsubsection{Poisson Problems on 3D Domains}

We adopt a suite of experiments similar to \cite{SBBRGS12}, which
evaluated the performance of the HyPre-BoomerAMG algebraic multigrid
on discretizations of several Poisson problems, also known as scalar
elliptic diffusion problems.

We consider systems of linear equations arising from Poisson problems
with Dirichlet boundary conditions on a domain $\Omega$ with a
spatially varying scalar diffusion
coefficient $\mu$, and consider both the case of
$\mu$ isotropic and anisotropic,

\begin{equation}
  \begin{cases}
    \operatorname{div}(\mu(x)\grad u (x)) = f(x)
    & \text{ for all } x \in
  \Omega \\
  u(x) = 0
  & \text{ for } x \in \partial \Omega
  \end{cases}
  \label{eq:poissoncontinuous}
\end{equation}

We adopt a finite volume method discretization approach similar to
\cite{BCKFH18}.

We consider a unit cube domain $\Omega = [0,1]^3$ using a standard 3D
7-point stencil. We discretize along three dimensions with $n_1$, $n_2$, and $n_3$ variables spaced
evenly along each dimension, leading to a 3D grid with $n = n_1\cdot n_2
\cdot n_3$ variables, indexed by $u_{i,j,k}$ with $(i,j,k) \in [n_1]
\times [n_2] \times [n_3]$, at positions 
\[
\left(\frac{i-1}{n_1-1}, \frac{j-1}{n_2-1}, \frac{k-1}{n_3-1} \right) \in \Omega
\]

When a point is on the boundary $\partial \Omega$, we introduce the equation
$u_{i,j,k} = 0$ (equivalently, we remove the corresponding matrix row and column).
When a point is not on the boundary of $\Omega$, we set
\begin{align*}
&-(a_{\text{north}}+a_{\text{south}}+a_{\text{east}}+a_{\text{west}}+a_{\text{up}}+a_{\text{down}})u_{i,j,k}
  \\
  &+
  a_{\text{north}}u_{i+1,j,k}+a_{\text{south}}u_{i-1,j,k}+a_{\text{east}}u_{i,j+1,k}+a_{\text{west}}u_{i,j-1,k}+a_{\text{up}}u_{i,j,k+1}+a_{\text{down}}u_{i,j,k-1}
  =
  f_{i,j,k}
\end{align*}
where we used $a_{\text{north}}$ to refer to the coefficient $\mu(x)$
at the domain location $x \in \Omega$ halfway between the domain locations
of $u_{i,j,k}$ and $u_{i+1,j,k}$ and so on.
This discretization scheme directly leads to an SDDM matrix.

Unlike \cite{BCKFH18}, for simplicity, when our domain $\Omega$ has
regions with varying diffusion coefficient $\mu$, we ensure these regions are aligned so
that there is always a variable on the boundary between regions and
the $x \in \Omega$ giving rise to each coefficient $a_{\text{direction}}$ falls entirely
in one region. This lets us avoid averaging over different $x$ to
obtain our coefficients $a_{\text{direction}}$.

\paragraph{Uniform coefficient 3D cube.}
In this experiment, we consider a unit cube domain $\Omega =
[0,1]^3$ with uniform coefficient $\mu(x) = 1$ everywhere and equally
fine-grained discretization along each dimension, i.e. $n_1 = 
n_2 = n_3$.
Results are shown in Table~\ref{tab:uniformgridtable}. We see that {\ACSM} has
slightly better iteration counts on the uniform grid than AC, suggesting that
indeed {\ACSM} produced more reliable facorization. ICC in contrast produced
the most unreliable factorization, causing its iteration counts to be significantly more
than other solvers. HyPre and PETSc have the best iteration counts. However, even though
they both use BoomerAMG as the preconditioner, HyPre has better runtime than PETSc.

\paragraph{High-contrast variable coefficient 3D cube.}
In this experiment, we consider a unit cube domain $\Omega =
[0,1]^3$, with equally
fine-grained discretization along each dimension, i.e. $n_1 = 
n_2 = n_3$.
We now adopt a checkboard pattern for the coefficient $\mu$,
dividing this into subregions with varying coefficients, similar to
experiments in \cite{SBBRGS12,BCKFH18}.

\begin{figure}[H]
  \centering
    \fbox{\includegraphics[width=0.3\textwidth]{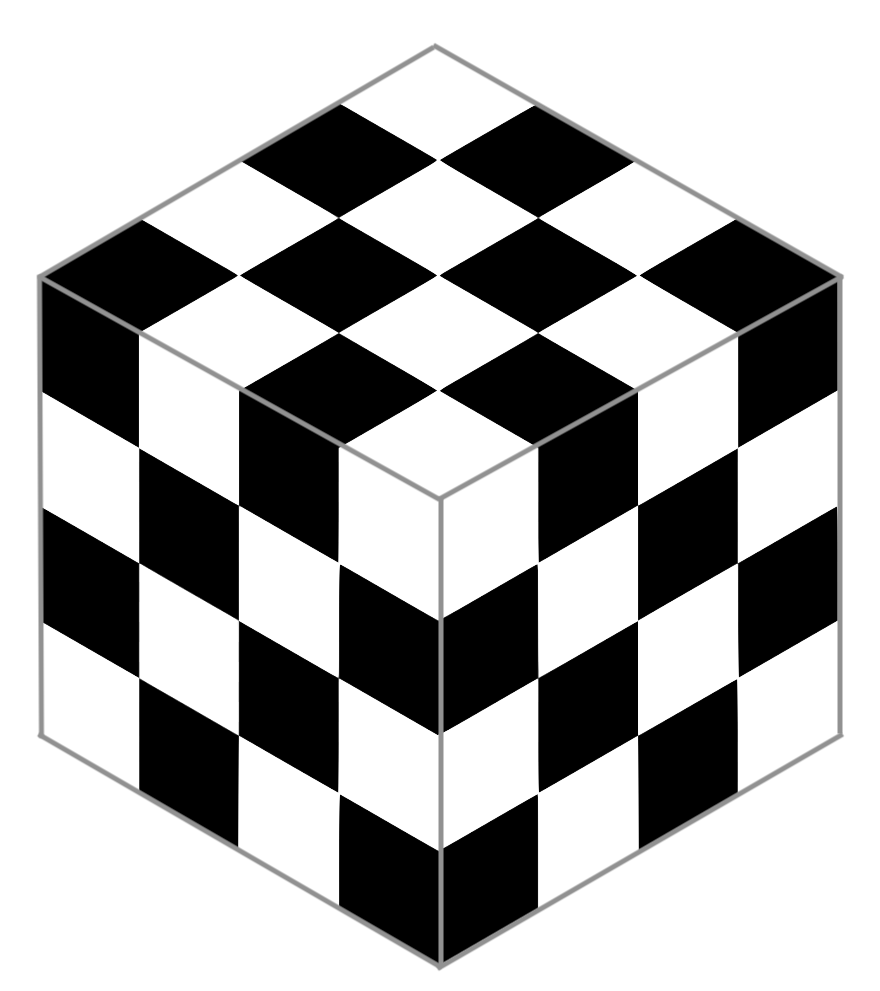}}
  \caption{Checkered cube domain $\Omega$ following
    Equation~\eqref{eq:checkeredcube} with $k = 4$. Black subregions have $\mu(x) = w$,
    white subregions have $\mu(x) = 1$.}
  \label{fig:checkeredcube}
\end{figure}

Formally, we divide the $[0,1]^3$ into $k$ intervals along each axis,
leading to $k^3$ smaller cubes, which we refer to as subregions.
Each subregion has a fixed $\mu$ value, and we pick these to form a
checkerboard pattern. Formally, the checkerboard pattern is given a 
pattern of for $(x,y,z) \in \Omega$ 
\begin{equation}
  \label{eq:checkeredcube}
  \mu(x,y,z) = 
  \begin{cases}
    1 & \text{ if } \floor{k\cdot x}+\floor{k\cdot z}+\floor{k\cdot
      y} \equiv 0 \pmod 2\\
    w & \text{ otherwise }
  \end{cases}
\end{equation}
Figure~\ref{fig:checkeredcube} shows an example with $k = 4$.

Results are shown in Table~\ref{tab:checkeredboard}. Experiments showed that high-contrast variable
coefficient 3D cube is more difficult than the uniform grid for our solvers. But it is not 
difficult enough for {\ACSM} to match AC in terms of the runtime. Again, ICC produced
the most unreliable factorization and consequently the worst iteration count. CMG and HyPre on the other
hand maintained roughly the same level of iteration count and runtime
as for the uniform grid.

\begin{landscape}%
  \newcommand{\lstablelen}{24cm}
  
  \begin{table}[H]	
    \centering
    \resizebox{\lstablelen}{!}{%
      \begin{tabular}{
        |S[table-format=6.0]
        |S[table-format=5.0]
        |S[table-format=2.0] S[table-format=1.3] S[table-format=1.3] S[table-format=1.3]
        |S[table-format=2.0] S[table-format=1.3] S[table-format=1.2] S[table-format=1.3]
        |S[table-format=2.0] S[table-format=1.3] S[table-format=1.3] S[table-format=1.3]
        |
        }
        \hline
        { nonzeros }&{ variables }
        &  \multicolumn{4}{c|}{AC}
        &  \multicolumn{4}{c|}{{\ACSM}}
        &  \multicolumn{4}{c|}{CMG}
        \\
        \hline
        {$\nnz$}  &  {$n$}
        &{ $N_{\text{iter}}$ }&{ $t_{\text{solve}} / \nnz$ }&{ $t_{\text{total}} / \nnz$ }&{ $\text{res. err.}$ }
                  &{ $N_{\text{iter}}$ }&{ $t_{\text{solve}} / \nnz$ }&{ $t_{\text{total}} / \nnz$ }&{ $\text{res. err.}$ }
                 &{ $N_{\text{iter}}$ }&{ $t_{\text{solve}} / \nnz$ }&{ $t_{\text{total}} / \nnz$ }&{ $\text{res. err.}$ }
                                                                           
        \\
        {$(K)$} & {$(K)$}
        &  & {$(\mu s)$} & {$(\mu s)$} & {$(10^{-8})$}
                  & & {$(\mu s)$} & {$(\mu s)$} & {$(10^{-8})$}
                & & {$(\mu s)$}& {$(\mu s)$} & {$(10^{-8})$}
                                                                    
        \\  \hline %
        1990&287&24&0.319&0.769&0.524&18&0.32&1.17&0.739&27&0.422&0.602&0.815\\
19900&2860&25&0.453&1.07&0.913&20&0.549&1.74&0.472&27&0.423&0.537&0.773\\
200000&28700&27&0.833&1.46&0.612&20&0.99&2.49&0.679&27&0.509&0.624&0.982\\
        \hline
      \end{tabular}
    } %
    \resizebox{\lstablelen}{!}{%
      \begin{tabular}{
        |S[table-format=6.0]
        |S[table-format=5.0]
        |S[table-format=1.0] S[table-format=1.3] S[table-format=1.3] S[table-format=1.3]
        |S[table-format=1.0] S[table-format=1.3] S[table-format=1.2] S[table-format=1.6]
        |S[table-format=3.0] S[table-format=1.3] S[table-format=1.3]S[table-format=1.3]
        |
        }
        \hline
        { nonzeros }&{ variables }
        &  \multicolumn{4}{c|}{HyPre}
        &  \multicolumn{4}{c|}{PETSc}
        &  \multicolumn{4}{c|}{ICC}    \\
        \hline
        {$\nnz$}  &  {$n$}
        &{ $N_{\text{iter}}$ }&{ $t_{\text{solve}} / \nnz$ }&{ $t_{\text{total}} / \nnz$ }&{ $\text{res. err.}$ }
                  &{ $N_{\text{iter}}$ }&{ $t_{\text{solve}} / \nnz$ }&{ $t_{\text{total}} / \nnz$ }&{ $\text{res. err.}$ }
                 &{ $N_{\text{iter}}$ }&{ $t_{\text{solve}} / \nnz$ }&{ $t_{\text{total}} / \nnz$ }&{ $\text{res. err.}$ }
                                                                           
        \\
        {$(K)$} & {$(K)$}
        &  & {$(\mu s)$} & {$(\mu s)$} & {$(10^{-8})$}
                  & & {$(\mu s)$} & {$(\mu s)$} & {$(10^{-8})$}
                & & {$(\mu s)$}& {$(\mu s)$} & {$(10^{-8})$}
        \\  \hline %
        1990&287&7&0.201&0.488&0.197&5&0.62&1.62&0.00402&61&0.694&0.765&0.917\\
19900&2860&8&0.259&0.572&0.113&5&0.771&2.02&0.00518&109&1.09&1.19&0.974\\
200000&28700&7&0.235&0.587&0.206&6&1.01&2.62&0.000333&187&2.08&2.22&0.973\\
        \hline
      \end{tabular}
    } %
    \caption{Uniform coefficient 3D cube}
    \label{tab:uniformgridtable}
  \end{table}

  \begin{table}[H]
    \centering
    \resizebox{\lstablelen}{!}{%
      \begin{tabular}{
        |S[table-format=3.0]
        |S[table-format=2.0] S[table-format=1.2] S[table-format=1.2] S[table-format=1.3]
        |S[table-format=2.0] S[table-format=1.2] S[table-format=1.2] S[table-format=1.3]
        |S[table-format=2.0] S[table-format=1.3] S[table-format=1.3] S[table-format=1.3]
        |
        }
        \hline    
        {axis intervals}
        &  \multicolumn{4}{c|}{AC}
        &  \multicolumn{4}{c|}{{\ACSM}}
        &  \multicolumn{4}{c|}{CMG}
   \\
        \hline
        {$k$}
        &{ $N_{\text{iter}}$ }&{ $t_{\text{solve}} / \nnz$ }&{ $t_{\text{total} / \nnz}$ }&{ $\text{res. err.}$ }
        &{ $N_{\text{iter}}$ }&{ $t_{\text{solve}} / \nnz$ }&{ $t_{\text{total} / \nnz}$ }&{ $\text{res. err.}$ }
        &{ $N_{\text{iter}}$ }&{ $t_{\text{solve}} / \nnz$ }&{ $t_{\text{total} / \nnz}$ }&{ $\text{res. err.}$ }
        \\ \hline %
        &  & {$(\mu s)$} & {$(\mu s)$} & {$(10^{-8})$}
        &  & {$(\mu s)$} & {$(\mu s)$} & {$(10^{-8})$}
        &  & {$(\mu s)$} & {$(\mu s)$} & {$(10^{-8})$}
        \\ \hline %
        2&41&1.28&1.88&0.742&30&1.48&3.01&0.987&28&0.422&0.539&0.588\\
4&47&1.47&2.1&0.733&35&1.74&3.14&0.786&28&0.436&0.553&0.781\\
8&52&1.61&2.24&0.836&38&1.87&3.26&0.985&27&0.425&0.551&0.745\\
16&52&1.57&2.19&0.915&46&2.2&3.57&0.707&28&0.446&0.572&0.834\\
32&59&1.73&2.33&0.802&45&2.06&3.38&0.967&28&0.399&0.534&0.867\\
64&57&1.59&2.16&0.925&48&2.02&3.23&0.762&23&0.383&0.511&0.952\\
128&53&1.39&1.94&0.98&35&1.31&2.46&0.828&26&0.403&0.541&0.907\\
        \hline
      \end{tabular}
} %
        \resizebox{\lstablelen}{!}{%
      \begin{tabular}{
        |S[table-format=3.0]
        |S[table-format=1.0] S[table-format=1.3] S[table-format=1.3] S[table-format=1.3]
        |S[table-format=1.0] S[table-format=1.3] S[table-format=1.3] S[table-format=1.6]
        |S[table-format=3.0] S[table-format=1.2] S[table-format=1.2] S[table-format=1.3]
        |
        }
        \hline    
        {axis intervals}
        &  \multicolumn{4}{c|}{HyPre}
        &  \multicolumn{4}{c|}{PETSc}
        &  \multicolumn{4}{c|}{ICC}
        \\
        \hline
       {k}
        &{ $N_{\text{iter}}$ }&{ $t_{\text{solve}} / \nnz$ }&{ $t_{\text{total} / \nnz}$ }&{ $\text{res. err.}$ }
        &{ $N_{\text{iter}}$ }&{ $t_{\text{solve}} / \nnz$ }&{ $t_{\text{total} / \nnz}$ }&{ $\text{res. err.}$ }
        &{ $N_{\text{iter}}$ }&{ $t_{\text{solve}} / \nnz$ }&{ $t_{\text{total} / \nnz}$ }&{ $\text{res. err.}$ }
        \\ \hline %
        &  & {$(\mu s)$} & {$(\mu s)$} & {$(10^{-8})$}
        &  & {$(\mu s)$} & {$(\mu s)$} & {$(10^{-8})$}
        &  & {$(\mu s)$} & {$(\mu s)$} & {$(10^{-8})$}
        \\ \hline %
        2&8&0.265&0.616&0.308&7&1.35&3.12&0.000145&175&1.94&2.09&0.968\\
4&8&0.264&0.656&0.738&7&1.65&3.79&0.000555&244&2.69&2.84&0.985\\
8&9&0.295&0.734&0.351&7&1.77&3.85&0.00142&323&3.54&3.68&0.982\\
16&9&0.289&0.752&0.184&8&1.95&3.84&0.000297&351&3.84&3.99&0.979\\
32&8&0.257&0.732&0.744&8&1.83&3.43&0.00117&343&3.77&3.92&0.952\\
64&8&0.258&0.64&0.513&9&1.43&2.35&0.000498&383&4.19&4.33&0.958\\
128&7&0.244&0.55&0.98&7&0.664&1.1&0.00106&420&4.61&4.76&0.994\\
        \hline
      \end{tabular}
} %
    
    \caption{High-contrast variable coefficient 3D cube.
      Unit cube domain with  $\Omega =[0,1]^3$, with equally
      fine-grained discretization along each dimension.
      28.7 M variables and 200 M non-zeros.
      Weight $w = 10^7$.
      Times are reported in seconds.
    }
    \label{tab:checkeredboard}
  \end{table}

\begin{table}[H]
    \centering
  \resizebox{\lstablelen}{!}{%
  \begin{tabular}{
    |S[table-format=4.3]
    |S[table-format=3.0]
    |S[table-format=2.1]
    |S[table-format=2.0] S[table-format=1.3] S[table-format=1.3] S[table-format=1.3]
    |S[table-format=2.0] S[table-format=1.3] S[table-format=1.2] S[table-format=1.3]
    |S[table-format=2.0] S[table-format=1.3] S[table-format=1.3] S[table-format=1.3]
    |
    }
\hline
    { aniso. stretch }& { nonzeros }&{ variables }&
                                                    \multicolumn{4}{c|}{AC}     & \multicolumn{4}{c|}{{\ACSM}}     
                                                    & \multicolumn{4}{c|}{CMG} %
                                                    \\
\hline
    { $\eta$ } & {$\nnz$}  &  {$n$} 
        &{ $N_{\text{iter}}$ }&{ $t_{\text{solve}} / \nnz$ }&{ $t_{\text{total}} / \nnz$ }&{ $\text{res. err.}$ }
        &{ $N_{\text{iter}}$ }&{ $t_{\text{solve}} / \nnz$ }&{ $t_{\text{total}} / \nnz$ }&{ $\text{res. err.}$ }
        &{ $N_{\text{iter}}$ }&{ $t_{\text{solve}} / \nnz$ }&{ $t_{\text{total}} / \nnz$ }&{ $\text{res. err.}$ }
    \\ \hline %
        & {$(M)$} & {$(M)$}  
        & & {$(\mu s)$} & {$(\mu s)$} & {$(10^{-8})$}
        & & {$(\mu s)$} & {$(\mu s)$} & {$(10^{-8})$}
        & & {$(\mu s)$} & {$(\mu s)$} & {$(10^{-8})$}
    \\
    \hline
    0.001&203&29.6&26&0.637&1.16&0.661&19&0.729&1.99&0.871&25&0.39&0.512&0.929\\
0.01&200&28.8&26&0.724&1.32&0.904&20&0.864&2.26&0.592&26&0.521&0.644&0.655\\
0.1&200&28.7&27&0.809&1.42&0.618&20&0.956&2.42&0.681&26&0.543&0.659&0.712\\
1&200&28.7&26&0.811&1.44&0.906&20&0.994&2.5&0.753&27&0.508&0.624&0.982\\
10&200&28.7&27&0.858&1.48&0.635&20&0.991&2.47&0.61&28&0.513&0.629&0.708\\
100&196&28.6&25&0.741&1.31&0.845&19&0.837&2.12&0.842&28&0.516&0.637&0.635\\
1000&189&29.9&23&0.494&0.956&0.733&18&0.502&1.39&0.454&15&0.283&0.416&0.333\\
    \hline
  \end{tabular}
} %
  \resizebox{\lstablelen}{!}{%
  \begin{tabular}{
    |S[table-format=4.3]
    |S[table-format=3.0]
    |S[table-format=2.1]
    |S[table-format=1.0] S[table-format=1.3] S[table-format=1.3] S[table-format=1.4]
    |S[table-format=1.0] S[table-format=1.3] S[table-format=1.2] S[table-format=1.6]
    |S[table-format=3.0] S[table-format=1.3] S[table-format=1.3] S[table-format=1.3]
    |
    }
\hline
    { aniso. stretch }& { nonzeros }&{ variables }
    & \multicolumn{4}{c|}{HyPre}
    &\multicolumn{4}{c|}{PETSc}     
                                                    & \multicolumn{4}{c|}{ICC}        
                                                    \\
\hline
    { $\eta$ } & {$\nnz$}  &  {$n$} 
        &{ $N_{\text{iter}}$ }&{ $t_{\text{solve}} / \nnz$ }&{ $t_{\text{total}} / \nnz$ }&{ $\text{res. err.}$ }
        &{ $N_{\text{iter}}$ }&{ $t_{\text{solve}} / \nnz$ }&{ $t_{\text{total}} / \nnz$ }&{ $\text{res. err.}$ }
        &{ $N_{\text{iter}}$ }&{ $t_{\text{solve}} / \nnz$ }&{ $t_{\text{total}} / \nnz$ }&{ $\text{res. err.}$ }
    \\ \hline %
        & {$(M)$} & {$(M)$}  
        & & {$(\mu s)$} & {$(\mu s)$} & {$(10^{-8})$}
        & & {$(\mu s)$} & {$(\mu s)$} & {$(10^{-8})$}
        & & {$(\mu s)$} & {$(\mu s)$} & {$(10^{-8})$}
    \\
    \hline
    0.001&203&29.6&7&0.222&0.589&0.192&5&0.677&1.62&0.00366&47&0.531&0.603&0.955\\
0.01&200&28.8&7&0.228&0.583&0.2&5&0.732&1.84&0.00456&81&0.892&0.985&0.959\\
0.1&200&28.7&7&0.234&0.581&0.204&5&0.808&2.15&0.00618&129&1.43&1.54&0.999\\
1&200&28.7&7&0.235&0.587&0.205&6&1.01&2.6&0.000371&160&1.77&1.92&0.987\\
10&200&28.7&8&0.261&0.614&0.103&5&0.638&1.77&0.00414&101&1.13&1.27&0.943\\
100&196&28.6&7&0.223&0.574&0.187&4&0.494&1.52&0.116&35&0.401&0.53&0.752\\
1000&189&29.9&7&0.194&0.476&0.0744&4&0.545&1.41&0.00211&12&0.153&0.285&0.651\\
    \hline
  \end{tabular}
} %
  \caption{Anisotropic coefficient 3D cube with variable discretization and fixed weight.}
  \label{tab:anisogridtable}
\end{table}

\begin{table}[H]
  \centering
\resizebox{\lstablelen}{!}{%
\begin{tabular}{
  |S[table-format=4.3]
  |S[table-format=2.0] S[table-format=1.3] S[table-format=1.3] S[table-format=1.3]
  |S[table-format=2.0] S[table-format=1.3] S[table-format=1.2] S[table-format=1.3]
  |S[table-format=2.0] S[table-format=1.3] S[table-format=1.3] S[table-format=1.3]
  |
  }
\hline
  { weight }& 
                                                  \multicolumn{4}{c|}{AC}    
                                   & \multicolumn{4}{c|}{{\ACSM}}     
                                   & \multicolumn{4}{c|}{CMG} 
                                    \\
\hline
  { $\eta$ } 
      &{ $N_{\text{iter}}$ }&{ $t_{\text{solve}} / \nnz$ }&{ $t_{\text{total}} / \nnz$ }&{ $\text{res. err.}$ }
      &{ $N_{\text{iter}}$ }&{ $t_{\text{solve}} / \nnz$ }&{ $t_{\text{total}} / \nnz$ }&{ $\text{res. err.}$ }
      &{ $N_{\text{iter}}$ }&{ $t_{\text{solve}} / \nnz$ }&{ $t_{\text{total}} / \nnz$ }&{ $\text{res. err.}$ }
  \\ \hline %
      & & {$(\mu s)$} & {$(\mu s)$} & {$(10^{-8})$}
      & & {$(\mu s)$} & {$(\mu s)$} & {$(10^{-8})$}
      & & {$(\mu s)$} & {$(\mu s)$} & {$(10^{-8})$}
  \\
  \hline
    0.001&39&1.09&1.67&0.902&25&1.09&2.48&0.657&24&0.486&0.601&0.952\\
0.01&33&0.951&1.53&0.944&23&0.997&2.41&0.963&24&0.492&0.608&0.775\\
0.1&26&0.778&1.37&0.992&20&0.925&2.39&0.645&24&0.492&0.608&0.776\\
1&26&0.814&1.52&0.952&20&0.991&2.5&0.846&27&0.503&0.619&0.972\\
10&20&0.571&1.21&0.818&15&0.63&1.95&0.413&25&0.512&0.638&0.757\\
100&17&0.399&0.916&0.464&12&0.378&1.51&0.612&27&0.551&0.674&0.893\\
1000&14&0.302&0.78&0.814&10&0.257&1.05&0.91&27&0.507&0.625&0.95\\

  \hline
\end{tabular}
} %
\resizebox{\lstablelen}{!}{%
\begin{tabular}{
  |S[table-format=4.3]
  |S[table-format=1.0] S[table-format=1.3] S[table-format=1.3] S[table-format=1.3]
  |S[table-format=1.0] S[table-format=1.3] S[table-format=1.3] S[table-format=1.6]
  |S[table-format=3.0] S[table-format=1.3] S[table-format=1.3] S[table-format=1.3]
  |
  }
\hline
  { weight }
  & \multicolumn{4}{c|}{HyPre}      
  & \multicolumn{4}{c|}{PETSc}     
                                                  & \multicolumn{4}{c|}{ICC}        
                                                  \\
\hline
  { $\eta$ } 
       &{ $N_{\text{iter}}$ }&{ $t_{\text{solve}} / \nnz$ }&{ $t_{\text{total}} / \nnz$ }&{ $\text{res. err.}$ }
      &{ $N_{\text{iter}}$ }&{ $t_{\text{solve}} / \nnz$ }&{ $t_{\text{total}} / \nnz$ }&{ $\text{res. err.}$ }
      &{ $N_{\text{iter}}$ }&{ $t_{\text{solve}} / \nnz$ }&{ $t_{\text{total}} / \nnz$ }&{ $\text{res. err.}$ }
  \\ \hline %
      & & {$(\mu s)$} & {$(\mu s)$} & {$(10^{-8})$}
      & & {$(\mu s)$} & {$(\mu s)$} & {$(10^{-8})$}
      & & {$(\mu s)$} & {$(\mu s)$} & {$(10^{-8})$}
  \\
  \hline
    0.001&7&0.297&0.662&0.165&6&0.621&1.16&0.000816&207&2.29&2.44&0.981\\
0.01&7&0.299&0.683&0.159&6&0.733&1.38&0.000871&189&2.09&2.24&0.967\\
0.1&8&0.321&0.772&0.134&7&2.76&4.92&0.00865&168&1.86&2.01&0.969\\
1&7&0.235&0.588&0.207&5&0.861&2.46&0.00581&174&1.93&2.08&0.995\\
10&7&0.323&0.625&0.618&5&1.09&2.56&0.0704&149&1.65&1.8&0.966\\
100&8&0.282&0.449&0.34&5&0.629&1.05&0.0317&77&0.867&1.02&0.928\\
1000&7&0.223&0.332&0.128&6&0.644&0.904&0.000392&30&0.346&0.495&0.892\\
  \hline
\end{tabular}
} %
\caption{
Anisotropic coefficient 3D cube with fixed discretization and variable weight. 28.7 M variables and 200 M non-zeros}
\label{tab:wgridtable}
\end{table}

\end{landscape}

\paragraph{Anisotropic coefficient 3D cube
  with variable discretization and fixed weight.}
In this experiment, we follow \cite{SBBRGS12}, and study anisotropic
coefficient diffusion problems on a unit cube domain $\Omega =
[0,1]^3$.
We include this experiment as \cite{SBBRGS12} highlighted it as a
problematic case for Algebraic Multridgrid methods.

We consider systems of linear equations arising from Poisson problems
with Dirichlet boundary conditions on a domain $\Omega$ with a
fixed \emph{matrix} coefficient $\mu$.
We choose $\mu \in R^{3 \times 3}$ as a diagonal matrix with $\mu(1,1)
  = w$ and $\mu(2,2) = \mu(3,3) = 1$.
\begin{equation}
  \begin{cases}
    \operatorname{div}(\mu \grad u (x)) = f(x)
    & \text{ for all } x \in
  \Omega \\
  u(x) = 0
  & \text{ for } x \in \partial \Omega
  \end{cases}
  \label{eq:poissoncontinuousaniso}
\end{equation}
This is equivalent to the continuous uniform scalar coefficient problem,
studied on a rectangular domain $\Omega =
[0,1/w] \times [0,1]^2$ with equally fine-grained discretization along
each dimension.

Our discrete matrix is a grid non-zero structure with uniform weights,
symmetric diagonally dominant with negative off-diagonals
(SDDM), but we vary the discretization level along each dimension.
Concretely, we stretch or shrink the grid along the first dimension by a factor
$\eta = 1/w$ so that the number of variables along each, denoted $n_1, n_2, n_3$
satisfy $n_1 = \eta n_2 = \eta n_3$.

Results are shown in Table~\ref{tab:anisogridtable}. Experiments suggest that
this is an easier problem than the high-contrast variable coefficient
3D cube for many of the solvers.
The performances of the solvers are roughly comparable to their performances 
on the uniform grid for different discretizations, with the noticeable exception 
of ICC and PETSc, which performed significantly better for grids with more extreme stretch.

\paragraph{Anisotropic coefficient 3D cube
  with fixed discretization and variable weight.}
We again study the continuous problem from
Equation~\eqref{eq:poissoncontinuousaniso}, as in the previous experiment.
But in this experiment, we discretize equally fine-grained along each
axis.
This results in a discrete 3D grid, with weights along the
first axis of $w$ and weights along the second two axes of $1$, and
all three axes having the same number of variables, denoted $n_1, n_2,
n_3$ and satisfying $n_1 = n_2 =n_3$.

Results are shown in Table~\ref{tab:wgridtable}. For this problem, both AC, {\ACSM} and ICC achieved 
noticeably lower iteration counts and runtime for weights 100 and 1000, but not for their reciprocal.
CMG, HyPre and PETSc again maintained roughly similar performances as for uniform grids.

%% file: scaling.tex
\subsection{Running Time Scaling with Problem Size}
\label{sec:scaling}

In this section, we summarize the running time scaling with problem size of AC and AC2 across all our experiments.
We measure total time including factorization and solve.

\begin{figure}[H]
  \centering
\includegraphics[width=0.7\textwidth]{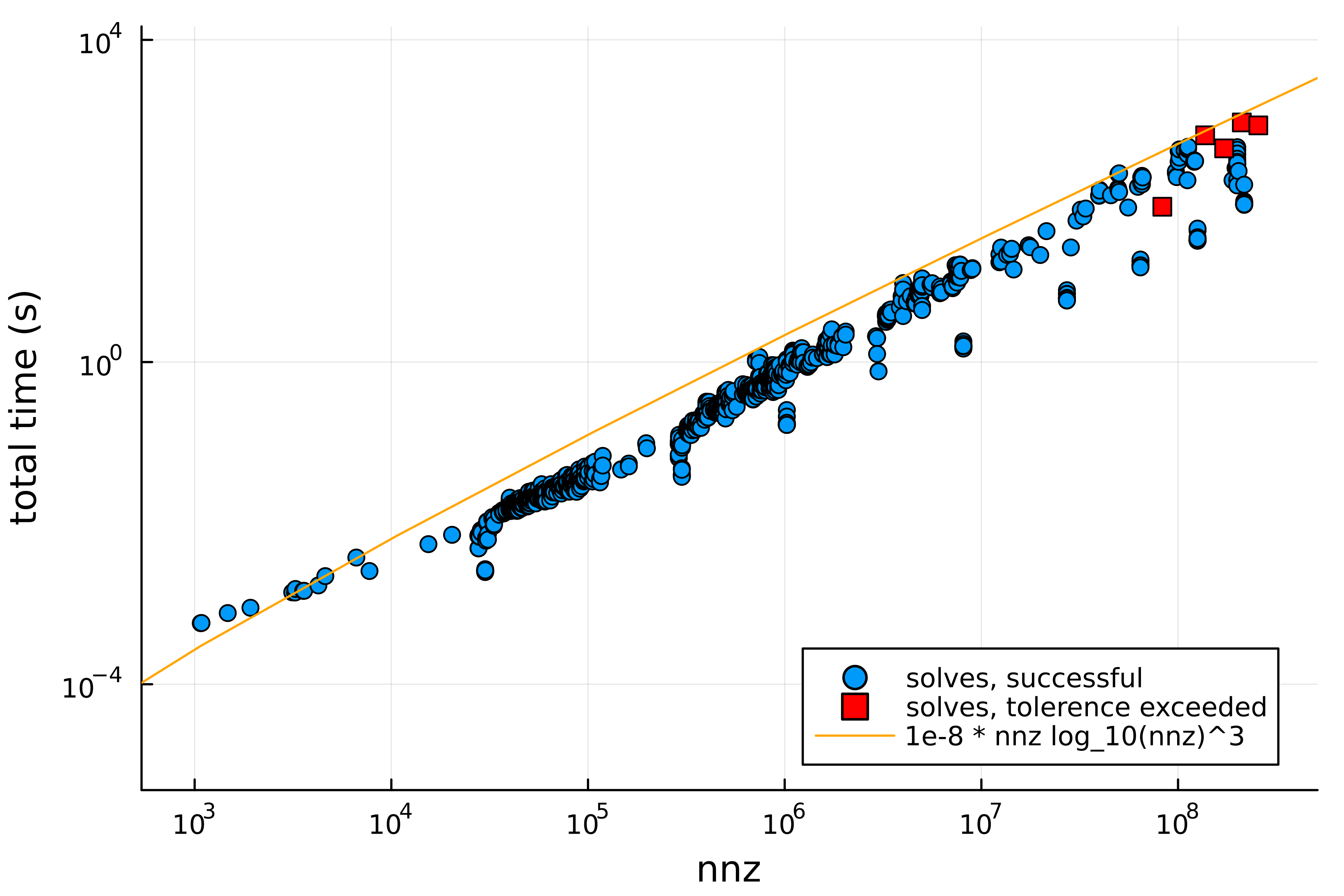}
  \caption{
  The plot shows data from all our experiments with the AC with 
  total running time in seconds on the y-axis plotted against 
  the number-of-non-zeros of the instance on the x-axis.
  Solves that reached the target tolerance are shown with a blue circle, while solves that failed to reach tolerance are shown with a red square. 
  An orange line shows the function $\nnz \log_{10}(\nnz)^3$ (seconds) where $\nnz$ is the non-zero count.
  }
  \label{fig:ac-nnz-vs-time}
\end{figure}

\begin{figure}[H]
  \centering
\includegraphics[width=0.7\textwidth]{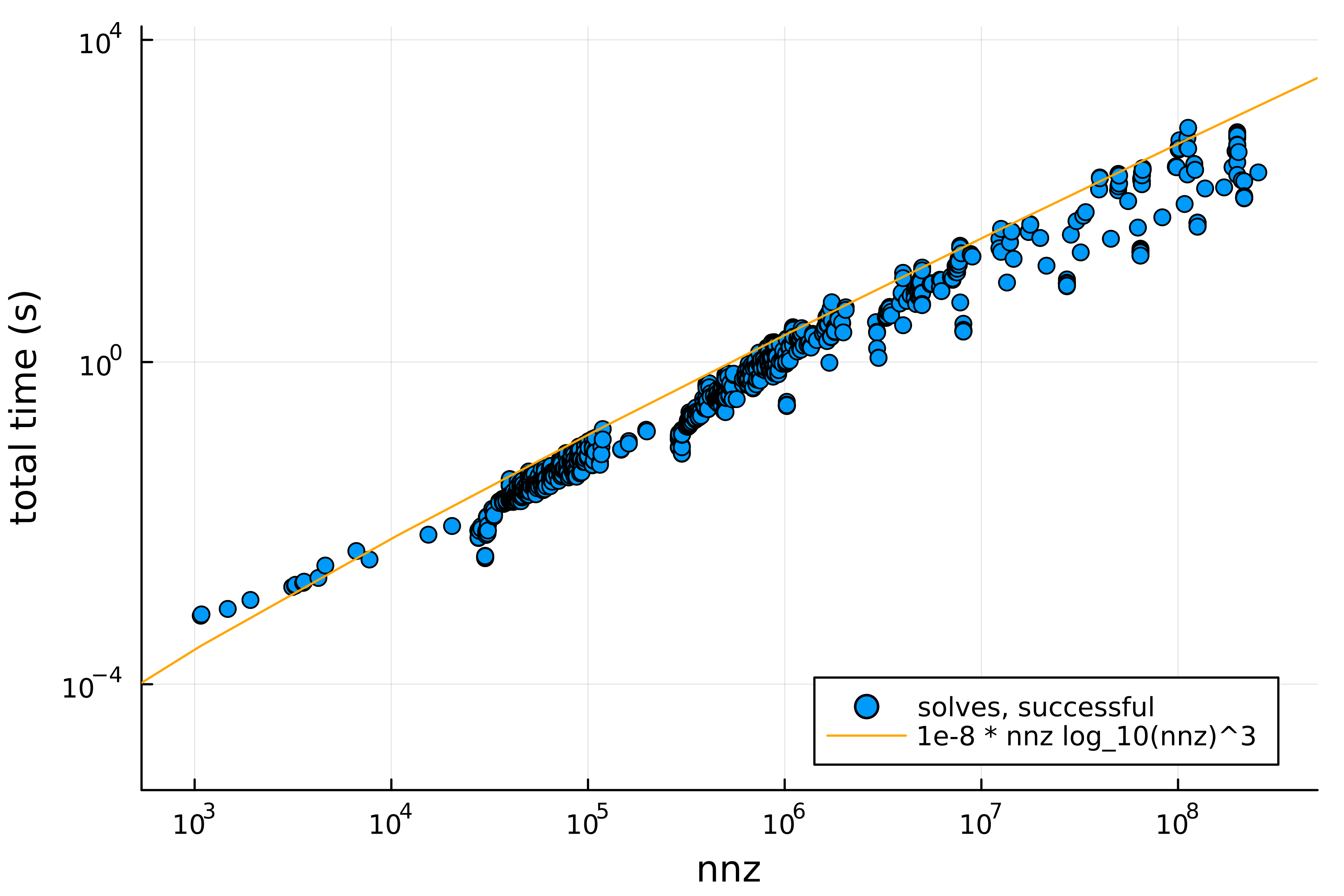}
  \caption{
  The plot shows data from all our experiments with the AC2 with 
  total running time in seconds on the y-axis plotted against 
  the number-of-non-zeros of the instance on the x-axis.
  Solves are shown with a blue circle. All solves reached the target tolerance.
  An orange line shows the function $\nnz \log_{10}(\nnz)^3$ (seconds) where $\nnz$ is the non-zero count.
  }
  \label{fig:ac2-nnz-vs-time}
\end{figure}

We show two different types of plots.
Our first plot type displays the total running time in seconds plotted against the number of non-zeros ($\nnz$) of the instances, both on a log scale. 
We also plot the function ${10^{-8} \nnz \log_{10}(\nnz)^3}$, which agrees with the asymptotic running time scaling predicted by theory from \cite{KS16}.
We caution that we do not have a theoretical analysis that establishes such a running time for AC or AC2.
Figure~\ref{fig:ac-nnz-vs-time} shows this plot for AC, while Figure~\ref{fig:ac2-nnz-vs-time} shows this plot for AC2.

\begin{figure}[H]
  \centering
\includegraphics[width=0.7\textwidth]{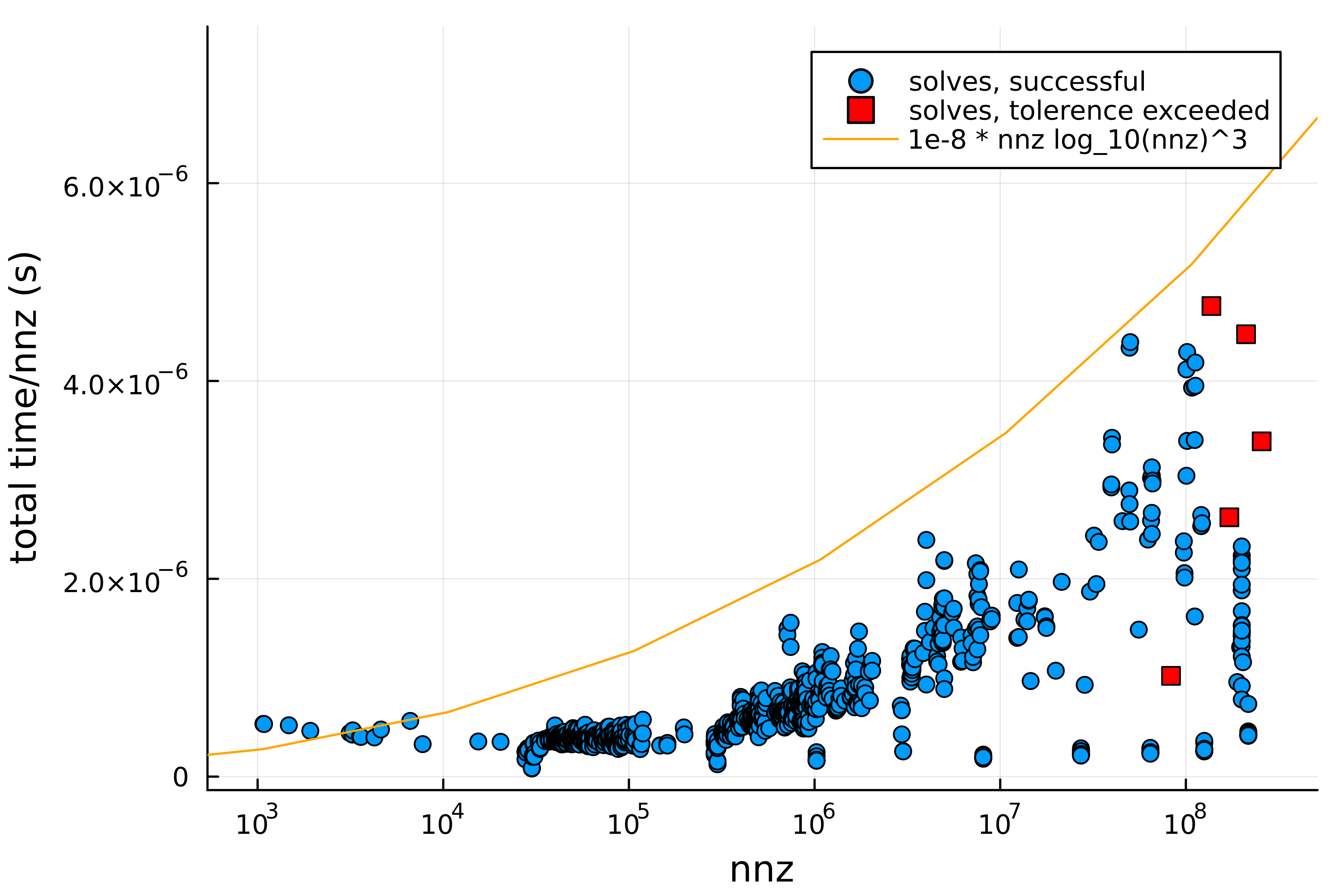}
  \caption{
  The plot shows data from all our experiments with the AC with 
  total running time \emph{per non-zero} in seconds on the y-axis plotted against 
  the number-of-non-zeros of the instance on the x-axis.
  Solves that reached the target tolerance are shown with a blue circle, while solves that failed to reach tolerance are shown with a red square. 
  An orange line shows the function $\nnz \log_{10}(\nnz)^3$ (seconds) where $\nnz$ is the non-zero count, and the same normalization is used as for the data (i.e. dividing by $\nnz$).
  }
  \label{fig:ac-nnz-vs-timebynnz}
\end{figure}

\begin{figure}[H]
  \centering
\includegraphics[width=0.7\textwidth]{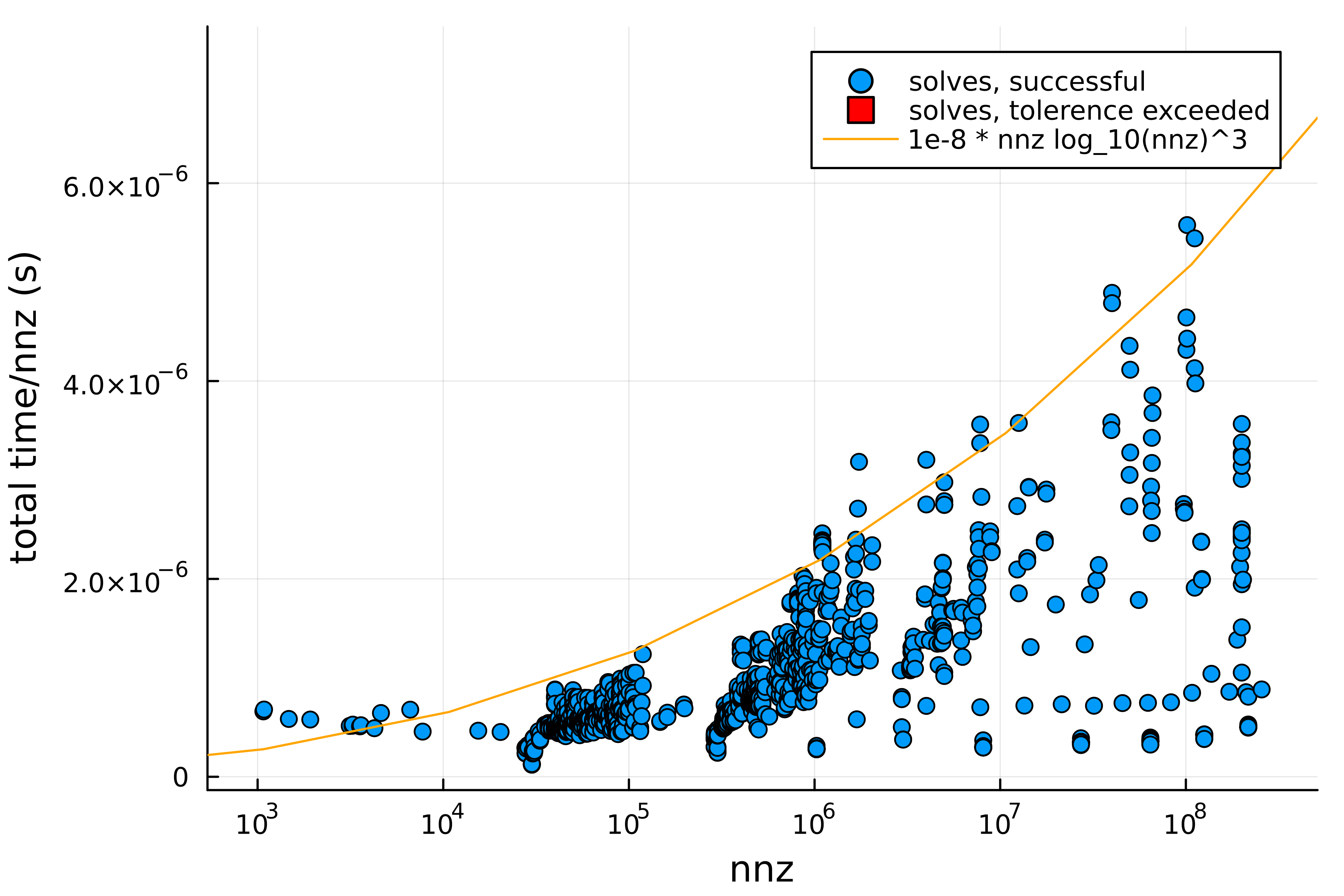}
  \caption{
  The plot shows data from all our experiments with the AC2 with 
  total running time \emph{per non-zero} in seconds on the y-axis plotted against 
  the number-of-non-zeros of the instance on the x-axis.
  Solves are shown with a blue circle. All solves reached the target tolerance.
  An orange line shows the function $\nnz \log_{10}(\nnz)^3$ (seconds) where $\nnz$ is the non-zero count, and the same normalization is used as for the data (i.e. dividing by $\nnz$).
  }
  \label{fig:ac2-nnz-vs-timebynnz}
\end{figure}

Our second plot type again shows running time in seconds \emph{measured relative to the non-zero count of the instance} plotted against non-zero count of the instances, with the x-axis on a log scale.
We also show the function $\frac{{10^{-8} \nnz \log_{10}(\nnz)^3}}{\nnz} = 10^{-8} \log_{10}(\nnz)^3$, i.e. the same running time fit as in the previous plots, but here normalized by non-zero count.
Figures~\ref{fig:ac-nnz-vs-timebynnz} and Figure~\ref{fig:ac2-nnz-vs-timebynnz} 
show these plots for AC and AC2 respectively.

\paragraph{AC running times.}
We first show the running time plots for AC as described in the previous section. 
We use a blue circle to show successful solves and a red square to show solves that failed to reach the target tolerance. 
All failed solves occurred on \emph{Sachdeva star} instances.

\paragraph{AC2 running times.} We now show running time plots for AC2 as described earlier. 
We use a blue circle to show successful solves. For AC2, all solves reached the tolerance (i.e. $10^{-8}$ relative residual error).

%% file: appendixTheory.tex
\section{Linear Algebra and Cholesky Background}
\label{sec:cholBackground}
We introduce some basic notation and definitions from linear algebra.

\paragraph{Moore-Penrose pseudo-inverse.} We use $\BB^\pinv$ to denote the Moore-Penrose pseudo-inverse of a
matrix $\BB$.
Let $\vecone \in \R^{n}$ denote the all ones vector, with dimension $n$
that will always be made clear in the context of its use.
Similarly, we let $\veczero$ denote the all zero vector or matrix,
depending on context.

The following fact is useful, since we often need to apply the
pseodo-inverse of a matrix.
\begin{restatable}[Pseudo-inverse of a product]{fact}{pinvproduct}
\label{fac:pinvproduct}
Suppose $\MM = \AA \BB \CC$ is square real matrix, where $\AA$ and
$\CC$ are non-singular.
Then
\[
\MM^{\pinv} = \proj_{\MM} \CC^{-1}\BB^{\pinv}\AA^{-1} \proj_{\MM^{\trp}}.
\]
\end{restatable}

\paragraph{LU-decomposition and Cholesky factorization of singular matrices.}
An LU-decomposition of a square matrix $\MM \in \R^{n \times n}$ is a
factorization $\MM = \matlow \matup$, where $\matlow$ is a
lower-triangular matrix and $\matup$ is an upper-triangular matrix,
both up to a permutation.

When the matrix $\MM$ is symmetrix and positive semi-definite,
study a special case of LU-decomposition, where $\matup =
\matlow^{\trp}$, i.e. $\MM = \matlow \matlow^{\trp}$.
This is known as a Cholesky Factorization.

When $\MM$ is non-singular, the linear equations $\matlow \yy = \bb$ and $\matup \xx = \yy$ and be
solved by forward and backward substitution algorithms respectively
(e.g. see \cite{TBI97}),
which run in time $O( \nnz(\matlow) )$ and $O(\nnz(\matup))$.
I.e. $\matlow^{-1}$ and $\matup^{-1}$ can be applied in time
  proportional to the number of non-zeros in $\matlow$ and $\matup$ respectively.
This means that if a decomposition $\matlow, \matup$ of $\MM$ is
known, then linear systems in $\MM$ be solved in time
$O(\nnz(\matlow)+\nnz(\matup))$, since $\MM \xx = \bb$
implies $\xx = \matup^{-1} \matlow^{-1} \bb$.

When $\MM$ is singular the same
forward and backward substitution algorithms can be used to compute
the pseudo-inverse, in some situations.
We focus on a special case where we can instead factor $\MM$ as
 $\MM = \matlow' \matdiag \matup'$,
where $\matdiag$ is singular and $\matlow'$, $\matup'$ are non-singular.

For the case of interest to us, Laplacians of connected
graphs, this slightly modified
factorization can be trivially obtained from an LU-decomposition, by
$\matlow'$ equal to $\matlow$ except $\matlow(n,n) = 1$
and
 $\matup'$ equal to $\matup$ except $\matup(n,n) = 1$
and taking $\matdiag$ to be the identity matrix, except $\matdiag(n,n)
= 0$.
For our approximate factorizations of Laplacians, the matrix $\MM =
\matlow' \matdiag \matup'$ is symmetric with $\matup' =
(\matlow')^{\trp}$, hence  $ \proj_{\MM^{\trp}}  = \proj_{\MM}$.
Now, by Fact~\ref{fac:pinvproduct},
\begin{equation}
  \label{eq:pinvproduct}
\MM^{\pinv} = \proj_{\MM} (\matlow')^{-\trp}\matdiag^{\pinv}\matlow^{-1}
\proj_{\MM}.
\end{equation}
The kernel of $\MM$ is precisely the span of the
indicator vectors of each connected component of the associated graph,
and hence $\proj_{\MM}$ can be applied by, for each connected
component, subtracting the mean value for that component from every
entry of the component.

%% file: appendixUnbiased.tex
\section{Unbiased approximate Cholesky factorization}
\label{sec:unbiased}

In this section, for completeness, we prove
Claim~\ref{clm:unbiasedchol}, which was shown in
\cite{KS16}\footnote{See also Lemma 4.1 in lecture notes at \url{http://kyng.inf.ethz.ch/courses/AGAO20/lectures/lecture9_apxgauss.pdf}.}.
\begin{claim}
  When $\textsc{ApproximateCholesky}$ is instantiated with an unbiased clique
  sampling routine \textsc{Clique\-Sample} so that
  \[
    \expct{}{\textsc{CliqueSample}(v,\SS)} = \cliq{\SS}{v}.
  \]
  Then letting $\matlow = \textsc{ApproximateCholesky}(\LL)$, we have
  \[
    \expct{}{\matlow \matlow^{\trp}} = \LL
    .
  \]
  In other words, the approximate Cholesky factorization equals the
  original matrix in expectation, which is in turn also equal, viewed
  as a linear operator, to any
  exact Cholesky factorization.
\end{claim}

\begin{proof}
Let $\SS_i$ denote $\SS$ in
Algorithm~\ref{alg:GeneralApproximateCholeskyFactorizationColwise}
after $i$ eliminations, and note that $\SS_0 = \LL$.
Let
\[
  \LL_i = \SS_i + \sum_{j=1}^i \ll_j\ll_j^\trp.
\]

Suppose the $i$th vertex to be eliminated is vertex $v$.
Conditional on the samples before the $i$th
elimination, we have
  \begin{align*}
    \expct{}{\LL_{i}}
     & = \expct{}{ \SS_i + \sum_{j=1}^i \ll_j\ll_j^\trp }. \\
     & = \expct{}{ \SS_{i-1}-\vstar{\SS_{i-1}}{v} + \textsc{CliqueSample}(v, \SS_{i-1})
       + \sum_{j=1}^i \ll_j\ll_j^\trp }\\
         & = \SS_{i-1}-\vstar{\SS_{i-1}}{v} + \expct{}{ \textsc{CliqueSample}(v, \SS_{i-1})}
       + \sum_{j=1}^i \ll_j\ll_j^\trp \\
     & = \SS_{i-1}-\vstar{\SS_{i-1}}{v} + \cliq{\SS_{i-1}}{v}
       + \sum_{j=1}^i \ll_j\ll_j^\trp\\
         & = \SS_{i-1} 
           + \sum_{j=1}^{i-1} \ll_j\ll_j^\trp\\
    & = \LL_{i-1}.
  \end{align*}
We can now chain together expectations to conclude that over all the
randomness of the algorithm
\[
  \expct{}{\matlow\matlow^{\trp}} =  \expct{}{\LL_{n}}
  =  \LL.
  \]
\end{proof}

%% file: appendixEquivalence.tex
\section{Equivalence of the output representations}
\label{sec:equivalence}

In this section, we prove Claim~\ref{clm:equivalence}, about the
equivalence of edgewise and standard Cholesky factorization.
We restate the lemma here for convenience:

\begin{claim*}
  \quad \\
Suppose we run both $\textsc{ApproximateCholesky}$ and 
$\textsc{ApproximateEdgewiseCholesky}$
\begin{itemize}
\item using the same elimination ordering,
\item and using $\textsc{CliqueTreeSample}$ as the clique
  sampling routine, and using the same outputs
  of \textsc{Clique\-Tree\-Sample}, i.e. the same outcomes of the random
  samples.
\end{itemize}
  Then the output factorizations $\matlow$ from
  $\textsc{ApproximateCholesky}$ and $\PPhi,
  \setof{\matlow^{(v)}_i}_{v,i}$ from
  \textsc{Approximate\-Edgewise\-Cholesky} will satisfy
  \[
    \matlow \matlow^{\trp} =
    \left(
   \Pi_{i = 1}^n \Pi_{j \sim i }
    \matlow_j^{(i)}
\right)
\PPhi
  \left(
   \Pi_{i = 1}^n \Pi_{j \sim i }
    \matlow_j^{(i)}
  \right)^{\trp}
  .
\]
In other words, the only difference between the two algorithms is in
the formatting of the output.
\end{claim*}

\begin{proof}
We can write a Laplacian $\LL$ with terms corresponding to edges of
the first vertex explicitly separated out as:
\begin{equation}
  \label{eq:lapuniformlayoutagain}
  \LL = 
\left(
\begin{array}{ccc}
d & -\aa^\trp \\
-\aa& \diag(\aa) + \LL_{-1}
\end{array} \right)
\end{equation}
Here $\LL_{-1}$ is the graph Laplacian corresponding to the induced
subgraph on the vertex set with vertex 1 removed.

From Section~\ref{sec:apxChol}, we know that it is possible to write
\begin{equation}
  \label{eq:singlevertexchol}
  \LL
  =
\begin{pmatrix}
0 & \veczero^\trp \\
  \veczero &  \SS
\end{pmatrix}
+
  \frac{1}{d}
  \begin{pmatrix}
    d \\
    -\aa
  \end{pmatrix}
    \begin{pmatrix}
    d \\
    -\aa
  \end{pmatrix}^\trp
\end{equation}
where $\SS = \schurto{\LL}{[n]\setminus\setof{1}}$.

If vertex $1$ has $k$ neighbors, then using the elimination procedure described in
Section~\ref{sec:rowform} to eliminate the first row and
column of the matrix will result in a partial factorization.
Thus,
\begin{equation}
  \label{eq:edgeelimfactoragain}
  \LL
  =
  \matlow_1 \matlow_2 \ldots \matlow_k
\left(
\begin{array}{ccc}
\phi & \veczero^\trp \\
  \veczero &  \SS
\end{array}
\right)
\matlow_k^{\trp}\ldots  \matlow_2^{\trp} \matlow_1^{\trp}
\end{equation}
where each
$\matlow_i^{\trp}$ is a lower-triangular matrix corresponding to a
single row-operation.

We can prove the following:

\begin{equation}
  \label{eq:schurunchanged}
\left(
\begin{array}{ccc}
 0  & \veczero^\trp \\
  \veczero &  \SS
\end{array}
\right)
  =
  \matlow_1 \matlow_2 \cdots \matlow_k
\left(
\begin{array}{ccc}
0 & \veczero^\trp \\
  \veczero &  \SS
\end{array}
\right)
\matlow_k^{\trp}\cdots  \matlow_2^{\trp} \matlow_1^{\trp}
\end{equation}
and
\begin{equation}
\label{eq:rank1partsagree}
  \frac{1}{d}
  \begin{pmatrix}
    d \\
    -\aa
  \end{pmatrix}
    \begin{pmatrix}
    d \\
    -\aa
  \end{pmatrix}^\trp
  =
  \matlow_1 \matlow_2 \cdots \matlow_k
\left(
\begin{array}{ccc}
\phi & \veczero^\trp \\
  \veczero &  \matzero
\end{array}
\right)
\matlow_k^{\trp}\cdots  \matlow_2^{\trp} \matlow_1^{\trp}
\end{equation}

To establish Equation~\eqref{eq:schurunchanged}, we observe
that each row operation factor takes the form
\[
  \matlow_i = \begin{pmatrix}
 \bb_i &  \veczero^\trp   \\
 \vrule & \II_{(n-1) \times (n-1)}
\end{pmatrix}
\]
for some vector  $\bb_i \in R^n$.
But, given any vector $\bb_i \in R^n$ and %
any matrix $\CC \in R^{(n-1) \times d}$, with any
number of columns $d$, we have
\[
\begin{pmatrix}
 \bb &  \veczero^\trp   \\
 \vrule & \II_{(n-1) \times (n-1)}
\end{pmatrix}
\begin{pmatrix}
 \veczero^\trp \\
 \CC
\end{pmatrix}
=
\begin{pmatrix}
 \veczero^\trp \\
 \CC
\end{pmatrix}
.
\]
Repeatedly applying this observation the right hand side of
Equation~\eqref{eq:schurunchanged} lets us conclude that it equals the
left hand side, proving the equation holds.

Now, by equating the two right hand sides of
Equation~\eqref{eq:singlevertexchol} and
Equation~\eqref{eq:edgeelimfactor}, and simplifying using
Equation~\eqref{eq:schurunchanged},
we conclude that Equation~\eqref{eq:rank1partsagree} holds.

The proof of Equation~\eqref{eq:rank1partsagree}
only relies on the row operation matrices having the form $\matlow_i = \begin{pmatrix}
 \bb_i &  \veczero^\trp   \\
 \vrule & \II_{(n-1) \times (n-1)}
\end{pmatrix}$ for some $\bb_i$ and does not require $\SS$ to be a
Schur complement.

Consquently, it holds for any $\SStil$ matrix that 
\begin{equation}
  \label{eq:apxschurunchanged}
\left(
\begin{array}{ccc}
 0  & \veczero^\trp \\
  \veczero &  \SStil
\end{array}
\right)
  =
  \matlow_1 \matlow_2 \cdots \matlow_k
\left(
\begin{array}{ccc}
0 & \veczero^\trp \\
  \veczero &  \SStil
\end{array}
\right)
\matlow_k^{\trp}\cdots  \matlow_2^{\trp} \matlow_1^{\trp}
\end{equation}
and hence 
\begin{equation}
\label{eq:apxformsequivalence}
  \begin{pmatrix}
0 & \veczero^\trp \\
  \veczero &  \SStil
\end{pmatrix}
+
  \frac{1}{d}
  \begin{pmatrix}
    d \\
    -\aa
  \end{pmatrix}
    \begin{pmatrix}
    d \\
    -\aa
  \end{pmatrix}^\trp
  =
  \matlow_1 \matlow_2 \cdots \matlow_k
\left(
\begin{array}{ccc}
\phi & \veczero^\trp \\
  \veczero &  \SStil
\end{array}
\right)
\matlow_k^{\trp}\cdots  \matlow_2^{\trp} \matlow_1^{\trp}
.
\end{equation}

Thus, by taking $\SStil = \LL_{-1} + \textsc{CliqueTreeSample}(v,\SS)$, 
we can conclude that the partial factorizations computed in one step of $\textsc{ApproximateCholesky}(\LL)$ and
$\textsc{EdgewiseApproximateCholesky}(\LL)$ are identical.
We can apply this equivalence repeatedly to conclude that the two
approximate factorizations returned by
$\textsc{ApproximateCholesky}(\LL)$ and
$\textsc{EdgewiseApproximateCholesky}(\LL)$ are identical, assuming we
\begin{itemize}
\item use the same elimination ordering for both,
\item and use $\textsc{CliqueTreeSample}$ as the clique
  sampling routine, and using the same outputs
  of \textsc{Clique\-Tree\-Sample}, i.e. the same outcomes of the random
  samples.
\end{itemize}
Thus the approximate factorizations returned by $\textsc{ApproximateCholesky}(\LL)$ and
\textsc{Edgewise\-Approximate\-Cho\-lesky}(\LL) are identical when
regarded as operators.

\end{proof}